\documentclass[11pt]{elsarticle}

\usepackage[a4paper,margin=.7in]{geometry}
\usepackage[T1]{fontenc}
\usepackage{amsmath,amssymb,amsthm,mathtools}
\usepackage{mathrsfs}
\usepackage{enumitem}
\usepackage{dsfont}
\usepackage[hidelinks]{hyperref}
\usepackage[ nameinlink, noabbrev, capitalise]{cleveref}
\allowdisplaybreaks[2]

\newtheorem{theorem}{Theorem}[section]
\newtheorem{lemma}[theorem]{Lemma}

\newtheorem{remark}[theorem]{Remark}
\newtheorem{definition}[theorem]{Definition}
\newtheorem{proposition}[theorem]{Proposition}

\newtheorem{example}[theorem]{Example}
\newcommand{\R}{\mathbb R}
\newcommand{\E}{\mathbb E}
\newcommand{\N}{\mathbb N}
\renewcommand{\P}{\mathbb P}
\newcommand{\U}{\mathcal U}
\newcommand{\LP}{\mathcal P}
\newcommand{\TN}{\mathcal T_N}
\newcommand{\HS}{\mathcal L_2}
\newcommand{\PN}{\mathcal N}
\newcommand{\tPN}{\widetilde{\mathcal N}}
\newcommand{\dd}{\mathrm d}
\newcommand{\diver}{\operatorname{div}}

\usepackage{orcidlink}
\crefname{equation}{equation}{equations}
\crefname{remark}{remark}{remarks}
\Crefname{remark}{Remark}{Remarks}
\newcounter{namedhypothesis}

\newcommand{\hypothesislabel}[2]{
  \begingroup
  \renewcommand{\thenamedhypothesis}{#1}
  \refstepcounter{namedhypothesis}
  \label[hyp#1]{#2}
  \endgroup
}
\crefformat{hypT}{#2\textup{(T)}#3}
\crefformat{hypL}{#2\textup{(L)}#3}
\crefformat{hypG}{#2\textup{(G)}#3}
\crefformat{hypT1}{#2\textup{\ensuremath{(\mathrm{T}_{\mathrm{loc}})}}#3}

\Crefformat{hypT}{#2\textup{(T)}#3}
\Crefformat{hypL}{#2\textup{(L)}#3}
\Crefformat{hypG}{#2\textup{(G)}#3}
\Crefformat{hypT1}{#2\textup{\ensuremath{(\mathrm{T}_{\mathrm{loc}})}}#3}
\begin{document}
\begin{frontmatter}

\title{Global Critical $L^3$ Well-Posedness for the Stochastic Tamed Navier--Stokes Equations on $\mathbb R^3$  with Wiener and Jump Noise }
\author{Bikram Podder \orcidlink{0009-0008-8765-3018}} \ead{bikrampoddar2@gmail.com} \author{Surendra Kumar* \orcidlink{0000-0001-6802-5894}} \ead{surendraiitr8@gmail.com}

\affiliation{organization={Department of Mathematics}, addressline={University of Delhi}, city={ Delhi}, postcode={110007}, country={India}}
\begin{abstract}

We study the global well-posedness for the three dimensional stochastic tamed Navier--Stokes equations on $\mathbb R^3$ with multiplicative Wiener and compensated Poisson noise.
When noise coefficients are globally Lipschitz on $L^3_\sigma(\mathbb R^3)$ and satisfy spatial $H^1_\sigma(\mathbb R^3)$ regularity assumptions, the system admits a  pathwise unique global,  probabilistically strong, critical mild solution for any  initial datum $u_0 \in L^0(\Omega, \mathcal F_0; L^3_\sigma(\mathbb R^3)) $.
Neither smallness nor moment or finite kinetic energy assumptions on $u_0$ are required.
We also establish that the  global solution depends continuously in probability on the initial datum in uniform-in-time $L^3$ and space-time $L^5$ norms on every finite time interval.
\end{abstract}

\begin{keyword}
Stochastic PDEs  \sep 
Tamed stochastic Navier--Stokes Equation  \sep Global strong solution

\MSC[2020] 60H15 \sep 35R60 \sep 35Q30 \sep 76D05

\end{keyword}

\end{frontmatter}
\tableofcontents

\section{Introduction}

The nonlinear damping in the three-dimensional tamed Navier--Stokes equations (NSE) provides an additional source of dissipation that is effective in energy and Sobolev estimates.
Its use for initial velocities belonging only to the scaling-critical space $L^3_\sigma(\R^3)$ [see \cref{Sec2.1}] raises a different question.
On the whole space, such velocities need not have finite kinetic energy, so the usual energy argument cannot be applied to the full solution from the initial time.
The problem we consider is whether the dissipation supplied by taming can nevertheless yield global well-posedness in the critical mild class when the forcing contains both a Wiener term and a compensated Poisson term, each with a velocity-dependent coefficient.

More precisely, we study
\begin{equation}
\label{eq:tnse}
\begin{cases}
\begin{aligned}
\dd u(t) ={}&\bigl[ \nu\Delta u(t) -\LP\diver(u(t)\otimes u(t)) -\LP\TN(u(t)) \bigr]\dd t \\
&+\Sigma(t,u(t-))\dd W(t) +\int_Z\mathcal G(t,u(t-),z)\, \tPN(\dd t,\dd z),
\end{aligned}
\\
\diver u(t)=0,\qquad (t,x)\in(0,\infty)\times\R^3, \\
u(0)=u_0,
\end{cases}
\end{equation}
where $\nu>0$ is the viscosity, $\LP$ denotes the Leray projection, and $\TN(u)=g_N(|u|^2)u$ is the taming term.
The initial velocity is a strongly $\mathcal F_0$-measurable $L^3_\sigma(\R^3)$-valued random variable.
The cylindrical Wiener process and Poisson random measure are prescribed on a fixed stochastic basis, and the noise coefficients take values in solenoidal spaces [see \cref{Sec2.2}].
The coefficients are evaluated at \(u(t-)\), as required for predictable stochastic integration . 
The natural solution space for \eqref{eq:tnse} is $$ D([0,T];L^3_\sigma(\R^3)) \cap L^5(0,T;L^5_\sigma(\R^3)).
$$ Its two defining norms are invariant under the classical Navier--Stokes rescaling $ u_\lambda(t,x)=\lambda u(\lambda^2t,\lambda x), $ with the corresponding rescaling of the time interval.
We use the term \emph{critical} in this functional sense.
Neither the fixed taming profile nor the stochastic coefficients are assumed to be invariant under this transformation.

The answer we seek in this article whether the coercive taming term in \cref{eq:tnse} suffices to extend the small data critical $L^3$ theory of Kukavica and Xu \cite{MR4998100} for Wiener driven stochastic NSE on $\mathbb T^3$ to arbitrary initial data on whole space  $\R^3$ for both Wiener and Jump noise.

The endpoint $L^3$ formulation requires particular care with the nonlinear drift.
The cubic growth of the taming term guarantees only $\TN(u)\in L^1(\R^3)$ for $u\in L^3(\R^3)$, whereas the Leray projection is not bounded on $L^1(\R^3)$.
Although the unprojected taming term is dissipative, this property does not directly determine the sign of its contribution to the projected $L^3$ duality pairing.
The space--time $L^5$ regularity offers a way to formulate the nonlinear terms without relying on such an endpoint projection estimate: $$ u\otimes u\in L^{5/2}((0,T)\times\R^3), \qquad \TN(u)\in L^{5/3}((0,T)\times\R^3).
$$ The corresponding heat convolutions can then be studied in the critical mild topology.

Local solvability in this topology, however, does not itself provide the energy control needed for global continuation.
A natural starting point is to separate the heat evolution of the initial velocity, $$ h(t)=e^{\nu t\Delta}u_0, \qquad R(t)=u(t)-h(t).
$$ The linear heat flow carries the possible infinite-energy part of the initial state, while $R$ starts from zero initial data.
Whether $R$ possesses finite-energy regularity is a separate analytical question.
It does not follow merely from this decomposition.
Moreover, any global Sobolev evolution obtained via dissipative estimates must still be identified with the original critical mild solution where the driving Poisson noise  has to respect both the stopping-time formulation and the jumps of the trajectories.

\subsection{Related literature} 
Kato's critical $L^3$ theory \cite{MR760047} for mild solutions of the incompressible NSE on $\R^3$ gives local existence and uniqueness for divergence-free initial velocities in $L^3(\R^3)$ and global solvability for sufficiently small data.
Koch and Tataru \cite{MR1808843} extended the small-data global well-posedness theory for the classical NSE to the critical space $\mathrm{BMO}^{-1}(\R^3)$.
At the endpoint of the classical Ladyzhenskaya--Prodi--Serrin regularity criterion, the theorem of Escauriaza, Seregin, and \v{S}ver\'ak \cite{MR1992563} implies that any Leray--Hopf weak solution of the three-dimensional incompressible Navier--Stokes Cauchy problem on $\R^3$ satisfying $ u\in L^\infty(0,T;L^3(\R^3)) $ belongs to $L^5((0,T)\times\R^3)$ and is therefore smooth and unique on $\R^3\times(0,T)$.
Their result identifies the boundedness of the scaling-critical $L^3$ norm as an endpoint regularity criterion for the classical three-dimensional NSE.
This highlights the importance of controlling the critical $L^3$ norm, since the usual energy estimates alone are insufficient to rule out finite-time singularities.

The stochastic Navier--Stokes theory developed initially in Hilbert-space settings, where energy estimates and compactness provide a natural starting point.
Bensoussan and Temam \cite{MR348841}, established existence and energy estimates for stochastic Navier--Stokes-type evolution equations with random initial data and forcing, employing measurable-selection arguments.
For the two-dimensional stochastic incompressible NSE with additive Wiener forcing, Menaldi and Sritharan \cite{MR1922695} constructed probabilistically strong solutions on both bounded and unbounded domains and proved their pathwise uniqueness, using local monotonicity and weighted energy estimates on a prescribed probability space.
Flandoli \cite{MR2459085} provided a systematic treatment of three-dimensional stochastic Navier--Stokes dynamics, including martingale solutions, stationary solutions, Markov selections, and related probabilistic properties, primarily for Gaussian forcing in the periodic setting.
While these works laid the probabilistic foundations for stochastics fluid equations, they work with different notions of solutions and regularity than the critical setting studied here 

Discontinuous stochastic forcing introduces additional difficulties, since the solution must retain its jumps while satisfying the appropriate energy or evolution estimates.
For the incompressible NSE subjected to Gaussian and L\'evy noise, Sakthivel and Sritharan \cite{MR3085233} established the existence of martingale solutions on bounded and unbounded domains in dimensions two and three, and they also obtained pathwise uniqueness in two dimensions and conditional uniqueness in three dimensions under additional integrability assumptions.
Motyl \cite{MR3034603} studied the existence of martingale solutions for the stochastic incompressible NSE on possibly unbounded domains in dimensions two and three, with a cylindrical Wiener component and a jump component represented by a compensated Poisson random measure.
For the three-dimensional stochastic nonhomogeneous incompressible NSE, Chen, Wang, and Wang \cite{MR3912800} proved the existence of martingale solutions on bounded domains with variable density, under Brownian forcing together with compensated small-jump and uncompensated large-jump Poisson perturbations.
Brze\'zniak, et al.
\cite{MR5073010} dealt with the existence of global martingale solutions to the incompressible stochastic NSE with purely discontinuous L\'evy forcing on $\R^2$ and $\R^3$. These results illustrate the several range of stochastic solution concepts available for incompressible fluid equations.
They do not, by themselves, give pathwise global continuation in a prescribed critical mild class.

For the classical, incompressible NSE on $\R^n$, $n=2,3$, Fernando, R\"udiger, and Sritharan \cite{MR3411979} established existence and uniqueness of local mild solutions in $L^p_\sigma(\R^n)$, $p>n$, under additive compensated Poisson forcing, using heat-semigroup estimates and Banach-space stochastic convolution inequalities.
Mohan and Sritharan \cite{MR3669657} investigated local mild solvability of the stochastic NSE on $\R^m$, $m\geq2$, with critical initial data in $L^m_\sigma(\R^m)$, under simultaneous additive Wiener and compensated Poisson perturbations. The significance of this result for the present problem is that the critical initial-data class and jump forcing can be treated together within a mild formulation.
Its local character, however, leaves global continuation as a separate question.

The interaction between low regularity and jump forcing has also been studied in other geometries and fluid models.
In the two-dimensional bounded-domain setting, Zhu, Brze\'zniak, and Liu \cite{MR3991628} constructed a unique global mild solution of the stochastic incompressible NSE with spatially irregular additive compensated Poisson forcing. For $0<\alpha<1/4$, their theorem allows initial data in $H^{-2\alpha,4}(D)$ and gives $u\in L^4(0,T;L^4_\sigma(D))$ almost surely.
Related local well-posedness theories for stochastic incompressible fluid equations with discontinuous forcing have also been developed beyond the classical Navier--Stokes framework.
For the stochastic incompressible Euler equations on $\R^d$, $d=2,3$, driven simultaneously by multiplicative Wiener and compensated Poisson noise, Kukavica and Xu \cite{MR3950967} showed the local pathwise existence and uniqueness of solutions with c\`adl\`ag trajectories in $W^{s,p}(\R^d)$ for $p>2$ and $s>1+d/p$.
For the three-dimensional stochastic incompressible NSE with hereditary viscosity, Mohan and Sritharan \cite{MR3942498} constructed a unique local pathwise strong solution whose paths are c\`adl\`ag in $H^s_\sigma(\R^3)$, $s>7/2$, under multiplicative Wiener and compensated Poisson forcing.
They further obtained global existence with probability arbitrarily close to one for sufficiently small initial data under additional structural conditions on the noise coefficients.
Although these equations differ from \eqref{eq:tnse}, the results make clear that the regularity of the initial state, the type of noise, and the desired notion of global existence must be distinguished.

For multiplicative Wiener forcing, maximal pathwise solution theories have been obtained in several Sobolev and Lebesgue settings.
For the stochastic incompressible NSE driven by Wiener noise, Glatt-Holtz and Ziane \cite{MR2502705} constructed unique maximal pathwise $H^1$-strong solutions on smooth bounded domains in dimensions two and three, established a blow-up alternative, and proved global existence almost surely in two dimensions under suitable integrability assumptions.
For the three-dimensional whole-space equation, Kim \cite{MR2865433} proved local existence and pathwise uniqueness of strong solutions with initial data in $H^{1/2+\alpha}_\sigma(\R^3)$, $0<\alpha<1/2$, allowing multiplicative Wiener forcing and obtained global existence for sufficiently small initial data under additional stabilizing conditions on the noise coefficients.
Recently, Kukavica, Xu, and Ziane \cite{MR4385406} developed local existence and pathwise uniqueness of maximal strong $L^p$-solutions, $p>5$, for the three-dimensional stochastic incompressible NSE on $\mathbb T^3$ with multiplicative Wiener noise.
They also proved global existence when the initial $L^p$ moment and the noise amplitude are sufficiently small.
For the three-dimensional stochastic incompressible NSE driven by multiplicative Wiener noise, Kukavica and Xu \cite{MR4552356} studied the local existence and pathwise uniqueness of strong solutions on $\mathbb T^3$ with initial data in $L^p_\sigma(\mathbb T^3)$, $p>3$, together with localized $L^p$-energy estimates.
Kukavica, Wang, and Xu \cite{MR4908978} subsequently extended this local $L^p$ theory to the whole space $\R^3$.
These works are closely related to the use of $L^p$ energy estimates and stopping times in stochastic fluid equations, while the critical endpoint $p=3$ requires a different treatment of the nonlinear drift.

The stochastic critical-space theory has advanced beyond the classical $L^p$ setting.
Critical Besov-space formulations for the stochastic incompressible NSE with multiplicative Wiener forcing were developed by Du and Zhang \cite{MR4002154}, who constructed a unique local strong solution on $\R^d$ for initial velocities in $\dot B^{d/p-1}_{p,r}(\R^d)$, $2\le p,r<\infty$, under their hypotheses on the noise coefficients. They obtained global existence when both the initial data and the stochastic forcing satisfy their smallness conditions.
For stochastic NSE with transport-type Wiener noise on $\mathbb T^d$, Agresti and Veraar \cite{MR4703457} established local well-posedness in critical Besov spaces $B^{d/q-1}_{q,p}(\mathbb T^d)$, together with instantaneous spatial and temporal regularization and stochastic Serrin-type blow-up criteria.
In particular, critical $L^3$ initial velocities are already admissible in an established stochastic well-posedness framework on $\mathbb T^3$.
For the stochastic incompressible NSE on $\R^d$ with transport-type Wiener noise, Agresti and Veraar \cite{MR4952170} proved that the problem is locally well posed in scaling-critical Besov spaces and that its solutions regularize. Their conditions on the parameters and coefficients allow initial velocities in $L^d_\sigma(\R^d)$.
These results show that the local, critical $L^3$ initial data theories are already well established. The main obstacle addressed here is to extend this local theory to global existence within this critical class.  

At the critical Sobolev regularity, Ayd{\i}n, Kukavica, and Xu \cite{AKXHhalf2025} studied the stochastic incompressible NSE on $\mathbb T^3$ with multiplicative cylindrical Wiener noise.
For sufficiently small initial data in $H^{1/2}_\sigma(\mathbb T^3)$ they established global existence together with critical Sobolev energy estimates, and obtained local small-data well-posedness without requiring the noise amplitude to be small.
Ayd{\i}n and Xu \cite{AXTransport2025} investigated the critical $H^{1/2}$ formulation with linear multiplicative transport noise.
Under smallness assumptions on the transport field and the initial data, they established pathwise unique strong solutions and global existence for $\mathbb T^3$ and $\R^3$.
Kukavica and Xu \cite{MR4998100} established almost global existence and pathwise uniqueness for the stochastic incompressible NSE with sufficiently small $L^3$ initial data and multiplicative cylindrical Wiener noise on $\mathbb T^3$. Under suitable smallness assumptions on the noise coefficient, their theorem  shows that the probability of global existence can be made arbitrarily close to one.
Ayd{\i}n, Kukavica, and Xu \cite{MR5107357} subsequently proved local solvability for arbitrary-size initial data in $L^\infty(\Omega;L^3_\sigma(\mathbb T^3))$ under suitable Lipschitz assumptions on the multiplicative Wiener coefficient.
Taken together, these results distinguish the local arbitrary-data problem from the global small-data and high-probability problems for the untamed equation.

Nonlinear damping changes the global analysis.
For the deterministic tamed three-dimensional NSE on $\R^3$, R\"ockner and Zhang \cite{MR2590154} proved global existence, uniqueness, and Sobolev regularity, and further proved subsequential convergence towards suitable weak solutions of the classical NSE as the taming parameter tends to infinity.
In their study of the stochastic tamed three-dimensional NSE driven by Wiener noise, R\"ockner and Zhang \cite{MR2520127} showed global existence and uniqueness of strong solutions with $H^1$ initial data on $\R^3$ and $\mathbb T^3$, allowing transport-type and state-dependent multiplicative noise.
For bounded domains in $\R^3$ with homogeneous Dirichlet boundary conditions, R\"ockner and T.
Zhang \cite{MR2852224} analyzed global strong solvability and uniqueness for the three-dimensional tamed NSE with $H^1$ initial data and multiplicative Wiener forcing.

Other noise structures have also been incorporated into the tamed equations.
Brze\'zniak and Dhariwal \cite{MR4085355} subsequently established existence and pathwise uniqueness for the stochastic tamed NSE on $\R^3$ with transport-type Wiener noise, and proved the existence of an invariant measure for an additionally linearly damped, time-homogeneous version of the equation.
For multiplicative jump forcing, Dong and Zhang \cite{MR4127295} constructed a unique global strong solution of the three-dimensional periodic tamed NSE for $H^1$ initial data. The jump term is integrated against a compensated Poisson random measure, and the solution has c\`adl\`ag paths in $H^1$.
More recently, Lu, Yang, and Li \cite{LuYangLi2025} investigated Wiener-driven stochastic tamed three-dimensional NSE on $\R^3$ and $\mathbb T^3$ under locally weak monotonicity and growth assumptions on the drift and diffusion coefficients. Under these hypotheses they obtained a unique global strong solution in the $H^1$ framework.

Related deterministic work concerns the three-dimensional incompressible NSE with nonlinear power-law damping of the form $\beta |u|^{r-1}u$, which exhibits cubic growth at the exponent $r=3$.
For the Cauchy problem on $\R^3$, Cai and Jiu \cite{MR2401535} established the existence of global weak solutions for divergence-free $L^2$ initial datawhen $r\geq1$.  For initial data in $H^1\cap L^{r+1}$, they obtained global strong solutions for $r\geq7/2$ and proved uniqueness for $7/2\leq r\leq5$.
Zhang, Wu, and Lu \cite{MR2754840} subsequently improved the strong-solution theory for the same three-dimensional NSE with nonlinear damping, extending global strong existence to every $r>3$ and uniqueness to $3<r\leq 5$.
For stochastic three-dimensional incompressible NSE with nonlinear damping on smooth bounded domains, Gao and Liu \cite{MR3996790} established global existence and uniqueness of strong solutions driven by multiplicative jump noise under suitable growth and coercivity assumptions, with applications to stochastic tamed NSE.
For the periodic three-dimensional NSE with power-law damping, Hajduk and Robinson \cite{MR3705676} studied global strong solvability for $r>3$. At the cubic exponent $r=3$, their strong-solvability result assumes $4\mu\beta\geq1$.
They further established the energy equality for every weak solution in the critical case $r=3$, without requiring this additional viscosity--damping restriction.
The common feature is that nonlinear absorption supplies coercivity, though its use depends on both the solution class and the precise diffusion and damping assumptions.

Both the deterministic and stochastic Brinkman--Forchheimer literature provides especially relevant comparisons because the cubic damping has a similar algebraic structure.
Mohan \cite{MR4314117} considered the deterministic convective Brinkman--Forchheimer equations (BFE) on $\R^d$ and analyzed local mild solvability and uniqueness in $L^p_\sigma(\R^d)$ for $ p>\max\left\{d,\frac{d(r-1)}{2}\right\}.
$ The author also proved local pathwise mild solvability for the corresponding stochastic equations driven simultaneously by additive Wiener and compensated Poisson noise, together with global solvability under additional damping, coefficient, and noise-integrability conditions.
In particular, the local solvability restriction reduces to $p>3$ when $d=r=3$.

For finite-energy initial data on bounded or periodic domains, Mohan \cite{MR4439993} constructed global probabilistically strong solutions of the same equations with multiplicative pure-jump forcing and proved pathwise uniqueness. In dimensions three and four, the result covers $r>3$ and also $r=3$ when $2\beta\mu\geq1$.

For multiplicative Gaussian forcing, Kinra and Mohan \cite{KinraMohanUnbounded2025} developed a global well-posedness theory for stochastic convective BFE on general unbounded domains, including the cubic case under $2\beta\mu\geq1$.
Their framework admits finite-energy initial data with additional moment integrability and yields pathwise unique probabilistically strong solutions with continuous $L^2$-valued trajectories.
Wang \cite{MR5088480} proved a well-posedness theorem on unbounded domains for a broad class of SPDEs whose drift is polynomial and whose nonlinear Wiener coefficient has critical growth.
As an application, the author proved existence and uniqueness for stochastic BFE on $\R^3$ with $L^2$ initial data and absorption exponent $r\geq3$, subject to structural monotonicity and coercivity conditions that balance the nonlinear drift and noise.

The abstract theory of stochastic evolution equations with locally monotone coefficients is also relevant to the analytical framework.
Liu and R\"ockner \cite{MR2719279} extended the classical variational approach to coercive stochastic evolution equations with locally monotone drift and diffusion coefficients, establishing existence and uniqueness of strong variational solutions under multiplicative Wiener forcing.
Their framework applies, in particular, to the two-dimensional stochastic incompressible NSE.

In their subsequent work \cite{MR2990049}, Liu and R\"ockner established local and global well-posedness results for nonlinear evolution equations under generalized coercivity conditions, together with an abstract stochastic theory for additive Wiener forcing.
They also applied the deterministic framework to the three-dimensional tamed NSE, obtaining global existence and uniqueness of strong solutions with $H^1$ initial data.

For stochastic evolution equations driven by L\'evy noise, Brze\'zniak, Liu, and Zhu \cite{MR3158475} established existence and pathwise uniqueness of strong variational solutions with c\`adl\`ag trajectories under coercivity, local monotonicity, and suitable growth and integrability assumptions on a Gelfand triple.
Their framework accommodates multiplicative Wiener noise together with compensated and uncompensated Poisson perturbations and applies, in particular, to two-dimensional stochastic incompressible NSE and other stochastic hydrodynamical systems.

Agresti and Veraar \cite{MR4459102,MR4437443} used temporal weights and stochastic maximal $L^p$-regularity to prove local well-posedness for nonlinear parabolic stochastic evolution equations in critical interpolation spaces. They also derived stochastic criteria for blow-up and proved  the instantaneous spatial and temporal regularization of that solutions. 
In their subsequent critical variational framework \cite{MR4716343}, they proved local and global existence and uniqueness under suitable coercivity and local Lipschitz conditions.
Their theory gives global well-posedness for the three-dimensional stochastic tamed NSE on $\R^3$ with $H^1$ initial data, allowing multiplicative Wiener forcing with gradient-dependent transport terms.
Bechtel, Germ, and Veraar \cite{MR5120871} extended the critical variational framework to stochastic evolution equations with both a Wiener term and a jump integral against a compensated Poisson random measure, establishing local well-posedness, blow-up criteria, and global existence and uniqueness under appropriate coercivity assumptions.
They explicitly identified the tamed NSE among the applications of the earlier framework that admit an extension to L\'evy forcing.
Their theory does not require compactness of the embedding between the underlying Hilbert spaces and therefore also applies to suitable equations on unbounded domains.

These variational developments are important precedents.
It shows that the global solvability of stochastic tamed equations is not a new question in itself.
The distinctions are mainly the initial-data topology and the form of the solution.
A global variational solution constructed from an $H^1$ or finite-energy state does not, without a further assumptions, give a global critical mild solution for an arbitrary initial velocity in $L^3_\sigma(\R^3)$ that does not belong to $L^2_\sigma$.

The linear stochastic analysis underlying these approaches includes the stochastic maximal $L^p$-regularity theorem of van Neerven, Veraar, and Weis \cite{MR2952092} for Wiener-driven stochastic convolutions associated with analytic semigroups whose generators admit a bounded $H^\infty$-functional calculus.
For compensated Poisson forcing, Brze\'zniak and Hausenblas \cite{MR2529441} established maximal regularity of stochastic convolutions in martingale-type Banach spaces, using analytic semigroups and real interpolation spaces.
Hausenblas \cite{MR2832576} established maximal inequalities for Banach-valued stochastic integrals with respect to compensated Poisson random measures, including higher-moment estimates under suitable integrability assumptions.
Dirksen \cite{MR3265175} subsequently obtained sharp two-sided It\^o isomorphisms for $L^q$-valued compensated Poisson stochastic integrals, identifying the relevant square-function and higher-order jump-integrability contributions.
These estimates provide relevant stochastic tools for controlling the quadratic and cubic jump terms arising in $L^3$-moment arguments.
Zhu, Brze\'zniak, and Hausenblas \cite{MR3634281} obtained maximal inequalities for convolutions against compensated Poisson random measures and showed that these convolutions admit c\`adl\`ag modifications in suitable Banach spaces, assuming appropriate smoothness of an equivalent Banach norm and contractivity of the underlying strongly continuous semigroup with respect to that norm.
Gy\"ongy and Wu \cite{MR4165650} established It\^o formulas for the $L^p$-norms, $p\ge 2$, of stochastic processes driven by Wiener and compensated Poisson noise, including equations with divergence-form drift.
Their results provide a basis for $L^p$-energy identities that explicitly retain the nonlinear Poisson jump contributions.

Deterministic theories of the incompressible NSE provide closely related precedents for combining the evolution of rough initial data with a finite-energy correction.
Calder\'on \cite{MR1018571} introduced an initial-data decomposition into a sufficiently small critical $L^3$ component and a finite-energy $L^2$ component, thereby constructing global weak solutions on $\R^3$ for arbitrary divergence-free initial data in $L^p(\R^3)$, $2<p<3$.
Seregin and \v{S}ver\'ak \cite{MR3614655} treated initial velocities of arbitrary large size, constructing global weak $L^3$ solutions by decomposing the solution into the linear heat evolution of the initial velocity and a finite-energy remainder.
They also established subsequential stability under weak convergence of the initial data and uniqueness under additional regularity conditions.
Barker, Seregin, and \v{S}ver\'ak \cite{MR3902173} extended this caloric-remainder formulation to arbitrary initial data in the critical Lorentz space $L^{3,\infty}(\R^3)$, proving global weak existence and subsequential stability under weak-$*$ convergence of the initial data.
These decomposition methods provide the deterministic precedent for the present approach.
Here the caloric remainder is used to establish global continuation and pathwise uniqueness in a prescribed critical mild class, together with stability under perturbations of the initial velocity in the presence of state-dependent Wiener and compensated Poisson forcing.
Lai, Tang, and Xu \cite{MR5136351} developed a related caloric-remainder approach for the three-dimensional incompressible magnetohydrodynamics equations, constructing global weak solutions for arbitrarily large $L^3$ initial velocity and magnetic field, with uniqueness under additional $L^3$ regularity and short-time proximity assumptions.

These decomposition methods are therefore established tools in deterministic incompressible fluid theory.
For the problem considered here, the decomposition is used to obtain global evolution and uniqueness in a prescribed critical mild class, together with stability under perturbations of the initial velocity in the presence of state-dependent stochastic forcing.

A particularly direct comparison within the same class of stochastic tamed NSE driven simultaneously by multiplicative Wiener and compensated Poisson noise is provided by Podder and Kumar \cite{PodderKumar2026,PKGlobal2026}.
In Part I, they studied maximal local well-posedness and pathwise uniqueness in $L^p_\sigma(\R^3)$ for $p>3$, without requiring finite initial $L^2$-energy.
In Part II, they proved almost-sure global well-posedness for initial data in $L^p_\sigma(\R^3)\cap L^2_\sigma(\R^3)$, $p>3$, under additional coercivity and noise regularity assumptions, together with finite-energy estimates, positive-time Sobolev regularization, and continuous dependence in probability.
Both results exclude the endpoint $p=3$, while the global argument in the second begins with an $L^2$ estimate for the full velocity.
Therefore these results do not  provide the global critical mild conclusion for arbitrary initial data in $L^3_\sigma(\R^3)$.

\subsection{The global critical problem and its scope}

The preceding literature establishes several complementary well-posedness theories, but their conclusions do not directly yield the global critical result sought here.
Critical-space theories for the stochastic incompressible NSE provide local mild solvability.
Global theories for stochastic tamed equations and related critical variational frameworks exploit nonlinear dissipation in energy or Sobolev settings, while deterministic caloric-remainder methods provide global weak solutions for rough initial velocities.
None of these results, taken individually or through their stated conclusions, provides global pathwise well-posedness in the critical $L^3$ mild topology for arbitrary-size, possibly infinite-energy initial velocities under simultaneous state-dependent Wiener and compensated Poisson forcing.

This problem is not resolved by a formal passage to the endpoint $p=3$ or by applying an existing global energy theory.
At the endpoint, the cubic taming drift falls outside the range of the standard projected $L^p$ estimates, while an arbitrary velocity in $L^3_\sigma(\R^3)$ need not possess finite kinetic energy.
The deterministic caloric decomposition offers a natural way to separate the rough initial evolution from an energy-controlled correction, but the decomposition alone does not establish the required regularity of that correction in the stochastic critical mild setting.
 The existence of a global weak or variational solution does not automatically prove  the pathwise-uniqueness of the  critical evolution constructed on the prescribed stochastic basis.

The essential issue is therefore to reconcile critical mild solvability with dissipative global continuation within a single stochastic solution framework.
This requires a finite-energy remainder with zero initial data, whose estimates are independent of the critical norms used to localize the mild solution.
It also requires a positive-time Sobolev continuation that can be identified with the original critical evolution. It should preserve the stopping-time restart structure and the post-jump values generated by the Poisson forcing.
Without this compatibility, global bounds in an auxiliary Sobolev class do not by itself exclude the finite-time breakdown of the original critical mild solution.

We establish global pathwise well-posedness in the critical mild topology for arbitrary strongly $\mathcal F_0$-measurable $L^3_\sigma(\R^3)$-valued initial velocities.
The global argument uses Sobolev regularity of the stochastic coefficient outputs.
For fixed coefficients and driving noises, the global solutions depend continuously in probability on the initial velocity in $L^3_\sigma$.
The assumptions permit non-small coefficient bounds and infinite jump activity.
The resulting contribution lies in extending global well-posedness to the scaling-critical initial-data and solution topology, while retaining the original pathwise solution class throughout the global argument.

\subsection{Main results and proof strategy}

Under \cref{hyp:L}, a globally Lipschitz profile $g_N$ with $g_N(0)=0$ denoted as \cref{hyp:T1}, suffices for maximal local well-posedness of \cref{eq:tnse} in $\mathbb X_T$.
For every strongly $\mathcal F_0$-measurable $L^3_\sigma$-valued datum the maximal critical mild solution is pathwise unique on the prescribed stochastic basis.
Taking $g_N\equiv0$ gives the same local conclusion for the stochastic NSE with both noises.
The maximal solution satisfies the continuation criterion
\begin{equation*}
\P\left( \tau_*<\infty,\quad \int_0^{\tau_*}\|u(t)\|_5^5\dd t<\infty \right)=0.
\end{equation*}
The construction rests on deterministic convolution estimates for the quadratic convection and cubic taming terms, together with an $L^3$ stochastic Stokes estimate for simultaneous Wiener and compensated Poisson forcing.
The stochastic integration uses the Banach-space theories of van Neerven, Veraar, and Weis \cite{MR2330977} and Dirksen \cite{MR3265175}.
A predictable cutoff of the accumulated space--time $L^5$ norm makes the deterministic nonlinearities contractive on a short interval, while the stochastic coefficients are left unchanged.
The resulting solutions are joined at stopping times using pathwise uniqueness and restart identities that retain the terminal Poisson jump.
Localization on $\mathcal F_0$-events extends the finite-moment construction to all measurable initial velocities in \cref{thm:main-local}.

For global continuation we impose the coercive taming hypothesis \cref{hyp:T} and the additional coefficient hypothesis \cref{hyp:G}.
The usual energy estimate for the full velocity cannot be used at the initial time, since $L^3_\sigma(\R^3)$ is not contained in $L^2_\sigma(\R^3)$.
We therefore separate the heat evolution of the initial datum, writing $$ h(t)=S(t)u_0,\qquad R(t)=u(t)-h(t).
$$ The first difficulty is to show that $R$ has finite energy and this does not follow merely from the critical mild formulation.
\Cref{lem:remainder-energy-entrance} establishes, on every closed critical solution interval $[0,\tau]$, that $$ R\in D([0,\tau];L^2_\sigma) \cap L^2(0,\tau;H^1_\sigma) \cap L^4(0,\tau;L^4_\sigma), \qquad \lim_{t\downarrow0}\|R(t)\|_2=0, $$ almost surely.
The proof identifies the remainder with deterministic and stochastic convolutions whose forcing terms have the required energy-space integrability.
The convection cancellations and shifted taming coercivity then yield the estimates of \cref{prop:uniform-energy}.
On events where the initial $L^3$ norm is bounded, these estimates are independent of the critical norms of the solution and of the auxiliary stopping thresholds.
This uniformity is essential for passing to the maximal interval of existence.

The finite-energy estimates yield $H^1_\sigma$ regularity of the remainder at almost every positive time, but they do not by themselves provide a Sobolev evolution that can be identified with the critical solution.
By \cref{lem:admissible-h1}, there exist deterministic times $s_j\downarrow0$ at which $R(s_j)$ admits an $\mathcal F_{s_j}$-measurable $H^1_\sigma$-valued version on $\{s_j<\tau_*\}$.
On the $\mathcal F_{s_j}$-measurable events $A_{j,K,M}$, where both $\|u_0\|_3$ and $\|R(s_j)\|_{H^1}$ are bounded, the restarted velocity $h(s_j)+R(s_j)$, extended by zero outside the event, belongs to $L^3(\Omega,\mathcal F_{s_j};L^3_\sigma)$.
This permits the shifted Sobolev evolution to be compared with the original critical mild solution on the prescribed stochastic basis.

The coercivity of the shifted taming term and the quadratic $H^1$ growth bounds in \cref{hyp:G} provide global Sobolev estimates for the shifted evolution.
The event-localized identification argument, together with pathwise uniqueness and the stochastic restart identities, then identifies this evolution with the original remainder before the maximal lifetime.
The stochastic energy estimates retain the squared-jump contribution as a nondecreasing integral against the Poisson random measure.
Only quadratic $H^1$ integrability of the jump coefficient is needed for these Sobolev estimates; the quadratic and cubic $L^3$ jump-integrability conditions of \cref{hyp:L} remain in force.
No fourth $H^1$ jump moment is assumed.

For $0<\delta<T<\infty$, \cref{prop:positive-time-h1-persistence} bounds $\|R(t)\|_{H^1}$ uniformly over $\delta\le t\le T$ with $t<\tau_*$ on $\{\delta<\tau_*\}$, even when $\tau_*\le T$.
The embedding $H^1(\R^3)\hookrightarrow L^3(\R^3)\cap L^5(\R^3)$ and the heat estimates for $h$ then control the uniform $L^3$ and space--time $L^5$ norms of $u$ over the same times.
With $\delta=s_j$ chosen inside the first closed critical interval, local critical regularity supplies both bounds on $[0,s_j]$.
Thus \cref{lem:remainder-critical-transfer} and the space--time continuation criterion of \cref{thm:main-local} give $\P(\tau_*=\infty)=1$.

The probabilistically strong solution of \cref{thm:global-critical-wellposedness} is pathwise unique in the critical class.
On a single event of probability one, for every finite $T>0$, $$ u\in D([0,T];L^3_\sigma)\cap L^5(0,T;L^5_\sigma), $$ and the finite-energy regularity of $R$ established above holds on $[0,T]$.

For solutions driven by the same noises with the same coefficient fields, \cref{prop:global-critical-stability} gives continuous dependence for measurable initial data.
If $\|u_0^n-u_0\|_3\to0$ in probability, then for every finite $T>0$, $$ \sup_{0\le t\le T}\|u^n(t)-u(t)\|_3 +\|u^n-u\|_{L^5(0,T;L^5)}\longrightarrow0 \qquad\text{in probability}.
$$ We first estimate the third moment of the stopped difference conditionally on $\mathcal F_0$.
On the event where the limiting solution has critical norm at most $R$, an early exit of the pair at threshold $2R+1$ forces a difference of at least one at the exit time.
The conditional estimate controls this event for data converging in probability, including when the exit occurs at a jump.

Finally, the examples in \cref{Section 7} illustrate the scope of the global theorem.
In \cref{ex:infinite-energy-data}, explicit divergence-free initial velocities in $L^3_\sigma\setminus L^2_\sigma$ generate solutions that remain outside $L^2_\sigma$ at every finite time.
For random scalar multiples, this assertion holds on the event where the amplitude is nonzero.
In this example, the additional regularity of the initial field makes the heat increment $S(t)u_0-u_0$ finite-energy, and therefore $u(t)-u_0$ is also $L^2_\sigma$-valued.
The coefficient construction in \cref{ex:smoothing-mixed-noise} satisfies \cref{hyp:L} and \cref{hyp:G} and combines a state-dependent Wiener operator of infinite rank at a suitable state with infinite-activity jump noise having finite second and third moments but an infinite fourth moment.
These examples show that the global theory applies to genuinely infinite-energy initial data under nontrivial multiplicative Wiener and compensated Poisson perturbations.
\subsection{Organization of the paper.}
\cref{Section 2} specifies the equation and standing assumptions.
It also states the principal results, defines the critical solution class, and records the structural estimates.
\cref{Section 3} develops the deterministic and stochastic critical evolution estimates, and \cref{Section 4} establishes maximal local well-posedness, pathwise uniqueness, and continuation.
\cref{Section 5} proves energy regularity and uniform localized estimates for the caloric remainder.
\cref{Section 6} establishes positive-time Sobolev persistence and global critical well-posedness.
\cref{Section 7} proves continuous dependence and presents the examples.
\section{Functional setting, hypotheses, and main results}
\label{Section 2} We formulate the equation in the critical space $L^3_\sigma$ and specify the coefficient hypotheses for local and global well-posedness.
\subsection{Spatial setting and notations}\label{Sec2.1}

Unless otherwise specified, all spatial function spaces are defined over $\R^3$ and consist of $\R^3$-valued functions.
For $1<q<\infty$, let $$
\begin{aligned}
L^q_\sigma &:=\{v\in L^q(\R^3;\R^3): \diver v=0\text{ in }\mathcal D'(\R^3)\}, \\
E&:=L^3_\sigma,\qquad H^m_\sigma:=H^m(\R^3;\R^3)\cap L^2_\sigma, \quad m=0,1,2.
\end{aligned}
$$ We write $\|v\|_q=\|v\|_{L^q}$ and use the Euclidean and Frobenius norms for vectors and matrices, respectively.
The space $H^1_\sigma$ is equipped with $$ \|v\|_{H^1}^2 =\|v\|_2^2+\|\nabla v\|_2^2.
$$

Fix a viscosity $\nu>0$ and a taming parameter $N>0$.
Let $S(t)=e^{\nu t\Delta}$ denote the componentwise heat semigroup, and let $\LP$ denote the Leray projection.
Its Fourier multiplier and the duality identities used below are recorded in \cref{lem:projection}.
We adopt the tensor convention
\begin{equation}
\label{eq:2.1_tensor-convention}
(a\otimes b)_{ij}=b_i a_j,\qquad (\diver F)_i=\sum_{j=1}^3\partial_jF_{ij}, \qquad (\nabla c)_{ij}=\partial_jc_i.
\end{equation}
Consequently, for smooth solenoidal $a$, $\diver(a\otimes b)=(a\cdot\nabla)b$.

\subsection{Probabilistic setting, initial data and hypothesis}\label{Sec2.2}
Let $(\Omega,\mathcal F,(\mathcal F_t)_{t\ge0},\P)$ be a complete filtered probability space satisfying the usual conditions.
Let $\U$ be a separable real Hilbert space and $W$ an $(\mathcal F_t)$-cylindrical Wiener process on $\U$.
Let $(Z,\mathcal Z,\mu)$ be a $\sigma$-finite measure space, and let $\PN$ be an $(\mathcal F_t)$-Poisson random measure on $[0,\infty)\times Z$ with compensator $\mathrm dt\,\mu(\mathrm dz)$.
Its compensated measure is $$ \tPN(\mathrm dt,\mathrm dz) :=\PN(\mathrm dt,\mathrm dz) -\mathrm dt\,\mu(\mathrm dz).
$$ 
The process $W$ and $\PN$ are assumed  independent. 
Write $L^0(\Omega,\mathcal F_0;E)$ for the strongly $\mathcal F_0$-measurable $E$-valued random variables, identified up to almost sure equality.
The pathwise theory assumes $u_0\in L^0(\Omega,\mathcal F_0;E)$.
The cubic moment condition is
\begin{equation}
\label{eq:initial}
u_0\in L^3(\Omega,\mathcal F_0;E), \qquad\text{equivalently}\qquad \E\|u_0\|_3^3<\infty .
\end{equation}
We impose \eqref{eq:initial} where a finite third moment is needed.

\paragraph{Taming hypothesis for local theory \textup{\ensuremath{(\mathrm{T}_{\mathrm{loc}})}}}
\hypothesislabel{T1}{hyp:T1} The local critical theory uses a function $g_N:[0,\infty)\to\R$ satisfying, for some $C_g<\infty$,
\begin{equation*}
g_N(0)=0,\qquad |g_N(r)-g_N(s)|\le C_g|r-s| \quad (r,s\ge0).
\end{equation*}
For the global theory we impose the following stronger condition.

\paragraph{Taming hypothesis for global theory \textup{(T)}}
\hypothesislabel{T}{hyp:T} We assume
\begin{equation*}
\begin{aligned}
g_N\in C^1([0,\infty))&,\qquad g_N(0)=0,\qquad 0\le g_N'(r)\le C_g\quad(r\ge0), \\
c_{N,\nu}:&= \sup_{r\ge0}\left(\frac r\nu-2g_N(r)\right)<\infty, \qquad C_g<\infty.
\end{aligned}
\end{equation*}
The constants may depend on $N$ and $\nu$.
In particular, these conditions imply
\begin{equation}
\label{eq:2.2_g-consequences}
c_{N,\nu}\ge0,\qquad 0\le g_N(r)\le C_g r,\qquad g_N(r)\ge\frac{r}{2\nu}-\frac{c_{N,\nu}}2.
\end{equation}
The uniform upper bound on $g_N'$ is assumed independently of the monotonicity and growth properties of $g_N$.

\paragraph{ Noise coefficient hypothesis for the local theory
\textup{(L)}} \hypothesislabel{L}{hyp:L} For a Banach space $X$, let $\gamma(\U;X)$ denote the space of $\gamma$-radonifying operators from $\U$ into $X$.
For a finite-rank operator $\Phi v=\sum_{k=1}^m(v,e_k)_\U x_k$, where $(e_k)$ is an orthonormal family in $\U$, the norm is defined by $$ \|\Phi\|_{\gamma(\U;X)}^2 =\E_\gamma\left\| \sum_{k=1}^m\gamma_kx_k \right\|_X^2, $$ where $(\gamma_k)$ are independent standard real Gaussian random variables.
The space $\gamma(\U;X)$ is obtained by completing the finite-rank operators with respect to this norm.

Set $$ \mathcal J_3 :=L^2(Z,\mu;E)\cap L^3(Z,\mu;E), \qquad \|G\|_{\mathcal J_3} :=\|G\|_{L^2(Z;E)}+\|G\|_{L^3(Z;E)}, $$ where the spaces on the right are Bochner spaces.

Let $\mathscr P$ denote the predictable $\sigma$-field on $\Omega\times[0,\infty)$.
We assume that the coefficient maps $$ \Sigma:\Omega\times[0,\infty)\times E \longrightarrow\gamma(\U;E), \qquad \mathcal G:\Omega\times[0,\infty)\times E\times Z \longrightarrow E $$ are respectively $\mathscr P\otimes\mathcal B(E)$-measurable and $\mathscr P\otimes\mathcal B(E)\otimes\mathcal Z$-measurable, with strong measurability in the indicated Banach spaces.

For every $T<\infty$, there exist finite deterministic constants $L_T,K_T$ such that, outside a common $\P\otimes\mathrm dt$-null set in $\Omega\times[0,T]$, the following inequalities hold for every $v,w\in E$:
\begin{align*}
&\|\Sigma(t,v)-\Sigma(t,w)\|_{\gamma(\U;E)} +\|\mathcal G(t,v,\cdot)-\mathcal G(t,w,\cdot)\|_{\mathcal J_3} \le L_T\|v-w\|_3, \\
&\|\Sigma(t,0)\|_{\gamma(\U;E)} +\|\mathcal G(t,0,\cdot)\|_{\mathcal J_3} \le K_T.
\end{align*}
The exceptional null set is independent of $v,w$.
Consequently,
\begin{equation}
\label{eq:local-growth}
\|\Sigma(t,v)\|_{\gamma(\U;E)} +\|\mathcal G(t,v,\cdot)\|_{\mathcal J_3} \le K_T+L_T\|v\|_3.
\end{equation}
In particular, there exists a finite deterministic constant $C_T$ such that
\begin{equation*}
\begin{aligned}
&\|\Sigma(t,v)\|_{\gamma(\U;E)}^3 +\left(\int_Z \|\mathcal G(t,v,z)\|_3^2 \,\mu(\mathrm dz)\right)^{3/2} +\int_Z\|\mathcal G(t,v,z)\|_3^3 \,\mu(\mathrm dz) \le C_T(1+\|v\|_3^3).
\end{aligned}
\end{equation*}
Indeed, the left-hand side is bounded by the cube of the left-hand side of \eqref{eq:local-growth}, and $(a+b)^3\le4(a^3+b^3)$ for $a,b\ge0$.
\paragraph{Additional noise coefficient hypothesis for the
global theory \textup{(G)}} \hypothesislabel{G}{hyp:G} For the global theory, we additionally assume that the coefficients admit compatible strongly $\mathscr P\otimes\mathcal B(E)$-measurable versions with values in $\HS(\U;H^1_\sigma)$ and $L^2(Z,\mu;H^1_\sigma)$, respectively.
The latter is represented by a strongly $\mathscr P\otimes\mathcal B(E)\otimes\mathcal Z$-measurable $H^1_\sigma$-valued map.
Here $\HS$ denotes the Hilbert--Schmidt class, with norm $$ \|\Phi\|_{\HS(\U;H^1_\sigma)}^2 :=\sum_k\|\Phi e_k\|_{H^1}^2, $$ where $(e_k)$ is an orthonormal basis of $\U$.

For every $T<\infty$, there exist finite deterministic constants $L_T^{(1)},K_T^{(1)}$ such that, outside a common $\P\otimes\mathrm dt$-null set in $\Omega\times[0,T]$, the following inequalities hold for all $v,w\in E$:
\begin{align*}
&\|\Sigma(t,v)-\Sigma(t,w)\|_{\HS(\U;H^1_\sigma)} +\|\mathcal G(t,v,\cdot)-\mathcal G(t,w,\cdot)\|_{L^2(Z;H^1_\sigma)} \le L_T^{(1)}\|v-w\|_3, \\
&\|\Sigma(t,0)\|_{\HS(\U;H^1_\sigma)} +\|\mathcal G(t,0,\cdot)\|_{L^2(Z;H^1_\sigma)} \le K_T^{(1)}.
\end{align*}
Compatibility means that, under the continuous embedding $H^1_\sigma\hookrightarrow E$, these versions agree with the coefficient maps in \cref{hyp:L} as elements of $\gamma(\U;E)$ and $L^2(Z,\mu;E)$, respectively.
For each $T<\infty$, this and the preceding inequalities hold for every state outside a common $\P\otimes\mathrm dt$-null set independent of that state.
For the jump coefficient, equality is understood $\mu$-almost everywhere in $z$.

By the triangle inequality, there exists a finite deterministic constant $C_T^{(1)}$ such that
\begin{equation}
\|\Sigma(t,v)\|_{\HS(\U;H^1_\sigma)}^2 +\int_Z\|\mathcal G(t,v,z)\|_{H^1}^2 \,\mu(\mathrm dz) \le C_T^{(1)}(1+\|v\|_3^2).
\end{equation}
The $L^3(Z;E)$ integrability condition in \cref{hyp:L} remains in force.

The additional regularity required in \cref{hyp:G} concerns the values of the stochastic coefficients and imposes no spatial $L^2$ condition on the initial datum $u_0$.
The cubic upper bounds used in the local construction follow from \cref{hyp:T1}.
The lower bounds for the taming drift use \cref{hyp:T}.

\subsection{Caloric reference and remainder.}
For a critical solution $u$ of \eqref{eq:tnse}, define
\begin{equation*}
h(t):=S(t)u_0,\qquad R(t):=u(t)-h(t),
\end{equation*}
where $R$ is defined on the interval of existence of $u$.
Subtracting the heat equation for $h$ from \eqref{eq:tnse} gives, on that interval and in the distributional or mild sense,
\begin{equation*}
\begin{cases}
\begin{aligned}
\mathrm dR =&\bigl[ \nu\Delta R -\LP\diver\bigl((h+R)\otimes(h+R)\bigr) -\LP\TN(h+R) \bigr]\dd t \\
&+\Sigma(t,h+R_-)\dd W(t) +\int_Z\mathcal G(t,h+R_-,z) \,\tPN(\mathrm dt,\mathrm dz),
\end{aligned}
\\
\diver R=0. \qquad R(0)=0,
\end{cases}
\end{equation*}

Since $h$ is continuous in $E$, we have $u(t-)=h(t)+R(t-)$ for every positive time in the interval of existence.
The decomposition gives $R(0)=0$ in $E$.
Under \cref{hyp:T} and \cref{hyp:G}, \cref{lem:remainder-energy-entrance} shows that $R$ also has a zero initial data in $L^2_\sigma$.
No spatial $L^2$ integrability is assumed for $h$.

For $1<q<\infty$, let $q'=q/(q-1)$ and $$ W^{-1,q}:=(W^{1,q'}(\R^3;\R^3))^*,\qquad H^{-1}=W^{-1,2}.
$$ Whenever a duality pairing is represented by an integrable product, we use the corresponding Lebesgue integral.
The notation $\langle\cdot,\cdot\rangle_{H^{-1},H^1}$ denotes distributional duality.
Denote by $D([0,T];E)$ the space of right-continuous $E$-valued functions admitting left limits.
Set
\begin{equation}
\label{eq:critical-space}
\mathbb X_T :=D([0,T];E)\cap L^5(0,T;L^5_\sigma), \qquad \|v\|_{\mathbb X_T} :=\sup_{0\le t\le T}\|v(t)\|_3 +\left(\int_0^T\|v(t)\|_5^5\dd t\right)^{1/5}.
\end{equation}
The designation \emph{critical} refers to the classical Navier--Stokes scaling of these two norms.
Scale invariance of the fixed taming function or of the stochastic coefficients is not assumed.
\subsection{Critical mild solutions and main results}

We start by the following definitions,
\begin{definition}
[Critical mild solutions and closed localization] \label{def:critical-solution} Assume \cref{hyp:T1} and \cref{hyp:L}, and let $u_0\in L^0(\Omega,\mathcal F_0;E)$.
Let $\tau$ be a bounded stopping time with $\P(\tau>0)=1$.
A \emph{closed local critical mild solution on $[0,\tau]$} is an adapted $E$-valued process $u$, defined for all $t\ge0$, satisfying the following conditions.

\begin{enumerate}
[label=\textup{(\roman*)},leftmargin=*] \item The process $u$ is c\`adl\`ag and stopped at $\tau$: $u(t)=u(t\wedge\tau)$ for all $t\ge0$, almost surely.
Moreover, $u(0)=u_0$ and
\begin{equation}
\sup_{0\le s\le T\wedge\tau}\|u(s)\|_3 +\left(\int_0^{T\wedge\tau}\|u(s)\|_5^5\dd s\right)^{1/5} <\infty \quad\text{almost surely, for every }T<\infty .
\end{equation}
The $L^5$ condition is understood in the Bochner sense on the stopped interval.
No $L^5$ condition is imposed on the constant continuation after $\tau$.

\item Set $u(0-)=u_0$ and $\chi_\tau(s):=\mathds1_{\{0<s\le\tau\}}$.
For each deterministic $t\ge0$, form
\begin{equation}
\label{eq:2.4_stop_Con}
\begin{aligned}
Y_\tau(t):=&S(t)u_0 -\int_0^t\chi_\tau(s)S(t-s) \LP\diver(u(s)\otimes u(s))\dd s -\int_0^t\chi_\tau(s)S(t-s)\LP\TN(u(s))\dd s \\
&+\int_0^t\chi_\tau(s)S(t-s) \Sigma(s,u(s-))\dd W(s) +\int_0^t\int_Z\chi_\tau(s)S(t-s) \mathcal G(s,u(s-),z)\, \tPN(\mathrm ds,\mathrm dz).
\end{aligned}
\end{equation}
We require an $E$-valued c\`adl\`ag modification of this expression such that, outside a single null set,
\begin{equation}
\label{eq:closed-mild-equality}
u(t\wedge\tau)=Y_\tau(t\wedge\tau) \qquad\text{for every }t\ge0 .
\end{equation}
The deterministic convolutions are interpreted initially in $W^{-1,5/2}$ and $L^{5/3}$, respectively, and therefore in $\mathcal S'(\R^3;\R^3)$.
The stochastic integrals are the $E$-valued It\^o integrals, defined by localization.
\end{enumerate}

A \emph{local critical mild solution with lifetime $\tau_*$} is an adapted process $u$ on $[0,\tau_*)$, where $\tau_*>0$ is a stopping time, admitting bounded stopping times $\tau_n\uparrow\tau_*$ almost surely, with $\tau_n<\tau_*$ on $\{\tau_*<\infty\}$, such that $(u(\cdot\wedge\tau_n),\tau_n)$ is a closed local critical mild solution for each $n$.

It is \emph{maximal} if it has no extension $(\bar u,\bar\tau)$ in this class, on the same stochastic basis and with the same initial datum and noises, satisfying $$ \bar\tau\ge\tau_*\ \text{almost surely},\qquad \P(\bar\tau>\tau_*)>0,\qquad \bar u(t)=u(t)\quad(0\le t<\tau_*)\ \text{almost surely}.
$$ Equality of solutions on a stochastic interval means equality simultaneously at all times there, outside one null set.
In particular, comparison of two closed local solutions includes the time $\tau\wedge\bar\tau$.
A solution with $\tau_*=\infty$ almost surely is called global.
\end{definition}

\paragraph{Probabilistically strong solutions and the mild formulation}
All solutions are constructed on the prescribed stochastic basis with the specified driving noises and coefficient fields, including any $\omega$-dependence of those fields and no auxiliary randomness is introduced.
The term \emph{probabilistically strong critical solution} thus refers to solutions in the sense of \cref{def:critical-solution}, with the corresponding mild identities.
It does not imply an analytically strong $E$-valued equation or $L^1(0,\tau;W^{2,3})$ regularity.
No finite-energy condition is imposed on the full velocity, since $L^3(\R^3)\not\subset L^2(\R^3)$.
For a c\`adl\`ag path $v$, we write $$ \tilde\Delta v(t):=v(t)-v(t-), \qquad t>0, $$ for its temporal jump.
The principal results are stated below.
Their proofs are given in \cref{Section 4}, \cref{Section 6} and \cref{Section 7} respectively.
\begin{theorem}
[Maximal critical local well-posedness] \label[theorem]{thm:main-local} Assume the standing probabilistic framework, \cref{hyp:T1} and \cref{hyp:L}, and let $u_0\in L^0(\Omega,\mathcal F_0;E)$.
Under these assumptions, \eqref{eq:tnse} admits a maximal critical mild solution $(u, \{\tau_n\}_{n \ge1}, \tau_*)$ on the prescribed stochastic basis. The solution is probabilistically strong in the sense of \cref{def:critical-solution}, and its maximal lifetime satisfies $\P(\tau_*>0)=1$.
Pathwise uniqueness holds in the critical class.
More precisely, if $(u,\tau)$ and $(v,\sigma)$ are closed local critical mild solutions with the same initial datum, coefficients, Wiener process, and Poisson random measure, then $$ \P\bigl( u(t)=v(t) \text{ for every }0\le t\le\tau\wedge\sigma \bigr)=1.
$$ Consequently, maximal solutions with identical data coincide, and their lifetimes agree almost surely.
The maximal lifetime satisfies $$ \P\left(\tau_*<\infty,\quad \int_0^{\tau_*}\|u(t)\|_5^5\,\dd t<\infty\right)=0.
$$ Every local solution $v$ with the same initial datum, coefficients and driving noises is a restriction of $u$.
If $\zeta$ denotes its lifetime, then almost surely on $\{\zeta<\infty,\ \int_0^\zeta\|v(t)\|_5^5\,\dd t<\infty\}$, $$ \zeta<\tau_*,\qquad \lim_{t\uparrow\zeta}v(t)=u(\zeta)\quad\text{in }E.
$$ The maximal solution continues from this state on the prescribed basis.
In particular, taking $g_N\equiv0$ gives maximal local well-posedness and the same continuation criterion for the stochastic NSE on $\R^3$ with the Wiener and compensated Poisson coefficients of \cref{hyp:L}.
\end{theorem}

\begin{theorem}
[Global existence and pathwise uniqueness in the critical class] \label[theorem]{thm:global-critical-wellposedness} Assume the standing probabilistic framework, the coefficient measurability conditions, \cref{hyp:T}, \cref{hyp:L} and \cref{hyp:G}, and $u_0\in L^0(\Omega,\mathcal F_0;E)$.
The maximal critical solution of \cref{thm:main-local} satisfies
\begin{equation}
\label{eq:6.5_global}
\P(\tau_*=\infty)=1.
\end{equation}
It is a probabilistically strong critical mild solution on the prescribed stochastic basis and is pathwise unique in the class of \cref{def:critical-solution}, for the same initial datum, driving noises, and coefficient fields.

Outside a single null set, for every $0<T<\infty$,
\begin{equation}
\label{eq:6.5_global_2}
u\in D([0,T];L^3_\sigma)\cap L^5(0,T;L^5_\sigma), \qquad u(0)=u_0.
\end{equation}
The remainder satisfies
\begin{equation}
\label{eq:6.5_global_3}
\begin{aligned}
R&\in D([0,T];L^2_\sigma) \cap L^2(0,T;H^1_\sigma) \cap L^4(0,T;L^4_\sigma), \\
R(0)&=0,
\end{aligned}
\end{equation}
and
\begin{equation}
\label{eq:6.5_global_4}
\lim_{t\downarrow0}\|u(t)-u_0\|_3=0, \qquad \lim_{t\downarrow0}\|R(t)\|_2=0.
\end{equation}
For every $0<\delta<T<\infty$,
\begin{equation}
\label{eq:6.5_global_5}
\begin{aligned}
R&\in D([\delta,T];H^1_\sigma) \cap L^2(\delta,T;H^2_\sigma), \qquad \int_\delta^T\int_{\R^3} |u(t)|^2|\nabla R(t)|^2\dd x\dd t<\infty.
\end{aligned}
\end{equation}
These assertions hold simultaneously for all such $T$ and $\delta$.

For $K\ge1$, set $\Omega_K=\{\|u_0\|_3\le K\}$.
For every $T>0$ and $A\in\mathcal F_0$ with $A\subseteq\Omega_K$,
\begin{equation}
\label{eq:6.5_global_6}
\E[\mathds1_A\mathcal E_R(T)]\le\Gamma_{T,K}\P(A).
\end{equation}
For $0<\delta<T$, define $$ \mathcal H_{\delta,T} =\sup_{\delta\le t\le T}\|R(t)\|_{H^1}^2 +\int_\delta^T\mathcal D_{h(t)}(R(t))\dd t .
$$ For the same $K$ and $A$ one has
\begin{equation}
\label{eq:6.5_global_positive}
\E[\mathds1_A\mathcal H_{\delta,T}] \le C_{\delta,T}(1+K^4)\P(A).
\end{equation}
If, in addition, \eqref{eq:initial} holds, then, with $Q_T=(0,T)\times\R^3$,
\begin{align}
\E\!\left[(\mathcal E_R(T))^{3/4}+\|u\|_{L^4(Q_T)}^3\right] &\le C_T\bigl(1+\E\|u_0\|_3^3\bigr),
\label{eq:6.5_global_fractional}\\
\E[\mathcal H_{\delta,T}^{3/4}] &\le C_{\delta,T}\bigl(1+\E\|u_0\|_3^3\bigr).
\label{eq:6.5_positive_fractional}
\end{align}
Here $C_T$ and $C_{\delta,T}$ are finite deterministic constants depending only on the displayed times and the fixed structural and coefficient bounds.
Event inclusions are understood modulo null sets.
For every $u_0\in L^0(\Omega,\mathcal F_0;E)$, the jump identity holds almost surely, simultaneously for all $t>0$,
\begin{equation}
\label{eq:6.5_global_7}
\tilde \Delta u(t)=\tilde\Delta R(t) =\int_Z\mathcal G(t,u(t-),z)\,\PN(\{t\},\dd z).
\end{equation}
The identity for $\tilde \Delta R(t)$ holds in $H^1_\sigma$, and the identity for $\tilde \Delta u(t)$ holds in $L^3_\sigma$.
No spatial $L^2$ or $H^1$ assumption and no smallness condition are imposed on $u_0$.
\end{theorem}
The quantities $\mathcal E_R$, $\Gamma_{T,K}$, and $\mathcal D_h$ are defined in \eqref{eq:5.4_rem_eng}, \eqref{eq:5.4_loc_eng1}, and \eqref{eq:6.1_h1_diss}, respectively.

Now for two such global solutions $v,w$, define, for $t\ge0$, $$
\begin{aligned}
\mathcal C_{v,w}(t) &:=\sup_{0\le r\le t} \bigl(\|v(r)\|_3+\|w(r)\|_3\bigr) +\left(\int_0^t \bigl(\|v(r)\|_5^5+\|w(r)\|_5^5\bigr)\dd r\right)^{1/5}, \\
D_t(v,w) &:=\sup_{0\le r\le t}\|v(r)-w(r)\|_3 +\left(\int_0^t\|v(r)-w(r)\|_5^5\dd r\right)^{1/5}.
\end{aligned}
$$ In particular, $D_0(v,w)=\|v_0-w_0\|_3$.
For $T,M>0$, set $$ \sigma_{M,T}(v,w) :=T\wedge\inf\{t\ge0:\mathcal C_{v,w}(t)\ge M\}, \qquad \inf\varnothing:=\infty.
$$ For a single global solution $u$, write $$ \mathcal C_u(t) :=\sup_{0\le r\le t}\|u(r)\|_3 +\left(\int_0^t\|u(r)\|_5^5\dd r\right)^{1/5}.
$$
\begin{proposition}
[Global continuous dependence on finite horizons] \label{prop:global-critical-stability} Assume \cref{hyp:T}, \cref{hyp:L}, and \cref{hyp:G}.
Let $v,w$ be the global critical mild solutions with respective data $v_0,w_0\in L^0(\Omega,\mathcal F_0;E)$ on the same stochastic basis, with the same coefficient fields, Wiener process and Poisson random measure.
For $T,M>0$, put $\sigma=\sigma_{M,T}(v,w)$ and $d=\|v_0-w_0\|_3$.
There is a deterministic $C_{T,M}\ge1$, depending only on $T,M$ and the fixed structural and coefficient bounds, such that almost surely
\begin{equation}
\label{eq:7.1_closed-1}
d^3\le \E\!\left[D_\sigma(v,w)^3\mid\mathcal F_0\right] \le C_{T,M}d^3 .
\end{equation}
In particular, $\E D_\sigma(v,w)^3<\infty$ if and only if $\E d^3<\infty$.
In that case $\E D_\sigma(v,w)^3\le C_{T,M}\E d^3$.
The supremum includes the terminal value, also when $\sigma=0$ or the threshold is reached at a jump.

For two such solutions $v,u$, let $a_\varepsilon=\min\{\varepsilon,1\}$.
For every $T,R,\varepsilon>0$,
\begin{equation}
\label{eq:7.1_prob_est}
\P\{D_T(v,u)>\varepsilon\} \le \P\{\mathcal C_u(T)>R\} +\E\!\left[\min\!\left\{1, \frac{C_{T,\,2R+1}\|v_0-u_0\|_3^3}{a_\varepsilon^3} \right\}\right].
\end{equation}
Consequently, if $u_0^n,u_0\in L^0(\Omega,\mathcal F_0;E)$ and $\|u_0^n-u_0\|_3\to0$ in probability, then the corresponding solutions with the same coefficients and driving noises satisfy, for every $T>0$,
\begin{equation}
\label{eq:7.1_global-convergence}
\sup_{0\le t\le T}\|u^n(t)-u(t)\|_3 +\|u^n-u\|_{L^5(0,T;L^5)} \longrightarrow0 \qquad\text{in probability}.
\end{equation}
\end{proposition}

\paragraph{Admissibility of the stopped mild formulation}
We record the measurability and integrability properties of the stopped processes appearing in \eqref{eq:2.4_stop_Con}, which justify the stochastic integrals at deterministic times.
The existence of a common c\`adl\`ag realization is established in \cref{Section 4}.
The processes $u_-$ and $\chi_\tau$ are predictable. On $(0,\infty)$, the latter is adapted and left-continuous.
By \eqref{eq:local-growth} and boundedness of c\`adl\`ag paths, for every $T<\infty$,
\begin{equation}
\int_0^T\chi_\tau(s)\left[ \|\Sigma(s,u(s-))\|_{\gamma(\U;E)}^2 +\sum_{p=2,3}\int_Z\|\mathcal G(s,u(s-),z)\|_3^p\,\mu(\mathrm dz) \right]\dd s<\infty \quad\text{a.s.}
\end{equation}
The space $E$ is a closed UMD subspace of $L^3$ with type $2$.
Consequently, \cite[Corollary 3.10(1)]{MR2330977} and \cite[Theorem 1.1, $p=q=3$]{MR3265175}, together with localization, provide the stochastic integrals at each deterministic $t$.
For the Poisson integral, set $$ \Psi(s,z):=\chi_\tau(s)\mathcal G(s,u(s-),z).
$$ Minkowski's inequality gives $$ \left\|\left(\int_0^T\!\int_Z|\Psi|^2\,\mu(\mathrm dz)\dd s\right)^{1/2}\right\|_3 \le\left(\int_0^T\!\int_Z\|\Psi\|_3^2\,\mu(\mathrm dz)\dd s\right)^{1/2}; $$ The remaining integrability term is $\int_0^T\!\int_Z\|\Psi\|_3^3\,\mu(\mathrm dz)\dd s$.
For each fixed $t$, the boundedness of $S(t-s)$ preserves these integrability conditions.
Similarly, \cref{lem:taming} and \cref{lem:projection} give $$ \chi_\tau u\otimes u\in L^{5/2}(0,T;L^{5/2}),\qquad \chi_\tau\TN(u)\in L^{5/3}(0,T;L^{5/3}), $$ which justify the stated distributional drift convolutions.
The definition requires a common $E$-valued c\`adl\`ag realization, whose existence is established in the subsequent construction.
The factor $\chi_\tau$ includes the Poisson jump at the terminal time.
The process $Y_\tau$ is defined at deterministic times and subsequently evaluated at $t\wedge\tau$; consequently, the kernel $S(t\wedge\tau-s)$ is not inserted into an It\^o integrand.
The definition of maximality does not incorporate a continuation criterion.
\subsection{Taming estimates and solenoidal cancellations}
\begin{lemma}
[Taming algebra]\label[lemma]{lem:taming} Under \cref{hyp:T} all the following assertions hold.
The size estimate $|\TN(a)|\le C_g|a|^3$, \eqref{eq:cubic-difference}, and the two upper norm bounds in part \textup{(iii)} also hold under \cref{hyp:T1} alone.

\begin{enumerate}
[label=\textup{(\roman*)},leftmargin=*] \item The map $\TN:\R^3\to\R^3$ is continuously differentiable, and for $a,\eta\in\R^3$,
\begin{equation}
D\TN(a)\eta =g_N(|a|^2)\eta+2g_N'(|a|^2)(a\cdot\eta)a .
\end{equation}
In particular,
\begin{align}
|\TN(a)|&\le C_g|a|^3,\qquad \|D\TN(a)\|_{\mathrm{op}}\le3C_g|a|^2,\label{eq:taming-growth} \\
\eta\cdot D\TN(a)\eta &\ge g_N(|a|^2)|\eta|^2 \ge\left(\frac{|a|^2}{2\nu} -\frac{c_{N,\nu}}2\right)|\eta|^2 .
\end{align}
The derivative matrix is symmetric and positive semidefinite.

\item For all $a,b\in\R^3$,
\begin{equation}
\label{eq:cubic-difference}
|\TN(a)-\TN(b)| \le\frac{3C_g}{2}(|a|^2+|b|^2)|a-b|.
\end{equation}
For $d=a-b$, the following identity holds:
\begin{equation}
\label{eq:monotonicity-identity}
\begin{aligned}
(\TN(a)-\TN(b))\cdot d =&\frac{g_N(|a|^2)+g_N(|b|^2)}2\,|d|^2 +\frac{g_N(|a|^2)-g_N(|b|^2)}2 (|a|^2-|b|^2).
\end{aligned}
\end{equation}
Consequently, $\TN$ is monotone:
\begin{equation}
\label{eq:monotone}
(\TN(a)-\TN(b))\cdot(a-b)\ge0.
\end{equation}
It also satisfies
\begin{equation}
\label{eq:quantitative-monotonicity}
\begin{aligned}
(\TN(a)-\TN(b))\cdot d &\ge\frac{|a|^2+|b|^2}{4\nu}|d|^2 -\frac{c_{N,\nu}}2|d|^2 \ge\frac{|d|^4}{8\nu}-\frac{c_{N,\nu}}2|d|^2 .
\end{aligned}
\end{equation}

\item For $3\le q<\infty$, the associated Nemytskii map sends $L^q$ into $L^{q/3}$ and is locally Lipschitz:
\begin{align}
\|\TN(v)\|_{q/3} &\le C_g\|v\|_q^3, \\
\|\TN(v)-\TN(w)\|_{q/3} &\le\frac{3C_g}{2} (\|v\|_q^2+\|w\|_q^2)\|v-w\|_q .
\end{align}
If $v,w\in L^4$ and $d=v-w\in L^2$, then
\begin{equation}
\label{eq:integrated-taming}
\begin{aligned}
\int_{\R^3}(\TN(v)-\TN(w))\cdot d\dd x &\ge\frac1{4\nu}\int_{\R^3}(|v|^2+|w|^2)|d|^2\dd x -\frac{c_{N,\nu}}2\|d\|_2^2 \\
&\ge\frac1{8\nu}\|d\|_4^4 -\frac{c_{N,\nu}}2\|d\|_2^2 .
\end{aligned}
\end{equation}
The integral on the left remains nonnegative under the sole assumption $v,w\in L^4$.
\end{enumerate}
\end{lemma}

\begin{proof}
Under \cref{hyp:T1}, $|g_N(r)|\le C_g r$, which gives $|\TN(a)|\le C_g|a|^3$.
Write $g_a=g_N(|a|^2)$.
For $a,b\in\R^3$ we have $$ 2\bigl(\TN(a)-\TN(b)\bigr) =(g_a+g_b)(a-b)+(g_a-g_b)(a+b).
$$ Since $|g_a-g_b|\le C_g|(a+b)\cdot(a-b)|$ and $|a+b|^2\le2(|a|^2+|b|^2)$, this proves \eqref{eq:cubic-difference}.
Under \cref{hyp:T}, the chain rule and \eqref{eq:2.2_g-consequences} give the differential bounds in part \textup{(i)}, also at $a=0$.
Identity \eqref{eq:monotonicity-identity} follows by expansion.
Both terms are nonnegative since $g_N$ is nonnegative and nondecreasing.
Moreover, \eqref{eq:2.2_g-consequences} bounds its first coefficient below by $(|a|^2+|b|^2)/(4\nu)-c_{N,\nu}/2$.
This proves \eqref{eq:monotone} and \eqref{eq:quantitative-monotonicity}, using $|a-b|^2\le2(|a|^2+|b|^2)$ for the last inequality.

H\"older's inequality with $1/(q/3)=1/(q/2)+1/q$, applied to \eqref{eq:taming-growth} and \eqref{eq:cubic-difference}, proves the Nemytskii bounds.
For $v,w\in L^4$, the taming difference belongs to $L^{4/3}$ and admits a duality pairing with $d\in L^4$.
Moreover, $(|v|^2+|w|^2)|d|^2\in L^1$.
Integrating \eqref{eq:monotone} proves nonnegativity, and, when also $d\in L^2$, integrating \eqref{eq:quantitative-monotonicity} gives \eqref{eq:integrated-taming}.
\end{proof}
\begin{lemma}
[Leray projection and solenoidal cancellations] \label{lem:projection} Let $1<q<\infty$ and $q'=q/(q-1)$.
The Fourier multiplier $$ p_{ij}(\xi) :=\delta_{ij}-\frac{\xi_i\xi_j}{|\xi|^2}, \qquad \xi\in\R^3\setminus\{0\}, $$ defines a bounded projection $\LP:L^q\to L^q_\sigma$.
It extends boundedly to $W^{1,q}$ and $W^{-1,q}$, and is the orthogonal projection onto $L^2_\sigma$ when $q=2$.

For $f\in L^q$ and $\varphi\in L^{q'}$,
\begin{equation}
\label{eq:leray-duality}
\int_{\R^3}\LP f\cdot\varphi\dd x = \int_{\R^3}f\cdot\LP\varphi\dd x.
\end{equation}
In particular, if $\varphi\in L^{q'}_\sigma$, then
\begin{equation}
\label{eq:remove-projection}
\int_{\R^3}\LP f\cdot\varphi\dd x = \int_{\R^3}f\cdot\varphi\dd x.
\end{equation}
Moreover, $\LP$ commutes with the heat semigroup and distributional differentiation: $$ \LP S(t)f=S(t)\LP f, \qquad \LP\partial_j f=\partial_j\LP f, $$ where the second identity holds in $W^{-1,q}$.
Consequently,
\begin{equation}
\label{eq:negative-duality}
\langle\LP F,\varphi\rangle_{H^{-1},H^1} = \langle F,\varphi\rangle_{H^{-1},H^1}, \qquad F\in H^{-1},\quad \varphi\in H^1_\sigma.
\end{equation}

For $F\in L^q(\R^3;\R^{3\times3})$, define $\LP^{(1)}F$ by $ (\LP^{(1)}F)_{\cdot j} :=\LP(F_{\cdot j}), \ j=1,2,3.$ Then
\begin{equation}
\label{eq:tensor-projection}
\LP\diver F = \diver(\LP^{(1)}F) \qquad\text{in }W^{-1,q}.
\end{equation}
Under the tensor convention \eqref{eq:2.1_tensor-convention}, smooth vector fields $a,b$ satisfy
\begin{equation}
\label{eq:tensor-product}
\diver(a\otimes b) = (a\cdot\nabla)b+b\,\diver a.
\end{equation}

If $a\in L^3_\sigma$ and $b,c\in H^1$, then $$ B(a,b,c) := \int_{\R^3}(a\cdot\nabla)b\cdot c\dd x $$ is absolutely convergent and satisfies
\begin{align}
|B(a,b,c)| &\le \|a\|_3\|\nabla b\|_2\|c\|_6 \le C_S\|a\|_3\|\nabla b\|_2\|\nabla c\|_2,
\label{eq:trilinear-bound}\\
B(a,b,c) &=-B(a,c,b), \qquad B(a,b,b)=0.
\label{eq:skew}
\end{align}
If, in addition, $c\in H^1_\sigma$, then
\begin{equation}
\label{eq:convection-pairing}
\begin{aligned}
\langle \LP\diver(a\otimes b),c \rangle_{H^{-1},H^1} &= -\int_{\R^3}(a\otimes b):\nabla c\dd x=B(a,b,c).
\end{aligned}
\end{equation}
\end{lemma}

\begin{proof}
The Mikhlin multiplier theorem gives the boundedness of $\LP$ on $L^q$ for every $1<q<\infty$.
The identities $p^2=p=p^{\mathsf T}$ and $\xi^{\mathsf T}p=0$ show that $\LP$ is idempotent, has divergence-free range, and satisfies \eqref{eq:leray-duality}, initially on Schwartz functions and then by density.
Conversely, if $v\in L^q_\sigma$, then $(I-p)\widehat v=0$ away from the origin.
Thus $(I-\LP)v$ has Fourier transform supported at $\{0\}$ and is therefore a polynomial.
Since $(I-\LP)v\in L^q$, this polynomial vanishes.
Hence $\LP v=v$.
This establishes the projection property and \eqref{eq:remove-projection}.
For $q=2$, the symmetry of $p$ gives orthogonality.
The Sobolev boundedness and the commutation properties follow from the multiplier representation.
The extension to negative-order spaces is defined by duality: $$ \langle\LP F,\varphi\rangle := \langle F,\LP\varphi\rangle, \qquad F\in W^{-1,q},\quad \varphi\in W^{1,q'}.
$$ This proves \eqref{eq:negative-duality}.
Applying the derivative commutation to each column of $F$ gives \eqref{eq:tensor-projection}.
Identity \eqref{eq:tensor-product} follows from the product rule and \eqref{eq:2.1_tensor-convention}.

For $a\in L^3_\sigma$ and $b,c\in H^1$, H\"older's inequality and the Sobolev embedding $H^1\hookrightarrow L^6$ give \eqref{eq:trilinear-bound}.
If $b,c$ are smooth and compactly supported, the distributional solenoidality of $a$ implies $$ 0 = \int_{\R^3}a\cdot\nabla(b\cdot c)\dd x = B(a,b,c)+B(a,c,b).
$$ Density in $H^1$, together with \eqref{eq:trilinear-bound}, proves \eqref{eq:skew} for general $b,c\in H^1$.

Finally, $a\otimes b\in L^2$.
The definition of distributional divergence and \eqref{eq:negative-duality} yield $$ \langle \LP\diver(a\otimes b),c \rangle_{H^{-1},H^1} = -\int_{\R^3}(a\otimes b):\nabla c\dd x.
$$ Integration by parts, using $\diver a=0$ and \eqref{eq:skew}, identifies the right-hand side with $B(a,b,c)$.
This proves \eqref{eq:convection-pairing}.
\end{proof}
For $q>3$, the two lemmas imply $$ \|\LP(\TN(v)-\TN(w))\|_{q/3} \le C_{q,\nu,N}(\|v\|_q^2+\|w\|_q^2)\|v-w\|_q.
$$ The argument does not extend to $q=3$, since $\LP$ is not bounded on $L^1$.
The critical formulation instead uses $u\in L^5$ to place $\TN(u)$ in $L^{5/3}$.
Likewise, for $u\in L^5_\sigma$, duality gives only $$ \int_{\R^3}\LP\TN(u)\cdot|u|u\dd x =\int_{\R^3}\TN(u)\cdot\LP(|u|u)\dd x.
$$ Since $|u|u$ need not be solenoidal, neither monotonicity nor removal of the Leray projection determines the sign of this $L^3$ energy pairing.

\begin{remark}
[Pressure in the critical mild formulation] The equation is formulated using the Leray projection, so no independent pressure estimate is needed in the global continuation argument.
Nevertheless, pressure can be recovered from the critical solution.

Indeed, for $u\in L^5((0,T)\times\R^3)$, define $$
\begin{aligned}
\pi_{\mathrm{conv}} &:=\sum_{i,j=1}^3 \mathcal R_i\mathcal R_j(u_i u_j), &
\pi_{\mathrm{tame}} &:=(-\Delta)^{-1}\diver\TN(u), \qquad \pi:=\pi_{\mathrm{conv}}+\pi_{\mathrm{tame}},
\end{aligned}
$$ where $\mathcal R_j=\partial_j(-\Delta)^{-1/2}$.
The Calder\'on--Zygmund and Hardy--Littlewood--Sobolev inequalities give $$
\begin{aligned}
\|\pi_{\mathrm{conv}}\|_{L^{5/2}(Q_T)} &\le C\|u\|_{L^5(Q_T)}^2, &
\|\pi_{\mathrm{tame}}\|_ {L^{5/3}(0,T;L^{15/4})} &\le CC_g\|u\|_{L^5(Q_T)}^3.
\end{aligned}
$$ Consequently, $$ -\Delta\pi = \sum_{i,j=1}^3\partial_i\partial_j(u_i u_j) +\diver\TN(u) $$ in space--time distributions, and the projected drift agrees with the corresponding velocity--pressure formulation.

The pressure need not vanish in an $L^3$ energy pairing, since $|u|u$ is generally not solenoidal.
In the present argument, the critical estimates retain the projected forcing, whereas the $L^2$ and $H^1$ estimates are obtained by testing the remainder equation against solenoidal fields.
This is why no separate pressure estimate enters the global continuation proof.
\end{remark}

\section{Critical evolution estimates}
\label{Section 3} We establish the evolution estimates needed to construct and compare solutions in the critical space $\mathbb X_T$.
The deterministic estimates control the heat extension of the initial datum and the convolution operators associated with the quadratic convection and cubic taming terms.
For Wiener and compensated Poisson forcing, we derive cubic moment bounds and critical space--time estimates.
We then establish conditional estimates and restart identities on stochastic intervals, which will be used in the local construction, pathwise uniqueness, and maximal continuation in \cref{Section 4}.

Throughout this section, $S(t)=e^{\nu t\Delta}$ with $\nu>0$, $E=L^3_\sigma$, and $\mathbb X_T$ is defined in \eqref{eq:critical-space}.
For $0<T\le\infty$, set $$ Q_T:=(0,T)\times\R^3,\qquad \|f\|_{L^q(Q_T)} :=\left(\int_0^T\|f(t)\|_q^q\dd t\right)^{1/q} \quad(1\le q<\infty).
$$ The first three lemmas are deterministic; the subsequent results treat predictable Wiener--Poisson inputs.
Neither an $L^2$ assumption on the initial datum nor the additional hypothesis \cref{hyp:G} is used.
Unless stated otherwise, constants are independent of $T$, $\nu$, and the functions being estimated and dependence on $\nu$ is displayed explicitly.
All integrals against $\PN$ and $\tPN$ with limits from $0$ to $t$ are understood over $(0,t]\times Z$.
The initial value is prescribed separately.

We begin by placing the heat extension of $L^3_\sigma$ initial data in the critical evolution class and recording its cubic dissipation estimate.
The following lemma also proves that its $L^5$ norm vanishes on shrinking initial intervals, uniformly over compact sets of initial data.
The additional $L^8(0,\infty;L^4)$ and positive-time smoothing estimates provide the caloric-field bounds used in \cref{Section 5} and \cref{Section 6}.
\begin{lemma}
[Heat extension of critical initial data] \label{lem:heat-extension} Let $a\in E$ and set $h(t)=S(t)a$.
Then $$ h\in C([0,\infty);E)\cap L^5(0,\infty;L^5_\sigma) \cap L^8(0,\infty;L^4_\sigma),\qquad h(0)=a, $$ and the following assertions hold.
\begin{enumerate}
[label=\textup{(\roman*)},leftmargin=*] \item For every $t\ge0$,
\begin{equation}
\label{eq:3.1_heat-cubic}
\|h(t)\|_3^3+\frac{8\nu}{3}\int_0^t \bigl\|\nabla(|h(s)|^{3/2})\bigr\|_2^2\dd s \le\|a\|_3^3,
\end{equation}
and
\begin{equation}
\label{eq:3.1_heat-critical}
\sup_{t\ge0}\|h(t)\|_3 +\nu^{1/5}\|h\|_{L^5(Q_\infty)} +\nu^{1/8}\|h\|_{L^8(0,\infty;L^4)} \le C\|a\|_3.
\end{equation}

\item The critical norms vanish on shrinking initial intervals:
\begin{equation}
\label{eq:3.1_heat-short}
\lim_{t\downarrow0}\left( \|h\|_{L^5(Q_t)}+\|h\|_{L^8(0,t;L^4)}\right)=0.
\end{equation}
For every compact set $K\subset E$,
\begin{equation}
\label{eq:3.1_heat-compact}
\lim_{t\downarrow0}\sup_{a\in K} \|S(\cdot)a\|_{L^5(Q_t)}=0,
\end{equation}
and if $a_n\to a$ in $E$, then
\begin{equation}
\label{eq:3.1_heat-sequence}
\lim_{t\downarrow0}\sup_{n\ge1} \|S(\cdot)a_n\|_{L^5(Q_t)}=0.
\end{equation}

\item For every integer $m\ge0$, every $3\le q\le\infty$, and every $t>0$,
\begin{equation}
\label{eq:3.1_heat-smoothing}
\|\nabla^m h(t)\|_q\le C_{m,q} (\nu t)^{-m/2-\frac32(\frac13-\frac1q)}\|a\|_3.
\end{equation}
For every $\delta>0$,
\begin{align}
\sup_{t\ge\delta}\|\nabla h(t)\|_\infty &\le C(\nu\delta)^{-1}\|a\|_3,
\label{eq:3.1_heat-gradient}\\
\sup_{t\ge\delta}\bigl\||h(t)|\,|\nabla h(t)|\bigr\|_2 &\le C(\nu\delta)^{-3/4}\|a\|_3^2.
\label{eq:3.1_heat-mixed}
\end{align}
\end{enumerate}
\end{lemma}

\begin{proof}
Gaussian convolution gives $L^3$ contraction, strong continuity, and preservation of distributional solenoidality.
Hence $h\in C([0,\infty);E)$ and $h(0)=a$.
For $a\in C_c^\infty(\R^3;\R^3)$, not necessarily solenoidal, testing $\partial_t h=\nu\Delta h$ with $|h|h$ yields
\begin{equation}
\frac13\frac{\mathrm d}{\mathrm dt}\|h(t)\|_3^3 +\nu\int_{\R^3}|h| \bigl(|\nabla h|^2+|\nabla|h||^2\bigr)\dd x=0.
\end{equation}
The chain rule and integration by parts are legitimate for such data; the identity extends across zeros of $h$ in the Sobolev sense.
Since $|\nabla|h||\le|\nabla h|$ almost everywhere and $$ |\nabla(|h|^{3/2})|^2=\frac94|h|\,|\nabla|h||^2, $$ we obtain $$ \frac{\mathrm d}{\mathrm dt}\|h(t)\|_3^3 +\frac{8\nu}{3}\|\nabla(|h|^{3/2})\|_2^2\le0.
$$ Integration proves \eqref{eq:3.1_heat-cubic} for smooth data.

For general $a\in L^3$, take $a_n\in C_c^\infty$ with $a_n\to a$ in $L^3$ and write $h_n=S(\cdot)a_n$.
Heat contraction gives $h_n\to h$ in $C([0,T];L^3)$ for every finite $T$.
The elementary power inequality
\begin{equation}
\label{eq:3.1_cubic-power}
\bigl\||f|^{3/2}-|g|^{3/2}\bigr\|_2 \le C\bigl(\|f\|_3^{1/2}+\|g\|_3^{1/2}\bigr) \|f-g\|_3
\end{equation}
then implies $|h_n|^{3/2}\to|h|^{3/2}$ in $C([0,T];L^2)$.
The energy estimate bounds $\nabla(|h_n|^{3/2})$ in $L^2(Q_T)$.
Distributional testing identifies every weak limit of these gradients with $\nabla(|h|^{3/2})$.
Weak lower semicontinuity and strong $L^3$ convergence at each time yield \eqref{eq:3.1_heat-cubic} for $a\in E$.
In particular, $|h|^{3/2}\in L^2(0,T;H^1)$ for each finite $T$. For $v\in L^3$ with $|v|^{3/2}\in H^1$, interpolation between $L^3$ and $L^9$ and the homogeneous Sobolev inequality give
\begin{equation}
\label{eq:3.1_critical}
\|v\|_5^5\le\|v\|_3^2\|v\|_9^3 =\|v\|_3^2\bigl\||v|^{3/2}\bigr\|_6^2 \le C_S^2\|v\|_3^2\bigl\|\nabla(|v|^{3/2})\bigr\|_2^2.
\end{equation}
Applying this at almost every time and using \eqref{eq:3.1_heat-cubic} gives $$ \int_0^\infty\|h(t)\|_5^5\dd t \le C\nu^{-1}\|a\|_3^5.
$$ Moreover, interpolation between $L^3$ and $L^5$ gives $\|h(t)\|_4^8\le\|h(t)\|_3^3\|h(t)\|_5^5$; therefore $$ \int_0^\infty\|h(t)\|_4^8\dd t \le C\nu^{-1}\|a\|_3^8.
$$ Together with $L^3$ contraction, this proves \eqref{eq:3.1_heat-critical}.
Absolute continuity of the integrals proves \eqref{eq:3.1_heat-short}.

To establish \eqref{eq:3.1_heat-compact}, note that linearity and \eqref{eq:3.1_heat-critical} imply $$ \|S(\cdot)(a-b)\|_{L^5(Q_t)} \le C\nu^{-1/5}\|a-b\|_3\qquad(t>0).
$$ A finite $L^3$-net for each compact $K\subset E$, combined with \eqref{eq:3.1_heat-short} at the finitely many net points, gives uniform convergence on $K$.
Since $\{a\}\cup\{a_n:n\ge1\}$ is compact whenever $a_n\to a$ in $E$, \eqref{eq:3.1_heat-sequence} follows.

Finally, Gaussian scaling gives, for $1\le r\le\infty$, $$ \|\nabla^mG_{\nu t}\|_r \le C_{m,r}(\nu t)^{-m/2-\frac32(1-1/r)}.
$$ Choosing $r$ by $1+1/q=1/r+1/3$ and applying Young's convolution inequality proves \eqref{eq:3.1_heat-smoothing}.
The choice $(m,q)=(1,\infty)$ yields \eqref{eq:3.1_heat-gradient}, whereas $$ \bigl\||h(t)|\,|\nabla h(t)|\bigr\|_2 \le\|h(t)\|_6\|\nabla h(t)\|_3 \le C(\nu t)^{-1/4}(\nu t)^{-1/2}\|a\|_3^2 $$ proves \eqref{eq:3.1_heat-mixed}.
\end{proof}

\paragraph{Uniform small-time control}
The compactness assumption in \eqref{eq:3.1_heat-compact} cannot be replaced by boundedness in $E$.
Indeed, choose nonzero $a\in C_c^\infty(\R^3;\R^3)\cap E$ and set $a_\lambda(x)=\lambda a(\lambda x)$.
Then $\|a_\lambda\|_3=\|a\|_3$, while heat scaling gives $$ \|S(\cdot)a_\lambda\|_{L^5(Q_t)}^5 =\|S(\cdot)a\|_{L^5(Q_{\lambda^2t})}^5.
$$ Taking $\lambda=t^{-1/2}$ yields the fixed positive value $\|S(\cdot)a\|_{L^5(Q_1)}^5$, demonstrating the failure of bounded-set uniform convergence.

For $u_0\in L^0(\Omega,\mathcal F_0;E)$, the deterministic heat estimates apply almost surely.
Under \eqref{eq:initial}, $$ \E\|S(\cdot)u_0\|_{L^5(Q_\infty)}^3 \le C\nu^{-3/5}\E\|u_0\|_3^3, \qquad \lim_{t\downarrow0} \E\|S(\cdot)u_0\|_{L^5(Q_t)}^3=0, $$ where the limit follows by dominated convergence.

For $F\in L^{5/2}(Q_T;\R^{3\times3})$, define
\begin{equation}
\label{eq:3.2_tensor}
(\mathcal K_TF)(t) :=\int_0^t S(t-s)\LP\diver F(s)\dd s =\diver\int_0^t S(t-s)\LP^{(1)}F(s)\dd s, \qquad 0\le t\le T,
\end{equation}
initially as an identity in $W^{-1,5/2}$.

\begin{lemma}
[Critical convection estimate and time continuity] \label{lem:3.2_critical-convection} Let $0<T<\infty$ and $F\in L^{5/2}(Q_T;\R^{3\times3})$.
Then $\mathcal K_TF$ has a unique representative in $$ C([0,T];E)\cap L^5(0,T;L^5_\sigma), \qquad (\mathcal K_TF)(0)=0.
$$ There is a constant $C$, independent of $T$ and $\nu$, such that
\begin{equation}
\label{eq:3.2_tdc}
\sup_{0\le t\le T}\|(\mathcal K_TF)(t)\|_3 +\nu^{1/5}\|\mathcal K_TF\|_{L^5(Q_T)} \le C\nu^{-3/5}\|F\|_{L^{5/2}(Q_T)}.
\end{equation}
Moreover,
\begin{equation}
\label{eq:3.2_tde}
\sup_{0\le t\le T}\|\mathcal K_TF(t)\|_3^3 +\nu\int_0^T \bigl\|\nabla(|\mathcal K_TF(t)|^{3/2})\bigr\|_2^2\dd t \le C\nu^{-9/5}\|F\|_{L^{5/2}(Q_T)}^3.
\end{equation}

For $u,v\in L^5(Q_T;\R^3)$, set $\mathcal B_T(u,v):=\mathcal K_T(u\otimes v)$, with the tensor convention \eqref{eq:2.1_tensor-convention}.
This defines a bounded bilinear map into $$ C([0,T];E)\cap L^5(0,T;L^5_\sigma)\subset\mathbb X_T, $$ with
\begin{equation}
\label{eq:3.2_convection}
\|\mathcal B_T(u,v)\|_{\mathbb X_T} \le C_\nu\|u\|_{L^5(Q_T)}\|v\|_{L^5(Q_T)}, \qquad C_\nu=C\bigl(\nu^{-3/5}+\nu^{-4/5}\bigr).
\end{equation}
The diagonal map satisfies
\begin{equation}
\label{eq:3.2_critical-diff}
\|\mathcal B_T(u,u)-\mathcal B_T(v,v)\|_{\mathbb X_T} \le C_\nu\bigl(\|u\|_{L^5(Q_T)}+\|v\|_{L^5(Q_T)}\bigr) \|u-v\|_{L^5(Q_T)}.
\end{equation}
\end{lemma}

\begin{proof}
By \cref{lem:projection}, $\LP^{(1)}$ is bounded on $L^{5/2}(\R^3;\R^{3\times3})$, and $\LP\diver F=\diver(\LP^{(1)}F)$ in $W^{-1,5/2}$.
Since $F\in L^1(0,T;L^{5/2})$, the tensor-valued convolution $$ H_F(t):=\int_0^t S(t-s)\LP^{(1)}F(s)\dd s $$ is well defined in $L^{5/2}$ for every $t\in[0,T]$.
Its divergence defines \eqref{eq:3.2_tensor} in $W^{-1,5/2}$.
Heat contraction and H\"older's inequality give
\begin{equation*}
\sup_{0\le t\le T}\|\mathcal K_TF(t)\|_{W^{-1,5/2}} \le C\int_0^T\|F(s)\|_{5/2}\dd s \le CT^{3/5}\|F\|_{L^{5/2}(Q_T)}.
\end{equation*}

First take $F\in C_c^\infty(Q_T;\R^{3\times3})$ and set $G=\LP^{(1)}F$, $w=\mathcal K_TF$.
Then
\begin{equation}
\label{eq:3.2_forced}
\partial_t w-\nu\Delta w=\diver G,\qquad w(0)=0.
\end{equation}
The $W^{k,3}$-boundedness of $\LP^{(1)}$ and heat-semigroup regularity yield $$ w\in C^1([0,T];L^3)\cap C([0,T];W^{k,3}) \quad\text{for every integer }k\ge0.
$$ Thus testing \eqref{eq:3.2_forced} with $|w|w\in W^{1,3/2}$ is justified by Sobolev duality and gives $$ \frac13\frac{\mathrm d}{\mathrm dt}\|w(t)\|_3^3 +\nu\int_{\R^3}|w|\bigl(|\nabla w|^2+|\nabla|w||^2\bigr)\dd x =-\int_{\R^3}G:\nabla(|w|w)\dd x.
$$ In particular, the projection remains on the forcing tensor; no solenoidality of $|w|w$ is invoked.
Since $$ |\nabla(|w|w)|\le2|w||\nabla w|, \qquad \left|\int_{\R^3}G:\nabla(|w|w)\dd x\right| \le\frac\nu2\int_{\R^3}|w||\nabla w|^2\dd x +\frac2\nu\int_{\R^3}|G|^2|w|\dd x, $$ and $|\nabla|w||\le|\nabla w|$ and $|\nabla(|w|^{3/2})|^2=\frac94|w||\nabla|w||^2$ almost everywhere, we obtain
\begin{equation}
\label{eq:3.2_cubic-inequality}
\frac{\mathrm d}{\mathrm dt}\|w(t)\|_3^3 +2\nu\bigl\|\nabla(|w(t)|^{3/2})\bigr\|_2^2 \le\frac6\nu\int_{\R^3}|G(t)|^2|w(t)|\dd x.
\end{equation}

Put $$
\begin{gathered}
A:=\sup_{0\le t\le T}\|w(t)\|_3^3,\qquad D:=\int_0^T\|\nabla(|w|^{3/2})\|_2^2\dd t, \\
M:=A+\nu D,\qquad X:=\|w\|_{L^5(Q_T)},\qquad f:=\|F\|_{L^{5/2}(Q_T)}.
\end{gathered}
$$ Integrating \eqref{eq:3.2_cubic-inequality}, using $w(0)=0$, and applying H\"older's inequality and boundedness of $\LP^{(1)}$ give
\begin{equation}
M\le\frac9\nu\int_{Q_T}|G|^2|w|\dd x\dd t \le C\nu^{-1}\|G\|_{L^{5/2}(Q_T)}^2 X \le C\nu^{-1}f^2X.
\end{equation}
The critical Sobolev interpolation inequality \eqref{eq:3.1_critical} yields
\begin{equation}
X^5\le CA^{2/3}D\le C\nu^{-1}M^{5/3}, \qquad X\le C\nu^{-1/5}M^{1/3}.
\end{equation}
Consequently, $M\le C\nu^{-6/5}f^2M^{1/3}$.
If $M=0$, the claim is immediate; otherwise, division by $M^{1/3}$ gives $$ M^{1/3}\le C\nu^{-3/5}f, \qquad X\le C\nu^{-4/5}f.
$$ This proves \eqref{eq:3.2_tdc}--\eqref{eq:3.2_tde} for smooth forcing, with constants independent of $T$ and $\nu$ except for the displayed powers of $\nu$.

For general $F\in L^{5/2}(Q_T;\R^{3\times3})$, choose $F_n\in C_c^\infty(Q_T;\R^{3\times3})$ such that $F_n\to F$ in $L^{5/2}(Q_T)$, and set $w_n=\mathcal K_TF_n$.
Applying the preceding bounds to $F_n-F_m$ shows that $(w_n)$ is Cauchy in $C([0,T];L^3)\cap L^5(Q_T)$.
Its limit $w$ belongs to $C([0,T];E)\cap L^5(0,T;L^5_\sigma)$ and satisfies $w(0)=0$.
Furthermore, $$ \sup_{0\le t\le T} \|\mathcal K_T(F_n-F)(t)\|_{W^{-1,5/2}} \le CT^{3/5}\|F_n-F\|_{L^{5/2}(Q_T)}\longrightarrow0.
$$ For every $t\in[0,T]$, the $L^3$ convergence of $w_n(t)$ and the $W^{-1,5/2}$ convergence of the distributional convolutions both imply convergence in $\mathcal S'$.
Uniqueness of distributional limits identifies $w(t)$ with \eqref{eq:3.2_tensor} for every $t$, including the endpoints.
The continuous $E$-valued representative is therefore unique.

The estimate \eqref{eq:3.2_tdc} passes to the limit directly.
By \eqref{eq:3.1_cubic-power}, $|w_n|^{3/2}\to|w|^{3/2}$ in $C([0,T];L^2)$.
The dissipation bound makes $\nabla(|w_n|^{3/2})$ bounded in $L^2(Q_T)$, and its weak limit is $\nabla(|w|^{3/2})$, as follows by distributional testing.
Weak lower semicontinuity proves \eqref{eq:3.2_tde}.

The construction is consistent under restriction to $(0,t)$, so \eqref{eq:3.2_tdc} also gives $$ \sup_{0\le r\le t}\|w(r)\|_3 \le C\nu^{-3/5}\|F\|_{L^{5/2}(Q_t)} \longrightarrow0\qquad(t\downarrow0).
$$ For $0\le s\le t\le T$, splitting the convolution at $s$ and using the semigroup property first yields, in $W^{-1,5/2}$, $$ w(t)=S(t-s)w(s) +\int_s^t S(t-r)\LP\diver F(r)\dd r.
$$ The translated convolution on $[s,t]$ has a continuous $E$-valued representative by the estimate already proved; with this interpretation, the identity holds in $E$.

Finally, the tensor convention \eqref{eq:2.1_tensor-convention} implies $|u\otimes v|=|u||v|$, and H\"older's inequality gives $$ \|u\otimes v\|_{L^{5/2}(Q_T)} \le\|u\|_{L^5(Q_T)}\|v\|_{L^5(Q_T)}.
$$ Combining this bound with \eqref{eq:3.2_tdc} and the definition of $\mathbb X_T$ in \eqref{eq:critical-space} proves \eqref{eq:3.2_convection}.
The bilinear identity $$ \mathcal B_T(u,u)-\mathcal B_T(v,v) =\mathcal B_T(u-v,u)+\mathcal B_T(v,u-v) $$ then gives \eqref{eq:3.2_critical-diff}.
\end{proof}

For $f\in L^{5/3}(Q_T;\R^3)$, define
\begin{equation}
\label{eq:3.3_zero_D}
(\mathcal H_Tf)(t):=\int_0^tS(t-s)\LP f(s)\dd s, \qquad 0\le t\le T,
\end{equation}
initially as an $L^{5/3}$-valued Bochner integral.

Under \cref{hyp:T1}, the cubic growth estimate in \cref{lem:taming} places $\TN(u)$ in $L^{5/3}(Q_T)$ whenever $u\in L^5(Q_T)$.
We therefore estimate the Duhamel operator for vector-valued forcing in this space, including continuity of its $E$-valued representative and its zero initial value.

\begin{lemma}
[Critical cubic taming estimate] \label{lem:critical-taming} Let $0<T<\infty$ and $f\in L^{5/3}(Q_T;\R^3)$.
Then $\mathcal H_Tf$ has a unique representative in $$ C([0,T];E)\cap L^5(0,T;L^5_\sigma), \qquad (\mathcal H_Tf)(0)=0, $$ and
\begin{equation}
\label{eq:3.3_zero-critical}
\sup_{0\le t\le T}\|(\mathcal H_Tf)(t)\|_3 +\nu^{1/5}\|\mathcal H_Tf\|_{L^5(Q_T)} \le C\nu^{-2/5}\|f\|_{L^{5/3}(Q_T)},
\end{equation}
where $C$ is independent of $T$ and $\nu$.

If, in addition, \cref{hyp:T1} holds, then for all $u,v\in L^5(Q_T;\R^3)$, $$ \mathcal H_T(\TN(u))\in C([0,T];E)\cap L^5(0,T;L^5_\sigma), \qquad \mathcal H_T(\TN(u))(0)=0, $$ and
\begin{align}
\|\mathcal H_T(\TN(u))\|_{\mathbb X_T} &\le C_{\nu,g}\|u\|_{L^5(Q_T)}^3,
\label{eq:3.3_tame1}\\
\|\mathcal H_T(\TN(u))-\mathcal H_T(\TN(v))\|_{\mathbb X_T} &\le C_{\nu,g} \bigl(\|u\|_{L^5(Q_T)}^2+\|v\|_{L^5(Q_T)}^2\bigr) \|u-v\|_{L^5(Q_T)},
\label{eq:3.3_tame2}
\end{align}
where one may take $C_{\nu,g}=CC_g(\nu^{-2/5}+\nu^{-3/5})$.
\end{lemma}

\begin{proof}
The $L^{5/3}$-boundedness of $\LP$ and the heat contraction give $$ \sup_{0\le t\le T}\|\mathcal H_Tf(t)\|_{5/3} \le C\int_0^T\|f(s)\|_{5/3}\dd s \le CT^{2/5}\|f\|_{L^{5/3}(Q_T)}.
$$ Thus \eqref{eq:3.3_zero_D} defines a continuous $L^{5/3}$-valued function, by strong continuity of $S(t)$.

First let $f\in C_c^\infty(Q_T;\R^3)$, and set $w=\mathcal H_Tf$ and $g=\LP f$.
Then $w_t-\nu\Delta w=g$ and $w(0)=0$.
Since $\LP$ is bounded on $W^{k,3}$ and commutes with spatial derivatives, standard heat-semigroup regularity justifies testing this equation with $|w|w\in W^{1,3/2}$, as in \cref{lem:3.2_critical-convection}.
This gives $$ \frac13\frac{\mathrm d}{\mathrm dt}\|w(t)\|_3^3 +\nu\int_{\R^3}|w| \bigl(|\nabla w|^2+|\nabla|w||^2\bigr)\dd x =\int_{\R^3}g\cdot|w|w\dd x.
$$ Using $|\nabla|w||\le|\nabla w|$ and $|\nabla(|w|^{3/2})|^2=\frac94|w||\nabla|w||^2$, we obtain
\begin{equation}
\label{eq:3.3_zero_C}
\frac{\mathrm d}{\mathrm dt}\|w(t)\|_3^3 +\frac{8\nu}{3}\|\nabla(|w(t)|^{3/2})\|_2^2 \le 3\int_{\R^3}|\LP f(t)|\,|w(t)|^2\dd x.
\end{equation}
The projection is retained in the forcing term; no solenoidality of $|w|w$ is assumed.

Put $$
\begin{gathered}
A:=\sup_{0\le t\le T}\|w(t)\|_3^3, \qquad D:=\int_0^T\|\nabla(|w|^{3/2})\|_2^2\dd t, \\
M:=A+\nu D,\qquad X:=\|w\|_{L^5(Q_T)}, \qquad b:=\|f\|_{L^{5/3}(Q_T)}.
\end{gathered}
$$ Integrating \eqref{eq:3.3_zero_C} and using $w(0)=0$ yield $$ A\le3I,\qquad \nu D\le\frac98I, \qquad I:=\int_{Q_T}|\LP f|\,|w|^2\dd x\dd t.
$$ By H\"older's inequality with exponents $5/3$ and $5/2$, and the $L^{5/3}$-boundedness of $\LP$, $$ M\le\frac{33}{8}I \le C\|\LP f\|_{L^{5/3}(Q_T)}X^2 \le CbX^2.
$$ The critical Sobolev interpolation estimate \eqref{eq:3.1_critical} gives $$ X^5\le CA^{2/3}D\le C\nu^{-1}M^{5/3}, \qquad X\le C\nu^{-1/5}M^{1/3}.
$$ Consequently $M\le C\nu^{-2/5}bM^{2/3}$.
If $M=0$ the conclusion is immediate; otherwise, division by $M^{2/3}$ gives $$ M^{1/3}\le C\nu^{-2/5}b, \qquad X\le C\nu^{-3/5}b.
$$ This proves \eqref{eq:3.3_zero-critical} for smooth $f$, with constants independent of $T$ and $\nu$ apart from the explicit powers of $\nu$.

For general $f\in L^{5/3}(Q_T;\R^3)$, choose $f_n\in C_c^\infty(Q_T;\R^3)$ converging to $f$ in this space.
By linearity and the preceding bound, $w_n:=\mathcal H_Tf_n$ is Cauchy in $C([0,T];L^3)\cap L^5(Q_T)$.
Let $w$ be its limit.
Since every $w_n$ is solenoidal and $w_n(0)=0$, we have $$ w\in C([0,T];E)\cap L^5(0,T;L^5_\sigma), \qquad w(0)=0.
$$ Moreover, the initial $L^{5/3}$ estimate yields $$ \sup_{0\le t\le T} \|\mathcal H_T(f_n-f)(t)\|_{5/3} \le CT^{2/5}\|f_n-f\|_{L^{5/3}(Q_T)}\longrightarrow0.
$$ At each $t\in[0,T]$, the $L^3$ and $L^{5/3}$ limits agree as distributions.
Hence $w(t)$ coincides with the convolution \eqref{eq:3.3_zero_D} for every $t$, including $T$.
Continuity gives uniqueness of this representative, and the linear estimate passes to the limit.
Consistency under restriction to $(0,t)$ additionally gives $$ \sup_{0\le r\le t}\|\mathcal H_Tf(r)\|_3 \le C\nu^{-2/5}\|f\|_{L^{5/3}(Q_t)} \longrightarrow0\qquad(t\downarrow0).
$$

Finally, under \cref{hyp:T1}, \cref{lem:taming} and space--time H\"older's inequality, using $\frac{1}{5/3}=\frac{1}{5/2}+\frac15$, give $$
\begin{aligned}
\|\TN(u)\|_{L^{5/3}(Q_T)} &\le C_g\|u\|_{L^5(Q_T)}^3, \\
\|\TN(u)-\TN(v)\|_{L^{5/3}(Q_T)} &\le\frac{3C_g}{2} \bigl(\|u\|_{L^5(Q_T)}^2+\|v\|_{L^5(Q_T)}^2\bigr) \|u-v\|_{L^5(Q_T)}.
\end{aligned}
$$ Applying \eqref{eq:3.3_zero-critical} and the definition of $\mathbb X_T$ in \eqref{eq:critical-space} proves \eqref{eq:3.3_tame1}--\eqref{eq:3.3_tame2}, with $C_{\nu,g}=CC_g(\nu^{-2/5}+\nu^{-3/5})$.
\end{proof}

We next consider the linear heat equation with predictable Wiener and compensated Poisson forcing.
The estimates in the results  retains the terminal value and uses quadratic and cubic integrability of the Poisson coefficient, consistently with the coefficient hypothesis \cref{hyp:L}.

\begin{proposition}
[Linear cubic It\^o estimate with Wiener and Poisson forcing] \label{prop:linear-cubic-ito} Let $0<T<\infty$ and $v_0\in L^3(\Omega,\mathcal F_0;E)$.
Let $\Phi$ be a strongly predictable $\gamma(\U;E)$-valued process and let $\Psi:\Omega\times(0,T]\times Z\to E$ be strongly $\mathscr P\otimes\mathcal Z$-measurable.
For a stopping time $0\le\tau\le T$, set
\begin{equation}
\label{eq:3.4_linear}
\begin{aligned}
\mathcal I _\tau(v_0,\Phi,\Psi) := \E\|v_0\|_3^3 & +\E\left(\int_0^\tau \|\Phi(s)\|_{\gamma(\U;E)}^2\dd s\right)^{3/2} \\
&+\E\left(\int_0^\tau\int_Z \|\Psi(s,z)\|_3^2\,\mu(\dd z)\dd s\right)^{3/2} +\E\int_0^\tau\int_Z \|\Psi(s,z)\|_3^3\,\mu(\dd z)\dd s.
\end{aligned}
\end{equation}
Assume $\mathcal I _T(v_0,\Phi,\Psi)<\infty$.
The stochastic convolution
\begin{equation}
\label{eq:3.4_stochastic}
\begin{aligned}
V(t)=S(t)v_0 +\int_0^tS(t-s)\Phi(s)\dd W(s) +\int_0^t\int_Z S(t-s)\Psi(s,z)\, \tPN(\mathrm ds,\mathrm dz)
\end{aligned}
\end{equation}
has an adapted $E$-valued c\`adl\`ag modification, unique up to indistinguishability, with $V(0)=v_0$.
Moreover, $$ |V|^{3/2} \in L^2\bigl(\Omega;L^2(0,T;H^1(\R^3))\bigr), $$ and, for every stopping time $0\le\tau\le T$,
\begin{equation}
\label{eq:3.4_stopped}
\E\sup_{0\le t\le\tau}\|V(t)\|_3^3 +\nu\E\int_0^\tau \bigl\|\nabla(|V(s)|^{3/2})\bigr\|_2^2\dd s \le C\mathcal I _\tau(v_0,\Phi,\Psi).
\end{equation}
The constant $C$ is independent of $T$, $\nu$, $\tau$, and the data.
The same estimate holds for differences of two triples of data and their convolutions.
If $\Psi=0$, the modification has continuous $E$-valued paths.
The supremum in \eqref{eq:3.4_stopped} includes the value at $\tau$.
\end{proposition}

\begin{proof}
Set $\mathcal E(v)=\|v\|_3^3$ on $L^3(\R^3;\R^3)$.
This functional is $C^2$, with
\begin{align*}
D\mathcal E(v)[b] &=3\int_{\R^3}|v|v\cdot b\dd x, \\
D^2\mathcal E(v)[b,c] &=3\int_{\R^3}|v|\,b\cdot c\dd x +3\int_{\R^3}\frac{(v\cdot b)(v\cdot c)}{|v|}\dd x.
\end{align*}
The quotient is zero on $\{v=0\}$.
The pointwise Hessian is Lipschitz, so H\"older's inequality gives continuity of $D^2\mathcal E$ in operator norm and
\begin{equation}
\label{eq:3.4_cubic-d}
\|D\mathcal E(v)\|_{(L^3)^*}\le3\|v\|_3^2, \qquad |D^2\mathcal E(v)[b,c]| \le6\|v\|_3\|b\|_3\|c\|_3.
\end{equation}
For $\mathcal R(v,b)=\mathcal E(v+b)-\mathcal E(v)-D\mathcal E(v)[b]$, convexity and Taylor's formula imply
\begin{equation}
\label{eq:3.4_cubic-j}
\begin{aligned}
0\le\mathcal R(v,b) &=\int_0^1(1-\theta)D^2\mathcal E(v+\theta b)[b,b]\dd\theta \le3\|v\|_3\|b\|_3^2+\|b\|_3^3.
\end{aligned}
\end{equation}

For an orthonormal basis $(e_k)$ of $\U$, define $\mathcal Q(v,\Phi)=\frac12\sum_k D^2\mathcal E(v)[\Phi e_k,\Phi e_k]$.
For finite-coordinate operators, the identity $$ \sum_kD^2\mathcal E(v)[\Phi e_k,\Phi e_k] =\E_\gamma D^2\mathcal E(v) \left[\sum_k\gamma_k\Phi e_k,\sum_k\gamma_k\Phi e_k\right] $$ and \eqref{eq:3.4_cubic-d} give
\begin{equation}
\label{eq:3.4_cubic-w}
0\le\mathcal Q(v,\Phi) \le3\|v\|_3\|\Phi\|_{\gamma(\U;E)}^2.
\end{equation}
Approximation in the $\gamma$-norm extends this definition to every $\Phi\in\gamma(\U;E)$ and proves basis independence.
The ideal property \cite[Proposition 2.3]{MR2330977} also yields
\begin{equation*}
\|D\mathcal E(v)\circ\Phi\|_{\HS(\U;\R)} \le3\|v\|_3^2\|\Phi\|_{\gamma(\U;E)}.
\end{equation*}

For any real locally square-integrable martingale $L$ starting at zero with continuous predictable bracket, Lenglart's inequality \cite[Theorem 1.1]{MR4284625}, applied after localization to $|L|^2$ and $\langle L\rangle$ with exponent $1/2$, gives
\begin{equation}
\label{eq:3.4_predic_b}
\E\sup_{t\le T}|L(t)| \le3\E\langle L\rangle_T^{1/2},
\end{equation}
since the constant supplied by that theorem is $2\sqrt2<3$.
The same inequality holds after stopping at any $\tau\le T$.

Initially choose bounded simple initial data and bounded elementary predictable coefficients with finitely many deterministic spatial values in $E\cap\bigcap_{k\ge0}W^{k,3}$.
Take finitely many Wiener coordinates and $\Psi(s,z)=\sum_{j=1}^J\mathds1_{B_j}(z)\psi_j(s)$, where the $B_j$ are disjoint measurable sets of finite $\mu$-measure.
Integration by parts in the finitely many scalar Wiener convolutions, together with the finite Poisson sum, shows that $V$ has c\`adl\`ag paths in every $W^{k,3}$ and $$ \E\sup_{t\le T}\|V(t)\|_{W^{k,3}}^r<\infty \qquad(k\ge0,\ 1\le r<\infty).
$$ It solves the strong $L^3$ equation $$ \mathrm dV=\nu\Delta V\dd t+\Phi\dd W +\int_Z\Psi(t,z)\,\tPN(\mathrm dt,\mathrm dz), \qquad V(0)=v_0.
$$ Interpret the regularized jump integrals as integrals over the finite mark space indexing the sets $B_j$, with intensity $\mu(B_j)$ at the mark $j$.
Theorem 2.1 of \cite{MR4165650} then applies with $p=3$, three components, and predictable drift $\nu\Delta V_-$, which equals $\nu\Delta V$ almost everywhere in time.
The regularity above verifies its space--time integrability requirements, and the finite mark space meets its countable-generation assumption.
Spatial integration by parts gives, almost surely for every $t\in[0,T]$,
\begin{equation}
\label{eq:3.4_reg_ito}
\begin{aligned}
\mathcal E(V(t)) &+3\nu\int_0^t\int_{\R^3}|V| \bigl(|\nabla V|^2+|\nabla|V||^2\bigr)\dd x\dd s=\mathcal E(v_0)+M(t) +\int_0^t\mathcal Q(V(s-),\Phi(s))\dd s+K(t),
\end{aligned}
\end{equation}
where
\begin{align*}
M(t):=&\int_0^tD\mathcal E(V(s-))\circ\Phi(s)\dd W(s)
 +\int_0^t\int_Z
D\mathcal E(V(s-))[\Psi(s,z)]\,\tPN(\mathrm ds,\mathrm dz), \\
K(t):=&\int_{(0,t]\times Z} \mathcal R(V(s-),\Psi(s,z))\,\PN(\mathrm ds,\mathrm dz).
\end{align*}
The Sobolev chain rule justifies the integrations by parts, including at zeros of $V$.
By \eqref{eq:3.4_cubic-j}, $K$ is increasing, and its integration interval retains the terminal jump.

Orthogonality of the continuous and purely discontinuous martingale parts gives
\begin{equation}
\label{eq:3.4_lin_bracket}
\begin{aligned}
\langle M\rangle_t =&\int_0^t \|D\mathcal E(V(s-))\circ\Phi(s)\|_{\HS(\U;\R)}^2\dd s+\int_0^t\int_Z \bigl|D\mathcal E(V(s-))[\Psi(s,z)]\bigr|^2 \,\mu(\dd z)\dd s \\
\le&9\int_0^t\|V(s-)\|_3^4 \left(\|\Phi(s)\|_{\gamma(\U;E)}^2 +\int_Z\|\Psi(s,z)\|_3^2\,\mu(\dd z)\right)\dd s.
\end{aligned}
\end{equation}
This proves that the predictable bracket is absolutely continuous in time.

Fix a stopping time $\tau\le T$ and set $$
\begin{aligned}
U_\tau:=\sup_{0\le t\le\tau}\|V(t)\|_3, & \qquad A_\tau:=\int_0^\tau \left(\|\Phi(s)\|_{\gamma(\U;E)}^2 +\int_Z\|\Psi(s,z)\|_3^2\,\mu(\dd z)\right)\dd s, \\
 B_\tau&:=\int_0^\tau\int_Z \|\Psi(s,z)\|_3^3\,\mu(\dd z)\dd s.
\end{aligned}
$$ Since $\|V(s-)\|_3\le U_\tau$ for $0<s\le\tau$, \eqref{eq:3.4_lin_bracket}, \eqref{eq:3.4_predic_b}, and Young's inequality yield, for every $\varepsilon>0$,
\begin{equation*}
\E\sup_{t\le\tau}|M(t)| \le C\E[U_\tau^2 A_\tau^{1/2}] \le\varepsilon\E U_\tau^3 +C_\varepsilon\E A_\tau^{3/2}.
\end{equation*}
Similarly, \eqref{eq:3.4_cubic-w} gives
\begin{equation}
\E\int_0^\tau\mathcal Q(V(s-),\Phi(s))\dd s \le3\E[U_\tau A_\tau] \le\varepsilon\E U_\tau^3 +C_\varepsilon\E A_\tau^{3/2}.
\end{equation}
The process $\mathds1_{\{0<s\le\tau\}}$ is predictable.
Compensation for the nonnegative integrand, followed by \eqref{eq:3.4_cubic-j}, therefore gives
\begin{equation*}
\begin{aligned}
\E K(\tau) &=\E\int_0^\tau\int_Z \mathcal R(V(s-),\Psi(s,z))\,\mu(\dd z)\dd s \le3\E[U_\tau A_\tau]+\E B_\tau \le\varepsilon\E U_\tau^3 +C_\varepsilon\E A_\tau^{3/2} +\E B_\tau.
\end{aligned}
\end{equation*}
Thus the estimate involves only the stated quadratic and cubic integrability terms associated with the jump coefficient $\Psi$. As in \cref{lem:heat-extension}, the dissipation in \eqref{eq:3.4_reg_ito} dominates $\frac{8\nu}{3}\int_0^t\|\nabla(|V|^{3/2})\|_2^2\dd s$.
Taking the supremum for the energy and evaluating the dissipation at $\tau$, we obtain
\begin{equation*}
\begin{aligned}
&\E U_\tau^3 +\nu\E\int_0^\tau\|\nabla(|V|^{3/2})\|_2^2\dd s \le C\Bigl( \E\|v_0\|_3^3+\E\sup_{t\le\tau}|M(t)| +\E\int_0^\tau\mathcal Q(V(s-),\Phi(s))\dd s +\E K(\tau)\Bigr).
\end{aligned}
\end{equation*}
Insert the preceding bounds and choose $\varepsilon$ small enough to absorb the terms containing $\E U_\tau^3$.
The inequality $(a+b)^{3/2}\le\sqrt2(a^{3/2}+b^{3/2})$ then proves \eqref{eq:3.4_stopped} for regular data, with $C$ independent of $T$, $\nu$, $\tau$, and the approximation.

Use the predictable subspaces of
\begin{equation*}
\begin{aligned}
\mathcal W_T &:=L^3\bigl(\Omega;L^2(0,T;\gamma(\U;E))\bigr), & \mathcal J_T &:=L^3\bigl(\Omega; L^2((0,T)\times Z,\dd t\,\mu;E)\bigr) \cap L^3\bigl(\Omega\times(0,T)\times Z;E\bigr),
\end{aligned}
\end{equation*}
with the sum norm on the intersection.
Truncation, elementary predictable approximation, finite Wiener-coordinate projection, and heat smoothing give regular triples converging in $L^3(\Omega,\mathcal F_0;E)\times\mathcal W_T\times\mathcal J_T$.
For the jump term, first truncate to bounded values and finite-measure mark sets.
On these sets, H\"older's inequality controls the mixed $L^3(\Omega;L^2)$ norm by the space--time $L^3$ norm, so elementary approximation can be performed simultaneously in both norms.
Heat smoothing preserves solenoidality and converges strongly on $E$.
Therefore,
\begin{equation}
\label{eq:3.4_lin_app}
\mathcal I _T(v_0^n-v_0,\Phi_n-\Phi,\Psi_n-\Psi) \longrightarrow0.
\end{equation}
For the corresponding regular convolutions, linearity and the estimate at $\tau=T$ imply
\begin{equation}
\label{eq:3.4_lin_cauchy}
\E\sup_{t\le T}\|V_n(t)-V_m(t)\|_3^3 \le C\mathcal I _T(v_0^n-v_0^m,\Phi_n-\Phi_m,\Psi_n-\Psi_m) \longrightarrow0.
\end{equation}
A subsequence converges almost surely uniformly on $[0,T]$ to an adapted c\`adl\`ag $E$-valued process $V$ with $V(0)=v_0$.
Fatou's lemma and \eqref{eq:3.4_lin_cauchy} give convergence of the whole sequence such that,
\begin{equation}
\label{eq:3.4_lin_uniform}
\E\sup_{t\le T}\|V_n(t)-V(t)\|_3^3\longrightarrow0.
\end{equation}
These suprema are measurable through rational times and $T$.
On $[0,\tau]$, use rational times strictly before $\tau$ together with $\tau$ and $0$.

For each deterministic $t$, the defining stochastic integrals in \eqref{eq:3.4_stochastic} exist in $L^3(\Omega;E)$.
Indeed, $E$ is UMD and of type $2$, so \cite[Corollary 3.10(1)]{MR2330977}, the $\gamma$-ideal property, and heat contraction control the Wiener term.
For the Poisson term, apply \cite[Theorem 1.1, $p=q=3$]{MR3265175} componentwise.
Minkowski's inequality and heat contraction give, pathwise,
\begin{equation*}
\begin{aligned}
&\left\|\left(\int_0^t\int_Z |S(t-s)\Psi(s,z)|^2\,\mu(\dd z)\dd s\right)^{1/2}\right\|_3 \\
&\le \left(\int_0^t\int_Z \|S(t-s)\Psi(s,z)\|_3^2\,\mu(\dd z)\dd s\right)^{1/2}\le \left(\int_0^t\int_Z \|\Psi(s,z)\|_3^2\,\mu(\dd z)\dd s\right)^{1/2},
\end{aligned}
\end{equation*}
and heat contraction also controls the cubic Poisson integrability term.
Applying these estimates to differences of the initial data and stochastic coefficients shows that the regular convolutions at $t$ converge in $L^3(\Omega;E)$ to \eqref{eq:3.4_stochastic}.
Together with \eqref{eq:3.4_lin_uniform}, this identifies $V$ as a modification of that convolution.
Two c\`adl\`ag modifications agree simultaneously at all rational times and at $T$, and hence everywhere by right continuity; uniqueness up to indistinguishability follows.

\cref{eq:3.1_cubic-power}, H\"older's inequality, and \eqref{eq:3.4_lin_uniform} imply $$ \E\sup_{t\le T} \bigl\||V_n(t)|^{3/2}-|V(t)|^{3/2}\bigr\|_2^2 \longrightarrow0.
$$ The regular estimates bound $\nabla(|V_n|^{3/2})$ in $L^2(\Omega\times Q_T)$.
Testing against products of bounded measurable functions of $(\omega,t)$ and smooth compactly supported spatial functions identifies every weak subsequential limit as $\nabla(|V|^{3/2})$.
This is the lower-semicontinuity argument of \cite[Lemma 3.2]{MR4908978}, applied to the scalar moduli; the displayed convergence and joint measurability suffice.

For a fixed stopping time $\tau\le T$, multiplication by $\mathds1_{\{s\le\tau\}}$ is a bounded operator on $L^2(\Omega\times Q_T)$ and preserves weak convergence.
Hence the stopped dissipation is lower semicontinuous.
The stopped cubic suprema converge in $L^1(\Omega)$ by \eqref{eq:3.4_lin_uniform}, and the corresponding terms in $\mathcal I_\tau$ converge by \eqref{eq:3.4_lin_app}, since restriction to $(0,\tau]$ is contractive in the Wiener and Poisson integrand norms.
Passing to the limit proves \eqref{eq:3.4_stopped} for the original data, including the value at $\tau$.
Moreover, $$ \E\int_0^T\bigl\||V(s)|^{3/2}\bigr\|_2^2\dd s \le T\E\sup_{s\le T}\|V(s)\|_3^3<\infty, $$ which completes the asserted $L^2(\Omega;L^2(0,T;H^1))$ regularity.
Applying the result to the difference of two triples proves the difference assertion.
If $\Psi=0$, choose $\Psi_n=0$.
The regular convolutions are continuous, and so is their almost sure uniform limit.
\end{proof}
\paragraph{Localization and compatibility with \cref{hyp:L}}
The proposition also yields a localized construction when the three noise integrals  in \eqref{eq:3.4_linear} are finite almost surely, without assuming the finiteness of their expectations.
Define
\begin{equation*}
\begin{aligned}
\rho_n:=\inf\Bigg\{t\in[0,T]: \int_0^t\Bigg( &\|\Phi(s)\|_{\gamma(\U;E)}^2 +\int_Z \bigl(\|\Psi(s,z)\|_3^2+\|\Psi(s,z)\|_3^3\bigr) \,\mu(\dd z) \Bigg)\dd s\ge n \Bigg\}\wedge T,
\end{aligned}
\end{equation*}
with $\inf\varnothing=\infty$.
The cumulative integral is continuous, and hence the stopped coefficients $$ \mathds1_{\{0<s\le\rho_n\}}\Phi(s), \qquad \mathds1_{\{0<s\le\rho_n\}}\Psi(s,z) $$ satisfy the integrability conditions of the proposition, with finite corresponding moments.
Their stochastic convolutions agree up to the smaller stopping time by the stopped difference estimate.
Since the cumulative integral is finite at $T$, almost every sample path satisfies $\rho_n=T$ for all sufficiently large $n$.
The convolutions therefore paste to an adapted c\`adl\`ag process on the closed interval $[0,T]$, and the estimate remains valid whenever its right-hand side is finite.
The stopping indicators include the terminal time throughout.
For $v_0\in L^0(\Omega,\mathcal F_0;E)$, apply this construction to the stochastic convolution with zero initial value and add $S(t)v_0$ pathwise.
The heat estimates of \cref{lem:heat-extension} apply to each realization of $v_0$.
The expectation estimate \eqref{eq:3.4_stopped} is finite whenever $\mathcal I_\tau(v_0,\Phi,\Psi)<\infty$.

For example, let $u,v$ be adapted c\`adl\`ag $E$-valued processes with finite cubic moments of their uniform norms, and define $$ \Phi_{u,v}(s) :=\Sigma(s,u(s-))-\Sigma(s,v(s-)), \qquad\Psi_{u,v}(s,z) :=\mathcal G(s,u(s-),z) -\mathcal G(s,v(s-),z).
$$ Hypothesis \cref{hyp:L} gives
\begin{equation*}
\mathcal I_T(0,\Phi_{u,v},\Psi_{u,v}) \le CL_T^3(T^{3/2}+T) \E\sup_{s\le T}\|u(s)-v(s)\|_3^3.
\end{equation*}
Indeed, the two quadratic integrability terms contribute $T^{3/2}$, whereas the cubic Poisson term contributes $T$.
The growth bound \eqref{eq:local-growth} similarly controls the corresponding noise integrals  for a single process.
Thus the estimate requires precisely the Wiener and Poisson integrability prescribed by \cref{hyp:L}, with the coefficients evaluated predictably at the left limits of the processes.

For stopping times $0\le\sigma\le\tau\le T$, we use the pathwise notation
\begin{equation*}
\|v\|_{\mathbb X[\sigma,\tau]} :=\sup_{\sigma\le t\le\tau}\|v(t)\|_3 +\left(\int_\sigma^\tau\|v(t)\|_5^5\dd t\right)^{1/5}.
\end{equation*}
Both endpoints are included in the supremum.
When $\sigma=\tau$, the integral vanishes.
This notation concerns only the restriction to the indicated interval and imposes no condition on a frozen extension after $\tau$.

For predictable stochastic integrands $\Phi$ and $\Psi$, define the nonnegative random variable
\begin{equation*}
\begin{aligned}
\mathcal A_{\sigma,\tau}(\Phi,\Psi) :=& \left( \int_\sigma^\tau \|\Phi(s)\|_{\gamma(\U;E)}^2\dd s \right)^{3/2}+\left( \int_\sigma^\tau\int_Z \|\Psi(s,z)\|_3^2 \,\mu(\dd z)\dd s \right)^{3/2} +\int_\sigma^\tau\int_Z \|\Psi(s,z)\|_3^3 \,\mu(\dd z)\dd s.
\end{aligned}
\end{equation*}
Therefore, $$ \mathcal I_\tau(v_0,\Phi,\Psi) =\E\|v_0\|_3^3 +\E\mathcal A_{0,\tau}(\Phi,\Psi).
$$

\begin{theorem}
[Critical Wiener--Poisson Stokes estimate] \label{thm:critical-stokes} Under the hypotheses of \cref{prop:linear-cubic-ito}, let $V$ denote the c\`adl\`ag modification of \eqref{eq:3.4_stochastic}.
Then $V\in\mathbb X_T$ almost surely.
For every stopping time $0\le\tau\le T$,
\begin{equation}
\label{eq:3.5_crit_stokes}
\E\left[ \sup_{0\le t\le\tau}\|V(t)\|_3^3 +\nu^{3/5} \left(\int_0^\tau\|V(t)\|_5^5\dd t\right)^{3/5} \right] \le C\mathcal I_\tau(v_0,\Phi,\Psi).
\end{equation}
Consequently,
\begin{equation}
\label{eq:3.5_crit-stokes2}
\E\|V\|_{\mathbb X[0,\tau]}^3 \le C_\nu\mathcal I_\tau(v_0,\Phi,\Psi), \qquad C_\nu=C(1+\nu^{-3/5}).
\end{equation}
The constant $C$ is independent of $T$, $\tau$, $\nu$, the initial datum, and the stochastic integrands.

For two triples $(v_{0,i},\Phi_i,\Psi_i)$, $i=1,2$, and their corresponding stochastic convolutions,
\begin{equation}
\label{eq:3.5_crit3}
\E\|V_1-V_2\|_{\mathbb X[0,\tau]}^3 \le C_\nu\mathcal I_\tau (v_{0,1}-v_{0,2}, \Phi_1-\Phi_2,\Psi_1-\Psi_2).
\end{equation}
The estimates remain valid after replacing the stochastic integrands $\Phi,\Psi$ by $\theta\Phi,\theta\Psi$, where $\theta$ is any predictable scalar process with $|\theta|\le1$, with the corresponding quantity $\mathcal I_\tau(v_0,\theta\Phi,\theta\Psi)$ on the right-hand side.

Finally, for fixed $v_0,\Phi,\Psi$ satisfying the hypotheses,
\begin{equation}
\label{eq:3.5_crit_small}
\lim_{t\downarrow0} \E\|V-S(\cdot)v_0\|_{\mathbb X[0,t]}^3=0.
\end{equation}
\end{theorem}

\begin{proof}
Fix a stopping time $\tau\le T$.
On the full-probability set provided by \cref{prop:linear-cubic-ito}, write $$ A_\tau:= \sup_{0\le t\le\tau}\|V(t)\|_3^3, \qquad D_\tau:= \int_0^\tau \bigl\|\nabla(|V(t)|^{3/2})\bigr\|_2^2\dd t.
$$ For almost every time, $|V(t)|^{3/2}\in H^1$.
The interpolation inequality \eqref{eq:3.1_critical} therefore gives
\begin{equation}
\label{eq:3.5_stokes_int}
\int_0^\tau\|V(t)\|_5^5\dd t \le C A_\tau^{2/3}D_\tau.
\end{equation}
In particular, $V\in L^5(0,T;L^5_\sigma)$ almost surely.
The $L^5$ realization is strongly measurable.
Together with the c\`adl\`ag property already proved, this shows $V\in\mathbb X_T$.

Raising \eqref{eq:3.5_stokes_int} to the power $3/5$ and using the weighted arithmetic--geometric mean inequality yields
\begin{equation*}
\begin{aligned}
\nu^{3/5} \left(\int_0^\tau\|V(t)\|_5^5\dd t\right)^{3/5} &\le C A_\tau^{2/5}(\nu D_\tau)^{3/5}\le C\left( \frac25A_\tau+\frac35\nu D_\tau \right).
\end{aligned}
\end{equation*}
Taking expectations and applying \eqref{eq:3.4_stopped} proves \eqref{eq:3.5_crit_stokes}.
The inequality $(a+b)^3\le4(a^3+b^3)$ then proves \eqref{eq:3.5_crit-stokes2}.

Linearity of the stochastic convolution gives \eqref{eq:3.5_crit3}.
Predictable multiplication by $\theta$ preserves the integrability classes of the stochastic integrands, so the same proof applies to the modified integrands.
In particular, predictable stopping indicators are covered without changing the constant.

The process $V-S(\cdot)v_0$ is the stochastic convolution with zero initial datum and the same Wiener and Poisson integrands $\Phi,\Psi$.
Apply \eqref{eq:3.5_crit-stokes2} on $[0,t]$.
Each term in $\mathcal A_{0,t}(\Phi,\Psi)$ tends to zero almost surely by absolute continuity of the relevant time integrals and is dominated by the corresponding integrable term at $T$.
Dominated convergence proves \eqref{eq:3.5_crit_small}.
\end{proof}

For stopping times $0\le\sigma\le\tau\le T$, set $$ \chi_{\sigma,\tau}(s) :=\mathds1_{\{\sigma<s\le\tau\}}, $$ and define, for deterministic $t\in[0,T]$, $$ \mathcal Z_{\sigma,\tau}(t) :=\int_0^t \chi_{\sigma,\tau}(s)S(t-s)\Phi(s)\dd W(s)+\int_0^t\int_Z \chi_{\sigma,\tau}(s)S(t-s)\Psi(s,z) \,\tPN(\mathrm ds,\mathrm dz). $$

We use its c\`adl\`ag modification.
For $\xi\in L^3(\Omega,\mathcal F_\sigma;E)$, write $$ \mathcal V_{\sigma,\tau}^{\xi}(t) :=S(t-\sigma)\xi+\mathcal Z_{\sigma,\tau}(t), \qquad \sigma\le t\le\tau.
$$

Whenever $F$ and $f$ are progressively measurable tensor- and vector-valued processes satisfying $$ F\in L^{5/2}(Q_T;\R^{3\times3}), \qquad f\in L^{5/3}(Q_T;\R^3) \quad\text{almost surely}, $$ define
\begin{equation}
Y(t):=S(t)v_0+\mathcal K_TF(t) +\mathcal H_Tf(t)+\mathcal Z_{0,T}(t),
\end{equation}
$$ Y^\tau(t):=S(t)v_0 +\mathcal K_T(\chi_{0,\tau}F)(t) +\mathcal H_T(\chi_{0,\tau}f)(t) +\mathcal Z_{0,\tau}(t).
$$

To use these estimates on local solution intervals, we require a restart formulation at stopping times.
The following lemma provides conditional critical estimates and restart identities that use the post-jump state at the lower endpoint and retain a possible jump at the upper endpoint.
The lemma also identifies the heat continuation after the deterministic forcing and stochastic integrands have been set to zero, as required by the closed mild formulation in \cref{def:critical-solution}.

\begin{lemma}
[Restart on stochastic intervals and terminal jumps] \label{lem:stopping-time-restart} Assume the hypotheses of \cref{thm:critical-stokes} and let $0\le\sigma\le\tau\le T$ be stopping times.
Then the following assertions hold.

\begin{enumerate}
[label=\textup{(\roman*)},leftmargin=*]

\item The process $\chi_{\sigma,\tau}$ is predictable, and $\mathcal Z_{\sigma,\tau}\in\mathbb X_T$ almost surely, with $ \mathcal Z_{\sigma,\tau}(t)=0, \, 0\le t\le\sigma, $ outside one null set.
Moreover,
\begin{equation}
\label{eq:3.6_post_stop}
\E\|\mathcal Z_{\sigma,\tau}\|_{\mathbb X_T}^3 \le C_\nu \E\mathcal A_{\sigma,\tau}(\Phi,\Psi).
\end{equation}

\item For every $\xi\in L^3(\Omega,\mathcal F_\sigma;E)$, the process $\mathcal V_{\sigma,\tau}^{\xi}$ is adapted and $E$-valued c\`adl\`ag on $[\sigma,\tau]$, with $\mathcal V_{\sigma,\tau}^{\xi}(\sigma)=\xi$.
Furthermore, almost surely,
\begin{equation}
\label{eq:3.6_con_rest}
\begin{aligned}
&\E\left[ \sup_{\sigma\le t\le\tau} \|\mathcal V_{\sigma,\tau}^{\xi}(t)\|_3^3 +\nu^{3/5} \left( \int_\sigma^\tau \|\mathcal V_{\sigma,\tau}^{\xi}(t)\|_5^5 \dd t \right)^{3/5} \,\middle|\,\mathcal F_\sigma \right] \le C\left( \|\xi\|_3^3+ \E\left[ \mathcal A_{\sigma,\tau}(\Phi,\Psi) \,\middle|\,\mathcal F_\sigma \right]\right).
\end{aligned}
\end{equation}
For two triples $(\xi_i,\Phi_i,\Psi_i)$, $i=1,2$, on the same stochastic basis with the same driving noises, the same inequality holds for the difference of their restarted processes, with $(\xi,\Phi,\Psi)$ replaced on the right by $(\xi_1-\xi_2,\Phi_1-\Phi_2,\Psi_1-\Psi_2)$.
\item For $F$ and $f$ satisfying the above integrability assumptions, outside one null set, simultaneously for all $\sigma\le t\le\tau$,
\begin{equation}
\label{eq:3.6_full_stop}
\begin{aligned}
Y(t)=&S(t-\sigma)Y(\sigma) +\mathcal K_T(\chi_{\sigma,\tau}F)(t)+\mathcal H_T(\chi_{\sigma,\tau}f)(t) +\mathcal Z_{\sigma,\tau}(t).
\end{aligned}
\end{equation}
In particular, for the process $V$ of \cref{thm:critical-stokes}
\begin{equation}
\label{eq:3.6_stoch_stop}
V(t)=S(t-\sigma)V(\sigma) +\mathcal Z_{\sigma,\tau}(t), \qquad \sigma\le t\le\tau.
\end{equation}

\item Outside one null set, simultaneously for all $0<t\le T$,
\begin{align}
\tilde\Delta Y(t) &=\int_Z\Psi(t,z)\,\PN(\{t\},\mathrm dz),
\label{eq:3.6_mild_Jump}\\
\tilde\Delta\mathcal Z_{\sigma,\tau}(t) &=\chi_{\sigma,\tau}(t) \int_Z\Psi(t,z)\,\PN(\{t\},\mathrm dz).
\label{eq:3.6_rest_jump}
\end{align}
The integrals on the right denote the jumps of the corresponding $E$-valued compensated Poisson integrals.
Moreover,
\begin{equation}
\label{eq:3.6_heat_continuation}
Y^\tau(t)=
\begin{cases}
Y(t),&0\le t\le\tau, \\
S(t-\tau)Y(\tau),&\tau\le t\le T,
\end{cases}
\end{equation}
outside one null set.

\end{enumerate}
All constants are independent of the stopping times and the deterministic horizon.
In particular, the restart datum is the post-jump value $Y(\sigma)$, while the subsequent stochastic increment excludes $\sigma$ and includes $\tau$.
\end{lemma}

\begin{proof}
The indicator $\chi_{\sigma,\tau}$ is adapted and left-continuous on $(0,T]$, hence predictable.
\cref{thm:critical-stokes} applied to the restricted stochastic integrands gives \eqref{eq:3.6_post_stop} and its stopped estimate at $\sigma$ gives $\E\sup_{t\le\sigma}\|\mathcal Z_{\sigma,\tau}(t)\|_3^3=0$.
This proves \textup{(i)}, including simultaneous vanishing.
All convolutions are formed at deterministic terminal times before evaluating their c\`adl\`ag versions at stopping times.

For deterministic $r$, the semigroup property and bounded-operator interchange with stochastic integrals give
\begin{equation*}
V(t)=S(t-r)V(r)+\mathcal Z_{r,T}(t), \qquad r\le t\le T.
\end{equation*}
Equality first holds at each deterministic time and then simultaneously by c\`adl\`ag paths.
Set $\sigma_n=(2^{-n}\lceil2^n\sigma\rceil)\wedge T\downarrow\sigma$.
For each possible value $r_j$, the event $\{\sigma_n=r_j\}$ belongs to $\mathcal F_{r_j}$ and may be taken through integrals over $(r_j,t]$.
Summing these finitely many identities gives
\begin{equation*}
V(t)=S(t-\sigma_n)V(\sigma_n) +\mathcal Z_{\sigma_n,T}(t), \qquad \sigma_n\le t\le T.
\end{equation*}
Absolute continuity and dominated convergence give $\E\mathcal A_{\sigma,\sigma_n}(\Phi,\Psi)\to0$; hence
\begin{equation}
\label{eq:3.6_cut_con}
\E\|\mathcal Z_{\sigma_n,T} -\mathcal Z_{\sigma,T}\|_{\mathbb X_T}^3 \longrightarrow0.
\end{equation}
Choose a subsequence converging almost surely uniformly in $E$.
Right continuity and heat contraction imply, uniformly for $t\ge\sigma_n$, $$
\begin{aligned}
\|S(t-\sigma_n)V(\sigma_n)-S(t-\sigma)V(\sigma)\|_3 &\le\|V(\sigma_n)-V(\sigma)\|_3+\|(I-S(\sigma_n-\sigma))V(\sigma)\|_3\longrightarrow0.
\end{aligned}
$$ Thus the restart identity holds for all $t>\sigma$ on one event of probability one, and also at $\sigma$ by \textup{(i)}.
The stopped difference estimate at $\tau$ identifies $\mathcal Z_{\sigma,T}$ with $\mathcal Z_{\sigma,\tau}$ on $[0,\tau]$, proving \eqref{eq:3.6_stoch_stop} with both endpoints included.

On each sample path satisfying the assumed forcing integrability, split the Bochner convolutions at an arbitrary $r$ and use the semigroup property.
The identities are first established in $W^{-1,5/2}$ for $\mathcal K_T$ and in $L^{5/3}$ for $\mathcal H_T$.
\cref{lem:3.2_critical-convection} and \cref{lem:critical-taming} identify all terms with their continuous $E$-valued representatives for every $0\le r\le t\le T$.
Taking $r=\sigma$ and restricting to $t\le\tau$ proves \eqref{eq:3.6_full_stop}.
Applying the same lemmas to the forcing terms $F,f$ extended by zero outside $(\sigma,\tau]$ gives
\begin{align}
&\sup_{\sigma\le t\le\tau} \|\mathcal K_T(\chi_{\sigma,\tau}F)(t)\|_3 +\nu^{1/5} \|\mathcal K_T(\chi_{\sigma,\tau}F)\|_ {L^5((\sigma,\tau)\times\R^3)} \le C\nu^{-3/5}\|F\|_{L^{5/2}((\sigma,\tau)\times\R^3)},
\label{eq:3.6_convection_forcing}
\\
&\sup_{\sigma\le t\le\tau} \|\mathcal H_T(\chi_{\sigma,\tau}f)(t)\|_3 +\nu^{1/5} \|\mathcal H_T(\chi_{\sigma,\tau}f)\|_ {L^5((\sigma,\tau)\times\R^3)} \le C\nu^{-2/5}\|f\|_{L^{5/3}((\sigma,\tau)\times\R^3)}.
\label{eq:3.6_cubic_forcing}
\end{align}
These estimates are pathwise and their constants are independent of the interval length.

For $A\in\mathcal F_\sigma$, the process $\mathds1_A\chi_{\sigma,\tau}$ is predictable.
Let $\mathcal Z^A$ be the stochastic convolution with integrands $\mathds1_A\chi_{\sigma,\tau}\Phi$ and $\mathds1_A\chi_{\sigma,\tau}\Psi$.
Then
\begin{equation}
\label{eq:3.6_restart_locality}
\mathcal Z^A(t)=\mathds1_A\mathcal Z_{\sigma,\tau}(t) \qquad(0\le t\le T)
\end{equation}
Indeed, first replace the lower cutoff by $\sigma_n$ and use $A\cap\{\sigma_n=r_j\}\in\mathcal F_{r_j}$ in the preceding finite decomposition, keeping the upper cutoff inside each integrand.
The corresponding differences of the stochastic integrands tend to zero as in \eqref{eq:3.6_cut_con}.
A uniformly convergent subsequence proves \eqref{eq:3.6_restart_locality}.
Applying \eqref{eq:3.5_crit_stokes} to $\mathcal Z^A$ and restricting its output gives $$
\begin{aligned}
\E\left[\mathds1_A\left\{\sup_{\sigma\le t\le\tau}\|\mathcal Z_{\sigma,\tau}(t)\|_3^3 +\nu^{3/5}\left(\int_\sigma^\tau\|\mathcal Z_{\sigma,\tau}(t)\|_5^5\dd t\right)^{3/5}\right\}\right] \le C\E[\mathds1_A\mathcal A_{\sigma,\tau}(\Phi,\Psi)].
\end{aligned}
$$ Testing all $A\in\mathcal F_\sigma$ yields the conditional inequality for $\mathcal Z_{\sigma,\tau}$.
\cref{lem:heat-extension}, applied pathwise after time translation, gives $$ \sup_{\sigma\le t\le\tau}\|S(t-\sigma)\xi\|_3^3 +\nu^{3/5}\left(\int_\sigma^\tau\|S(t-\sigma)\xi\|_5^5\dd t\right)^{3/5} \le C\|\xi\|_3^3.
$$ The triangle inequality in $L^5$ and $(a+b)^3\le4(a^3+b^3)$ prove \eqref{eq:3.6_con_rest}.
Adaptedness follows because, on $\{\sigma\le t\}$, both $\sigma$ and $\xi$ are $\mathcal F_t$-measurable.
Strong continuity of $S$ gives the initial value $\xi$, and previous arguments gives $\mathcal Z_{\sigma,\tau}(\sigma)=0$.
Linearity proves the difference assertion.

Passing to the regular approximations of \cref{prop:linear-cubic-ito} gives, for every $\varphi\in C_c^\infty(\R^3;\R^3)$, $$
\begin{aligned}
\langle\mathcal Z_{0,T}(t),\varphi\rangle =&\nu\int_0^t\langle\mathcal Z_{0,T}(s),\Delta\varphi\rangle\dd s +\int_0^t\langle\Phi(s)\dd W(s),\varphi\rangle+\int_0^t\int_Z\langle\Psi(s,z),\varphi\rangle\,\tPN(\mathrm ds,\mathrm dz).
\end{aligned}
$$ Uniform $L^3$ convergence passes the drift to the limit, and the scalar martingales converge uniformly in probability by their quadratic estimates.
The drift and Wiener term are continuous, so the compensated Poisson jump formula identifies every scalar jump.
A countable family of smooth test functions dense in $L^{3/2}$ separates the jumps of the $E$-valued c\`adl\`ag processes.
The heat and deterministic Duhamel terms are continuous, hence \eqref{eq:3.6_mild_Jump} follows, and applying the same argument to the restricted stochastic integrands gives \eqref{eq:3.6_rest_jump}.

The deterministic convolutions and the stopped stochastic difference estimate give $Y^\tau=Y$ simultaneously on $[0,\tau]$.
Restart $Y^\tau$ at $\tau$ using \eqref{eq:3.6_full_stop}.
All subsequent forcing terms and stochastic integrands vanish, so $Y^\tau(t)=S(t-\tau)Y(\tau)$ for $\tau\le t\le T$.
This proves \eqref{eq:3.6_heat_continuation} and the post-jump restart convention.
\end{proof}

\paragraph{Short random intervals and the nonlinear mild formula}
Suppose $\tau-\sigma\le h$ almost surely, where $h>0$ is deterministic.
Let $u,v$ be adapted c\`adl\`ag $E$-valued processes with finite cubic moments of their uniform norms, and use the coefficient differences $\Phi_{u,v},\Psi_{u,v}$ defined above.
Hypothesis \cref{hyp:L} gives, pathwise, $$ \mathcal A_{\sigma,\tau}(\Phi_{u,v},\Psi_{u,v}) \le CL_T^3(h^{3/2}+h) \sup_{\sigma\le s\le\tau}\|u(s)-v(s)\|_3^3.
$$ For $\sigma<s\le\tau$, the pre-jump difference is bounded by the supremum on the right, which justifies this estimate also when either endpoint is a jump time.
Writing $\mathcal Z_{\sigma,\tau}[\Phi,\Psi]$ to indicate the stochastic integrands, the conditional estimate consequently yields
\begin{equation}
\label{eq:3.6_random_interval_noise_lipschitz}
\begin{aligned}
&\E\left[ \|\mathcal Z_{\sigma,\tau}[\Phi_{u,v},\Psi_{u,v}] \|_{\mathbb X[\sigma,\tau]}^3 \,\middle|\,\mathcal F_\sigma\right] \le C_\nu L_T^3(h^{3/2}+h) \E\left[ \sup_{\sigma\le s\le\tau}\|u(s)-v(s)\|_3^3 \,\middle|\,\mathcal F_\sigma\right].
\end{aligned}
\end{equation}

For a closed local solution $(u,\eta)$ in the sense of \cref{def:critical-solution}, we extend the nonlinear drift and noise coefficients to $[0,T]$ by setting

$$
\begin{aligned}
F&=-\chi_{0,\eta\wedge T}(u\otimes u), &\qquad f&=-\chi_{0,\eta\wedge T}\TN(u), \\
\Phi&=\chi_{0,\eta\wedge T}\Sigma(\cdot,u_-), &\qquad \Psi&=\chi_{0,\eta\wedge T}\mathcal G(\cdot,u_-,\cdot).
\end{aligned}
$$ Each of these four terms vanish after $\eta$.
The critical norm is required only on the closed solution interval $[0,\eta]$ and no critical integrability of a frozen continuation is assumed.
That norm and \cref{hyp:L} give the required integrability of the drift forcing terms and stochastic integrands almost surely.
Apply \eqref{eq:3.6_full_stop} to the zero-initial stochastic convolution after localization of its integrability quantities.
Restart the heat term by the semigroup identity and the deterministic convolutions by splitting their time integrals.
This gives the pathwise restart formula for every $u_0\in L^0(\Omega,\mathcal F_0;E)$.
For $0\le\sigma\le\tau\le\eta\wedge T$, the resulting formula restarts the solution from its post-jump value $u(\sigma)$ and retains a possible jump at $\tau$.

\section{Local critical theory}
\label{Section 4}

Throughout this section we assume \cref{hyp:T1} and \cref{hyp:L}.
We construct solutions with $u_0\in L^3(\Omega,\mathcal F_0;E)$ and use localization on $\mathcal F_0$ to obtain pathwise solvability and conditional stability for $u_0\in L^0(\Omega,\mathcal F_0;E)$.
We first solve an auxiliary equation with a cutoff depending on the cumulative critical norm and recover the original equation on a strictly positive stochastic interval.
A localized stability estimate then yields pathwise uniqueness in the full critical mild solution class, allowing consistent continuation at stopping times and the construction of a maximal solution.

For $h>0$, let $\mathcal X_h$ be the space of adapted $E$-valued c\`adl\`ag processes with paths in $\mathbb X_h$ almost surely and finite norm
\begin{equation*}
\|v\|_{\mathcal X_h} :=\left(\E\|v\|_{\mathbb X_h}^{3}\right)^{1/3}.
\end{equation*}
Processes are identified up to indistinguishability.
This definition uses measurable path norms and does not require Bochner measurability with values in the space of c\`adl\`ag paths equipped with the uniform norm.

Fix $\vartheta\in C^1([0,\infty);[0,1])$ such that $$ \vartheta=1\ \text{on }[0,1],\qquad \vartheta=0\ \text{on }[2,\infty),\qquad \ell_\vartheta:=\|\vartheta'\|_\infty<\infty.
$$ For $v\in\mathcal X_h$ and $\varepsilon>0$, define
\begin{equation*}
q_v(t):=\left(\int_0^t\|v(s)\|_5^5\dd s\right)^{1/5}, \qquad (\mathcal Q_\varepsilon v)(t) :=\vartheta\bigl(q_v(t)/\varepsilon\bigr)v(t).
\end{equation*}
The spatial $L^5$ norm is interpreted as an extended-valued Borel function on $E$; its values at exceptional times do not affect the time integral.
The cumulative norm $q_v$ has continuous adapted paths and $q_v(0)=0$.
Consequently, $\vartheta(q_v/\varepsilon)$ is predictable.
Since this factor depends only on the time variable, $\mathcal Q_\varepsilon$ preserves the divergence constraint.

For subsequent estimates, put
\begin{equation*}
\eta(h):=(h^{3/2}+h)^{1/3}.
\end{equation*}

The cumulative cutoff makes the deterministic nonlinear map globally Lipschitz, with a Lipschitz constant that becomes small as $\varepsilon$ decreases.
Together with the stochastic convolution estimates on short intervals, this permits an adapted contraction argument on a deterministic interval whose length is independent of the initial datum.
The following proposition establishes the resulting auxiliary solution, its cubic moment bound, and its stability with respect to the initial datum.

\begin{proposition}
[Adapted contraction with a cumulative cutoff] \label{prop:adapted-critical-construction} Fix $T_0>0$.
There exist deterministic numbers $\varepsilon>0$ and $h\in(0,T_0]$, depending only on $\nu,C_g,\vartheta,T_0$, and $L_{T_0}$, such that the following holds for every $u_0\in L^3(\Omega,\mathcal F_0;E)$.

There is a unique $U\in\mathcal X_h$ satisfying, almost surely simultaneously for all $0\le t\le h$,
\begin{equation}
\label{eq:4.1_auxiliary-critical-equation}
\begin{aligned}
U(t)=&S(t)u_0 -\mathcal K_h\bigl((\mathcal Q_\varepsilon U) \otimes(\mathcal Q_\varepsilon U)\bigr)(t) -\mathcal H_h\bigl(\TN(\mathcal Q_\varepsilon U)\bigr)(t) \\
&+\int_0^tS(t-s)\Sigma(s,U(s-))\dd W(s)+\int_0^t\int_Z S(t-s)\mathcal G(s,U(s-),z) \,\tPN(\mathrm ds,\mathrm dz).
\end{aligned}
\end{equation}
Here $U(0-)=u_0$.
The cutoff acts only on the deterministic nonlinear terms.

The solution satisfies
\begin{equation}
\label{eq:4.1_aux_crit}
\E\|U\|_{\mathbb X_h}^3 \le C_\nu\left(\E\|u_0\|_3^3 +K_{T_0}^3(h^{3/2}+h)\right).
\end{equation}
For two initial data $u_0,v_0$, with the same $\varepsilon,h$, stochastic coefficients, and driving noises, the corresponding auxiliary solutions satisfy
\begin{equation}
\label{eq:4.1_auxiliary-data-stability}
\E\|U^{u_0}-U^{v_0}\|_{\mathbb X_h}^3 \le C_\nu\E\|u_0-v_0\|_3^3.
\end{equation}
The construction assumes only a finite cubic moment of the initial datum.
\end{proposition}

\begin{proof}
Let $(v_n)$ be Cauchy in $\mathcal X_h$ and choose a subsequence such that $\sum_n\|v_{n+1}-v_n\|_{\mathcal X_h}<\infty$.
Tonelli's theorem and $\E X\le(\E X^3)^{1/3}$ give $\sum_n\|v_{n+1}-v_n\|_{\mathbb X_h}<\infty$ almost surely.
Thus this subsequence converges uniformly in $E$ and in $L^5(Q_h)$.
The limits agree in distributions, and the uniform limit is adapted and c\`adl\`ag.
Completeness of the filtration permits a common exceptional set to be discarded, and Fatou's lemma gives convergence in $\mathcal X_h$.
Hence $\mathcal X_h$ is complete.

The extended-valued map $v\mapsto\|v\|_5$ is lower semicontinuous on $E$, by Fatou's lemma after extraction of an almost everywhere convergent subsequence, and therefore it is Borel measurable.
Progressive measurability of adapted c\`adl\`ag processes and pathwise $L^5$ integrability imply that $q_v$ is continuous and adapted.
The scalar cutoff is predictable, and $\mathcal Q_\varepsilon v$ is adapted and c\`adl\`ag.
Moreover, H\"older's inequality and \cref{hyp:T1} give, almost surely, $$ (\mathcal Q_\varepsilon v)\otimes(\mathcal Q_\varepsilon v) \in L^{5/2}(Q_h;\R^{3\times3}), \qquad \TN(\mathcal Q_\varepsilon v) \in L^{5/3}(Q_h;\R^3).
$$ Both nonlinear drift terms are progressively measurable.
These are precisely the integrability classes required by \cref{lem:3.2_critical-convection} and \cref{lem:critical-taming}, respectively.

We claim
\begin{align}
\|\mathcal Q_\varepsilon v\|_{L^5(Q_h)} &\le 2\varepsilon,
\label{eq:4.1_cutoff-uniform-bound}\\
\|\mathcal Q_\varepsilon v-\mathcal Q_\varepsilon w\|_{L^5(Q_h)} &\le C_\vartheta\|v-w\|_{L^5(Q_h)}, \qquad C_\vartheta:=1+4\ell_\vartheta.
\label{eq:4.1_cutoff-global-lipschitz}
\end{align}
Put $$ r_v:=h\wedge\inf\{t\in[0,h]:q_v(t)\ge2\varepsilon\}, \qquad d:=\|v-w\|_{L^5(Q_h)}, $$ and define $r_w$ similarly, with $\inf\varnothing=\infty$.
Continuity of $q_v$ gives $\|v\|_{L^5(Q_{r_v})}\le2\varepsilon$.
The cutoff vanishes on $(r_v,h]$, proving \eqref{eq:4.1_cutoff-uniform-bound}.
The reverse triangle inequality gives
\begin{equation*}
|q_v(t)-q_w(t)|\le d,\qquad 0\le t\le h.
\end{equation*}
Suppose $r_v\le r_w$.
On $(0,r_v)$, write $$ \mathcal Q_\varepsilon v-\mathcal Q_\varepsilon w =\vartheta(q_v/\varepsilon)(v-w) +[\vartheta(q_v/\varepsilon)-\vartheta(q_w/\varepsilon)]w.
$$ Since $\|w\|_{L^5(Q_{r_w})}\le2\varepsilon$, its norm on this interval is at most $(1+2\ell_\vartheta)d$.
On $(r_v,r_w)$, the first cutoff vanishes and $|\vartheta(q_w/\varepsilon)| \le\ell_\vartheta d/\varepsilon$, giving a further bound $2\ell_\vartheta d$.
Both cutoffs vanish after $r_w$, which proves \eqref{eq:4.1_cutoff-global-lipschitz}.

Set $\mathcal C_{N,h}(a):=\mathcal H_h(\TN(a))$.
Choose $b_\nu,c_{\nu,g}$ as constants in the bounds \eqref{eq:3.2_convection}--\eqref{eq:3.2_critical-diff} and \eqref{eq:3.3_tame1}--\eqref{eq:3.3_tame2}, respectively; these constants are independent of $h$.
For $$ \mathcal D_\varepsilon(v) :=-\mathcal B_h(\mathcal Q_\varepsilon v, \mathcal Q_\varepsilon v) -\mathcal C_{N,h}(\mathcal Q_\varepsilon v), $$ the cutoff bounds imply, pathwise,
\begin{align*}
\|\mathcal D_\varepsilon(v)\|_{\mathbb X_h} &\le4b_\nu\varepsilon^2+8c_{\nu,g}\varepsilon^3, & \|\mathcal D_\varepsilon(v)-\mathcal D_\varepsilon(w)\|_{\mathbb X_h} &\le\Lambda_\varepsilon\|v-w\|_{L^5(Q_h)},
\end{align*}
where $\Lambda_\varepsilon :=C_\vartheta(4b_\nu\varepsilon +8c_{\nu,g}\varepsilon^2).$ $\mathcal D_\varepsilon$ is $\Lambda_\varepsilon$-Lipschitz on $\mathcal X_h$ and vanishes at the zero process.
\cref{lem:3.2_critical-convection} and \cref{lem:critical-taming} provide continuous $E$-valued representatives vanishing at time zero.
For each $t$, the defining integral representations of $\mathcal K_h$ and $\mathcal H_h$ involve only the restrictions of their respective integrands to $(0,t)$ and thus we obtain the adaptedness of these representatives.
The pathwise cutoff bound above ensures finite cubic moments.

Let $\mathcal W_h(v)$ denote the sum of the Wiener and compensated Poisson stochastic convolutions in \eqref{eq:4.1_auxiliary-critical-equation}, with integrands $\Sigma(\cdot,v_-)$ and $\mathcal G(\cdot,v_-,\cdot)$, where $v(0-):=v(0)$.
Predictability and hypothesis \cref{hyp:L} give
\begin{equation*}
\begin{aligned}
&\mathcal I_h\bigl(0,\Sigma(\cdot,v_-), \mathcal G(\cdot,v_-,\cdot)\bigr)\le C\eta(h)^3 \left(K_{T_0}^3+L_{T_0}^3 \E\sup_{0\le s\le h}\|v(s)\|_3^3\right).
\end{aligned}
\end{equation*}
Before taking expectations, the quadratic Wiener and Poisson integrability terms are bounded by $h^{3/2}(K_{T_0}+L_{T_0}\sup_{s\le h}\|v(s)\|_3)^3$, whereas the cubic Poisson integrability term has the corresponding bound with $h$ in place of $h^{3/2}$.
\cref{thm:critical-stokes} therefore gives $\mathcal W_h(v)\in\mathcal X_h$, with zero initial value, and
\begin{align*}
\|\mathcal W_h(v)-\mathcal W_h(w)\|_{\mathcal X_h} &\le C_\nu L_{T_0}\eta(h)\|v-w\|_{\mathcal X_h}, & \|\mathcal W_h(0)\|_{\mathcal X_h} &\le C_\nu K_{T_0}\eta(h).
\end{align*}

Choose $\varepsilon>0$ and $h\in(0,T_0]$ so that $$ \Lambda_\varepsilon\le\tfrac14, \qquad C_\nu L_{T_0}\eta(h)\le\tfrac14.
$$ \cref{lem:heat-extension} and the preceding estimates show that $\mathcal F(v):=S(\cdot)u_0+ \mathcal D_\varepsilon(v)+\mathcal W_h(v)$ maps $\mathcal X_h$ into itself and satisfies
\begin{equation*}
\|\mathcal F(v)-\mathcal F(w)\|_{\mathcal X_h} \le\frac12\|v-w\|_{\mathcal X_h}.
\end{equation*}
The Banach fixed-point theorem, starting from the adapted iterate $U^{(0)}=0$, gives a unique $U\in\mathcal X_h$.
Equality in this space is indistinguishability, so the equation holds simultaneously at all times and $U(0)=u_0$.
Since $\mathcal D_\varepsilon(0)=0$, $$ \|U\|_{\mathcal X_h} \le\tfrac12\|U\|_{\mathcal X_h} +C_\nu\bigl((\E\|u_0\|_3^3)^{1/3} +K_{T_0}\eta(h)\bigr).
$$ Absorption and cubing prove \eqref{eq:4.1_aux_crit}.
For two initial data, the same contraction gives $$ \|U^{u_0}-U^{v_0}\|_{\mathcal X_h} \le\tfrac12\|U^{u_0}-U^{v_0}\|_{\mathcal X_h} +C_\nu(\E\|u_0-v_0\|_3^3)^{1/3}, $$ which proves \eqref{eq:4.1_auxiliary-data-stability}.
\end{proof}

The auxiliary equation coincides with \eqref{eq:tnse} while the cumulative critical norm remains below the level at which the cutoff becomes active.
Since this cumulative norm is continuous and vanishes initially, the corresponding stopping time is strictly positive almost surely.
The following theorem uses this observation to construct a closed local critical mild solution.

\begin{theorem}
[Local critical existence without smallness of the initial datum] \label[theorem]{thm:local-critical-existence} Assume \cref{hyp:T1}, \cref{hyp:L}, and $u_0\in L^3(\Omega,\mathcal F_0;E)$.
Fix $T_0>0$, and choose $\varepsilon,h$ as in \cref{prop:adapted-critical-construction}.
There exist a stopping time $\tau$, with $$ \P(0<\tau\le h)=1, $$ and a closed local critical mild solution $(u,\tau)$ of \eqref{eq:tnse} in the sense of \cref{def:critical-solution}.
The solution is constructed on the prescribed stochastic basis and satisfies
\begin{align*}
\int_0^\tau\|u(s)\|_5^5\dd s &\le\varepsilon^5 \quad\text{almost surely}, \\
\E\left[ \sup_{0\le s\le\tau}\|u(s)\|_3^3 +\left(\int_0^\tau\|u(s)\|_5^5\dd s\right)^{3/5} \right] &\le C_\nu\left( \E\|u_0\|_3^3 +K_{T_0}^3(h^{3/2}+h) \right).
\end{align*}
Moreover, at small times it satisfies
\begin{equation}
\label{eq:4.2_loc_crit}
\lim_{t\downarrow0} \E\|u-S(\cdot)u_0\|_{\mathbb X[0,t\wedge\tau]}^3=0.
\end{equation}
There exist strictly positive bounded stopping times $\tau_n\uparrow\tau$, with $\tau_n<\tau$ almost surely, such that each $(u(\cdot\wedge\tau_n),\tau_n)$ is a closed local critical mild solution.
Consequently, the restriction of $u$ to $[0,\tau)$ defines a local critical mild solution with lifetime $\tau$ in the sense of \cref{def:critical-solution}.

The jump identity holds simultaneously for all $t>0$:
\begin{equation}
\label{eq:4.2_loc_sol}
\tilde\Delta u(t) =\mathds1_{\{0<t\le\tau\}} \int_Z\mathcal G(t,u(t-),z)\, \PN(\{t\},\mathrm dz).
\end{equation}

\end{theorem}

\begin{proof}
Let $U$ denote the auxiliary solution constructed in \cref{prop:adapted-critical-construction}, and define
\begin{equation*}
\tau:=h\wedge \inf\{t\in[0,h]:q_U(t)\ge\varepsilon\}, \qquad \inf\varnothing:=\infty.
\end{equation*}
Since $q_U$ is continuous and adapted, with $q_U(0)=0$, the random time $\tau$ is a strictly positive stopping time.
Moreover,
\begin{equation*}
q_U(t)\le\varepsilon,\qquad \mathcal Q_\varepsilon U(t)=U(t), \qquad 0\le t\le\tau.
\end{equation*}
Define
\begin{equation*}
u(t):=U(t\wedge\tau).
\end{equation*}
Then $u$ is adapted and $E$-valued c\`adl\`ag, and $u(t)=u(\tau)$ for every $t\ge\tau$.
By the definition of $q_U$, $$ \int_0^\tau\|u(s)\|_5^5\dd s =q_U(\tau)^5\le\varepsilon^5.
$$ The asserted moment estimate follows from \eqref{eq:4.1_aux_crit} by restricting the critical norm of $U$ to $[0,\tau]$.

To verify the closed mild formulation, apply \cref{lem:stopping-time-restart} to the deterministic drift terms $$ F_\varepsilon :=-(\mathcal Q_\varepsilon U) \otimes(\mathcal Q_\varepsilon U), \qquad f_\varepsilon :=-\TN(\mathcal Q_\varepsilon U), $$ and the predictable stochastic integrands $$ \Phi_U(s):=\Sigma(s,U(s-)), \qquad \Psi_U(s,z):=\mathcal G(s,U(s-),z).
$$ The construction in \cref{prop:adapted-critical-construction} establishes, almost surely, $$ F_\varepsilon \in L^{5/2}(Q_h;\R^{3\times3}), \qquad f_\varepsilon \in L^{5/3}(Q_h;\R^3).
$$ These are precisely the integrability conditions required for the deterministic convolution operators $\mathcal K_h$ and $\mathcal H_h$, respectively, by \cref{lem:3.2_critical-convection} and \cref{lem:critical-taming}.
The stochastic integrands $\Phi_U,\Psi_U$ satisfy the Wiener and Poisson integrability conditions of \cref{prop:linear-cubic-ito}, as established in the auxiliary construction.

Retaining the initial term $S(t)u_0$ and multiplying the drift terms and stochastic integrands by $\chi_{0,\tau}$ defines the process through the corresponding mild integral representation.
By \cref{lem:stopping-time-restart}, this process coincides with $U$ simultaneously on $[0,\tau]$, outside one null set.
Since $\mathcal Q_\varepsilon U=U$ on this interval and $u(s-)=U(s-)$ for $0<s\le\tau$, the restricted drift terms and stochastic integrands coincide with those in the stopped mild formulation \eqref{eq:2.4_stop_Con}.

Extend the restricted drift terms and stochastic integrands by zero beyond $h$ to each larger deterministic time interval.
By the semigroup continuation identity for the truncated drift terms and stochastic integrands \eqref{eq:3.6_heat_continuation}, the corresponding $E$-valued c\`adl\`ag modifications agree on overlapping deterministic time intervals and satisfy $$ Y_\tau(t)=S(t-\tau)U(\tau), \qquad t\ge\tau, $$ and $$ u(t\wedge\tau)=Y_\tau(t\wedge\tau), \qquad t\ge0.
$$ Taking a countable increasing sequence of deterministic time horizons and intersecting the corresponding full-probability events yields these identities simultaneously on $[0,\infty)$.
Consequently, $(u,\tau)$ satisfies \cref{def:critical-solution}.

Set $$ V:=U-S(\cdot)u_0.
$$ Then $V\in\mathcal X_h$ and $V(0)=0$.
The right continuity of its $E$-valued paths and the absolute continuity of the integral $\int_0^t\|V(s)\|_5^5\dd s$ imply $$ \lim_{t\downarrow0} \|V\|_{\mathbb X[0,t]}=0 \qquad\text{almost surely}.
$$ Since $$ \|u-S(\cdot)u_0\|_{\mathbb X[0,t\wedge\tau]}^3 \le\|V\|_{\mathbb X_h}^3 $$ and $\E\|V\|_{\mathbb X_h}^3<\infty$, dominated convergence proves \eqref{eq:4.2_loc_crit}.

The deterministic Duhamel convolutions are continuous as $E$-valued functions by \cref{lem:3.2_critical-convection} and \cref{lem:critical-taming}.
Consequently, the jump identity from \cref{lem:stopping-time-restart}, together with the definition of the stopped process, gives $$ \tilde\Delta u(t) =\mathds1_{\{0<t\le\tau\}}\tilde\Delta U(t) \qquad (0<t\le h).
$$ The process $u$ is constant after $\tau$ and hence has no subsequent jumps.
On $[0,\tau]$, the predictable Poisson integrand is $\mathcal G(t,u(t-),z)$.
These observations prove \eqref{eq:4.2_loc_sol}, including the terminal value at $t=\tau$.

Finally, for $n\ge1$ let $a_n=n/(n+1)$ and define
\begin{equation*}
\tau_n:=a_nh\wedge \inf\{t\in[0,h]:q_U(t)\ge a_n\varepsilon\}.
\end{equation*}
Since $q_U$ is continuous and adapted, each $\tau_n$ is a strictly positive bounded stopping time, and the sequence $(\tau_n)$ is increasing.
If $\tau=h$, the deterministic upper bounds $a_nh$ imply $\tau_n<\tau$.
If $\tau<h$, continuity of $q_U$ and the strict inequality $a_n\varepsilon<\varepsilon$ imply that the level $a_n\varepsilon$ is reached strictly before $\tau$.
Thus $\tau_n<\tau$ almost surely in either case.
Also for every $t<\tau$, $q_U(t)<\varepsilon.$

Since $a_n\uparrow1$, eventually $a_nh>t$ and $a_n\varepsilon>q_U(t)$, so $\tau_n>t$.
It follows that $$ \tau_n\uparrow\tau \qquad\text{almost surely}.
$$ Applying the same semigroup continuation argument for the truncated drift terms and stochastic integrands at each $\tau_n$ shows that $(u(\cdot\wedge\tau_n),\tau_n)$ is a closed local critical mild solution.
The corresponding restriction to $[0,\tau)$ is therefore a local critical mild solution with lifetime $\tau$.
\end{proof}

\paragraph{Dependence of the existence interval on the initial datum}
The auxiliary equation is solved on the deterministic time interval $[0,h]$, whereas the solution of the original tamed stochastic NSEn is constructed on the random interval $[0,\tau]$.
Strict positivity of $\tau$ follows from the continuity of the cumulative critical norm $q_U$ and the identity $q_U(0)=0$.
The construction therefore does not require uniform small-time $L^5$ estimates for heat extensions of initial data belonging to bounded subsets of $L^3_\sigma(\R^3)$.

We next compare arbitrary closed local critical mild solutions in the sense of \cref{def:critical-solution}, defined on the same stochastic basis with the same coefficient fields and driving noises.

\begin{proposition}
[Pathwise uniqueness and closed localized stability] \label{prop:critical-uniqueness-stability} Assume \cref{hyp:T1} and \cref{hyp:L}.
Let $(u,\tau_u)$ and $(v,\tau_v)$ be any two closed local critical mild solutions of \eqref{eq:tnse} in the sense of \cref{def:critical-solution}, with respective initial data $u_0,v_0\in L^0(\Omega,\mathcal F_0;E)$.
Both solutions are defined on the same prescribed stochastic basis, with the same coefficient fields $\Sigma,\mathcal G$, Wiener process $W$, and Poisson random measure $\PN$.
Fix $T,M>0$ and set $\lambda:=T\wedge\tau_u\wedge\tau_v$.
Define $$ C_{u,v}^{\lambda}(t) :=\sup_{0\le s\le t\wedge\lambda} \bigl(\|u(s)\|_3+\|v(s)\|_3\bigr)+\left( \int_0^{t\wedge\lambda} \bigl(\|u(s)\|_5^5+\|v(s)\|_5^5\bigr) \dd s \right)^{1/5}, $$ $$ \rho:=\lambda\wedge \inf\{t\ge0:C_{u,v}^{\lambda}(t)\ge M\}.$$

Write $d=\|u_0-v_0\|_3$.
There is a deterministic constant $1\le C_{T,M}<\infty$, depending only on $T,M,\nu,C_g$ and the coefficient bounds on $[0,T]$, such that almost surely $$ d^3\le \E\!\left[\|u-v\|_{\mathbb X[0,\rho]}^3\mid\mathcal F_0\right] \le C_{T,M}d^3 .
$$ Consequently, whenever $\E d^3<\infty$,
\begin{equation}
\label{eq:4.3_close_loc}
\E\|u-v\|_{\mathbb X[0,\rho]}^3\le C_{T,M}\E d^3 .
\end{equation}
The supremum in this estimate includes the value at $\rho$, including when $\rho$ coincides with a jump time.

If $u_0=v_0$ almost surely, then $u=v$ simultaneously on $[0,\tau_u\wedge\tau_v]$ almost surely.
Any two local critical solutions with the same initial datum, coefficients, and driving noises agree on the intersection of their stochastic intervals of existence.
The analogous uniqueness assertion holds for solutions restarted at a stopping time from the same measurable state.
\end{proposition}

\begin{proof}

The process $C_{u,v}^{\lambda}$ is adapted, nondecreasing and c\`adl\`ag, so $\rho$ is a stopping time.

\emph{Data with finite third moments.} First suppose $u_0,v_0\in L^3(\Omega,\mathcal F_0;E)$.
Continuity of its integral component gives
\begin{equation*}
\int_0^\rho \bigl(\|u(s)\|_5^5+\|v(s)\|_5^5\bigr) \dd s \le M^5, \qquad \|u(s-)\|_3+\|v(s-)\|_3\le M \quad(0<s\le\rho).
\end{equation*}
The inequality involving the left limits is vacuous when $\rho=0$ and does not impose a bound on the values of $u$ and $v$ at $\rho$ after a possible jump.

The convective tensor fields and taming terms, restricted to $(0,\rho]$, have respective $L^{5/2}((0,\rho)\times\R^3)$- and $L^{5/3}((0,\rho)\times\R^3)$-norms bounded by $M^2$ and $C_gM^3$.
Hypothesis \cref{hyp:L} bounds the corresponding Wiener and Poisson integrability quantities by $C(T^{3/2}+T)(K_T+L_TM)^3$.
The mild representation with the deterministic drift terms and stochastic integrands restricted to $(0,\rho]$, \cref{lem:heat-extension}--\cref{lem:critical-taming}, and \cref{thm:critical-stokes} therefore give
\begin{equation}
\label{eq:4.3_comp}
\E\|u\|_{\mathbb X[0,\rho]}^3 +\E\|v\|_{\mathbb X[0,\rho]}^3<\infty.
\end{equation}
This estimate includes the terminal value at $\rho$ and the case $\rho=0$, in which only the initial data contribute.

Set $w=u-v$ and consider stopping times $0\le a\le b\le\rho$ satisfying
\begin{equation}
\label{eq:4.3_comp_s}
b-a\le h, \qquad \int_a^b \bigl(\|u(s)\|_5^5+\|v(s)\|_5^5\bigr) \dd s \le\delta^5 \quad\text{almost surely}.
\end{equation}
\cref{lem:stopping-time-restart} yields, outside one null set and simultaneously for all $a\le t\le b$,
\begin{equation*}
\begin{aligned}
w(t) =&S(t-a)w(a) -\mathcal K_T\bigl( \chi_{a,b}(u\otimes u-v\otimes v) \bigr)(t) -\mathcal H_T\bigl( \chi_{a,b}(\TN(u)-\TN(v)) \bigr)(t) +\mathcal Z_{a,b} [\Phi_{u,v},\Psi_{u,v}](t).
\end{aligned}
\end{equation*}
The Wiener and Poisson coefficient differences are evaluated at the left limits of the processes.
The initial state $w(a)$ includes a possible jump at $a$, whereas the subsequent stochastic integrals are taken over $(a,b]$.

Choose $b_\nu,c_{\nu,g}$ as in the proof of \cref{prop:adapted-critical-construction}.
By \eqref{eq:3.6_convection_forcing}-- \eqref{eq:3.6_cubic_forcing} and \cref{lem:taming}, the two deterministic convolution differences have combined norm bounded by $ d_\delta\|w\|_{\mathbb X[a,b]}, \ d_\delta=2b_\nu\delta+2c_{\nu,g}\delta^2.
$ The heat semigroup contribution is bounded by $C_\nu\|w(a)\|_3$.
\cref{eq:3.6_random_interval_noise_lipschitz} bounds the conditional third-moment norm of the stochastic convolution difference by $$ C_\nu L_T\eta(h) \left( \E[ \|w\|_{\mathbb X[a,b]}^3 \mid\mathcal F_a] \right)^{1/3}.
$$ Choose $\delta>0$ and $h\in(0,T]$ such that $d_\delta\le\frac14, \ C_\nu L_T\eta(h)\le\frac14.$ 
Conditional Minkowski's inequality and absorption give
\begin{equation}
\label{eq:4.3_cond_stb}
\E[ \|w\|_{\mathbb X[a,b]}^3 \mid\mathcal F_a] \le C_*\|w(a)\|_3^3 \quad\text{almost surely},
\end{equation}
for a deterministic constant $C_*\ge1$.
All conditional expectations are finite by \eqref{eq:4.3_comp}.

Set $$ A(t):=\int_0^{t\wedge\rho} \bigl(\|u(s)\|_5^5+\|v(s)\|_5^5\bigr) \dd s, \qquad t_0=0, $$ and define $$t_{j+1}:= \rho\wedge(t_j+h)\wedge \inf\{t\ge t_j:A(t)-A(t_j)\ge\delta^5\}.$$

Continuity and adaptedness of $A$ imply that these random times are stopping times.
Each interval satisfies \eqref{eq:4.3_comp_s}.
Before reaching $\rho$, each interval either has length $h$ or the cumulative $L^5$ integral increases by $\delta^5$.
Since $\rho\le T$ and $A(\rho)\le M^5$, $$ J:=2+\lceil T/h\rceil+\lceil M^5/\delta^5\rceil \quad\text{satisfies}\quad t_J=\rho.
$$ After reaching $\rho$, the subsequent stopping times are taken equal to $\rho$.
Iteration of \eqref{eq:4.3_cond_stb}, including the terminal value of each interval, gives $$ \E\|w\|_{\mathbb X[t_j,t_{j+1}]}^3 \le C_*^{j+1}\E\|w(0)\|_3^3.
$$ The cube of the supremum norm over the entire interval is bounded by the sum of the cubes of the suprema over the partition intervals.
Moreover, subadditivity of $r\mapsto r^{3/5}$ bounds the cubic power of the full $L^5$ norm by the sum of the corresponding quantities over the partition intervals.
Using $(x+y)^3\le4(x^3+y^3)$ proves \eqref{eq:4.3_close_loc}, with $C_{T,M}=4\sum_{j=1}^{J}C_*^j$.

Now let $u_0,v_0\in L^0(\Omega,\mathcal F_0;E)$, put $d=\|u_0-v_0\|_3$, and set $B=\{\|u_0\|_3+\|v_0\|_3<M\}$.
On $B^c$ the threshold is reached at time zero, so $\rho=0$ and $\|u-v\|_{\mathbb X[0,\rho]}=d$.
For $H\in\mathcal F_0$ with $H\subseteq B$, set $$ u^H=\mathds1_Hu,\qquad v^H=\mathds1_Hv,\qquad \Sigma^H=\mathds1_H\Sigma,\qquad \mathcal G^H=\mathds1_H\mathcal G .
$$ The restricted data are bounded by $M$ in $E$, and the new coefficient fields satisfy \cref{hyp:L} with the same deterministic bounds.
With their original closed stopping times, $u^H,v^H$ solve the equation with these fields.
Indeed, $(\mathds1_Hu)\otimes(\mathds1_Hu)=\mathds1_H(u\otimes u)$, $\TN(\mathds1_Hu)=\mathds1_H\TN(u)$, and $\Sigma^H(t,u^H(t-))=\mathds1_H\Sigma(t,u(t-))$, with the same identity for $\mathcal G^H$.
Equation \eqref{eq:3.6_restart_locality} applied after localization at initial time zero gives the corresponding identities for the localized stochastic integrals.
All identities hold up to and including the closed endpoints.

On $H$ the stopping time for the restricted pair equals $\rho$.
On $H^c$ the difference vanishes throughout its common interval.
The finite-moment estimate just proved therefore gives $$ \E[\mathds1_H\|u-v\|_{\mathbb X[0,\rho]}^3] \le C_{T,M}\E[\mathds1_Hd^3].
$$ Testing all such $H$ proves the conditional upper bound on $B$.
On $B^c$ the conditional expectation equals $d^3$ because $B$ and $d$ are $\mathcal F_0$-measurable.
The conditional expectation is understood in the nonnegative sense.
Its lower bound follows from $0\in[0,\rho]$.
Integration gives \eqref{eq:4.3_close_loc} for every pair with $\E d^3<\infty$.

\emph{Pathwise uniqueness.} If $u_0=v_0$, the stability estimate implies that $u=v$ on $[0,\rho]$.
Almost sure finiteness of the critical norms on the closed solution intervals implies $\rho=\lambda$ for all sufficiently large integers $M$.
Taking countable unions of exceptional null sets and then letting integer $T\uparrow\infty$ proves uniqueness on the common closed interval of existence.
Applying this result to pairs of increasing sequences of stopping times approaching the respective lifetimes gives uniqueness on the corresponding open intervals of existence.

The same conditional argument applies from a bounded stopping time $a$ when the restart state belongs to $L^3(\Omega,\mathcal F_a;E)$.
For an almost surely finite restart state $\xi$, multiply the restarted difference identity by the indicators of the events $$ \{k-1\le\|\xi\|_3<k\}\in\mathcal F_a.
$$ \cref{eq:3.6_restart_locality} justifies this localization of the stochastic integrals, and the restricted initial states have finite cubic moments.
For a potentially unbounded restart time, also restrict to $\{a\le\ell\}\in\mathcal F_a$ and bounded deterministic time intervals.
Intersecting the full-probability events for the countably many localizations and then exhausting the events on which the restart time and restart state are bounded proves the asserted uniqueness on \(\{a<\infty\}\).
\end{proof}

Local existence and pathwise uniqueness now allow successive solutions to be joined at stopping times into a single maximal critical solution.
Under the globally Lipschitz hypothesis \cref{hyp:L}, the construction shows that a finite maximal lifetime forces divergence of the critical spacetime integral, almost surely.
The following theorem also identifies the terminal $E$-valued state of a local solution whose critical norms remain finite at a finite lifetime, thereby permitting continuation on the prescribed stochastic basis.

\begin{theorem}
[Maximal strong critical solution and continuation] \label[theorem]{thm:maximal-critical-solution} Assume \cref{hyp:T1}, \cref{hyp:L}, and $u_0\in L^3(\Omega,\mathcal F_0;E)$.
There is a pathwise unique maximal probabilistically strong critical solution $(u,\{\tau_n\}_{n\ge1},\tau_*)$ in the sense of \cref{def:critical-solution}, with $\tau_*>0$ almost surely.
Each $\tau_n$ is bounded and the solution on $[0,\tau_n]$ includes its terminal value.

The maximal solution is independent, up to indistinguishability, of the auxiliary cutoff and of the choices of construction intervals.
Every other local critical solution with the same initial datum, coefficients, and driving noises is its restriction, and its lifetime is at most $\tau_*$ almost surely.
The combined intrinsic continuation criterion is
\begin{equation}
\label{eq:4.4_comb_crit}
\P\left( \tau_*<\infty,\quad \sup_{0\le t<\tau_*}\|u(t)\|_3<\infty,\quad \int_0^{\tau_*}\|u(t)\|_5^5\dd t<\infty \right)=0.
\end{equation}
In fact, under the globally Lipschitz hypothesis \cref{hyp:L} used here, the construction gives the stronger spacetime criterion
\begin{equation}
\label{eq:4.4_crit_space}
\P\left( \tau_*<\infty,\quad \int_0^{\tau_*}\|u(t)\|_5^5\dd t<\infty \right)=0.
\end{equation}

For any local solution $(v,\{\kappa_n\},\zeta)$ with the same initial datum, coefficients, and driving noises, on the event
\begin{equation*}
B:=\left\{\zeta<\infty,\quad \sup_{t<\zeta}\|v(t)\|_3<\infty,\quad \int_0^\zeta\|v(t)\|_5^5\dd t<\infty\right\},
\end{equation*}
one has $\zeta<\tau_*$ and a terminal state
\begin{equation}
\label{eq:4.4_crit_term}
v_\zeta:=\lim_{t\uparrow\zeta}v(t)=u(\zeta) \quad\text{in }E.
\end{equation}
The solution continues from this state on the prescribed stochastic basis.
No finite-energy assumption on $u_0$ and no hypothesis \cref{hyp:G} are used.
\end{theorem}

\begin{proof}

Fix $H>0$.
Choose $\varepsilon$ and $h_H\in(0,H]$ by \cref{prop:adapted-critical-construction}, using the coefficient bounds on $[0,H]$.
The same $\varepsilon$ works for all $H$, since only the drift estimates enter its choice.
For a stopping time $\sigma\le H$ and $\xi\in L^3(\Omega,\mathcal F_\sigma;E)$, put $\beta=(\sigma+h_H)\wedge H$ and consider the complete space of adapted $E$-valued c\`adl\`ag processes on $[\sigma,\beta]$ with norm $$ \|V\|_{\mathscr X[\sigma,\beta]} :=(\E\|V\|_{\mathbb X[\sigma,\beta]}^3)^{1/3}.
$$ Processes are identified up to indistinguishability on the stochastic interval $[\sigma,\beta]$.
Extend these processes by zero before $\sigma$ and keep them constant after $\beta$; the $L^5$ condition is imposed only on $[\sigma,\beta]$.
Completeness follows by the preceding uniform-$E$ and $L^5$ convergence argument.
Define the continuous adapted cumulative $L^5$-norm process and its cutoff by $$ q_{\sigma,V}(s) =\left(\int_0^s\chi_{\sigma,\beta}(r)\|V(r)\|_5^5\dd r\right)^{1/5}, \qquad \mathcal Q_{\varepsilon,\sigma}V(s) =\vartheta(q_{\sigma,V}(s)/\varepsilon)V(s).
$$ For $\sigma\le t\le\beta$, consider
\begin{equation*}
\begin{aligned}
\mathcal F_{\sigma,\xi}(V)(t) :=&S(t-\sigma)\xi -\mathcal K_H\bigl(\chi_{\sigma,\beta} (\mathcal Q_{\varepsilon,\sigma}V) \otimes(\mathcal Q_{\varepsilon,\sigma}V)\bigr)(t) \\
&\qquad-\mathcal H_H\bigl(\chi_{\sigma,\beta} \TN(\mathcal Q_{\varepsilon,\sigma}V)\bigr)(t) +\mathcal Z_{\sigma,\beta} [\Sigma(\cdot,V_-),\mathcal G(\cdot,V_-,\cdot)](t).
\end{aligned}
\end{equation*}
The stochastic convolutions are formed at deterministic terminal times before restriction to $[\sigma,\beta]$.
The cutoff estimates hold pathwise after translating the initial time.  \cref{lem:stopping-time-restart} supplies the deterministic convolution estimates and conditional stochastic convolution estimates on the random interval.
Thus the map is adapted and has Lipschitz constant $\Lambda_\varepsilon+C_\nu L_H\eta(h_H)\le1/2$.
The heat semigroup term and the stochastic convolution evaluated at the zero process have finite cubic moments, giving a unique fixed point $U^{\sigma,\xi}$.
Conditional Minkowski's inequality and absorption, subtracting the stochastic convolution evaluated at the zero process and using that the cutoff drift vanishes at zero, give
\begin{equation}
\label{eq:4.4_shif_const}
\E\left[ \|U^{\sigma,\xi}\|_{\mathbb X[\sigma,\beta]}^3 \mid\mathcal F_\sigma\right] \le C_\nu\bigl(\|\xi\|_3^3+K_H^3\eta(h_H)^3\bigr).
\end{equation}
Set
\begin{equation*}
\sigma^+:=\beta\wedge \inf\{t\ge\sigma:q_{\sigma,U^{\sigma,\xi}}(t) \ge\varepsilon\}.
\end{equation*}
On $\{\sigma<H\}$, continuity gives $\sigma^+>\sigma$.
The cutoff is one on the closed interval $[\sigma,\sigma^+]$; the argument using the truncated drift terms and stochastic integrands in \cref{thm:local-critical-existence} therefore gives the original equation with initial value $\xi$ and stochastic integrals over $(\sigma,\sigma^+]$.
Moreover,
\begin{equation*}
\int_\sigma^{\sigma^+}\|U^{\sigma,\xi}(s)\|_5^5\dd s \le\varepsilon^5,
\end{equation*}
with equality if $\sigma^+<\beta$.
Starting from $s_0=0$, $\xi_0=u_0$, iterate the construction with $\sigma=s_n$, $\xi_n=u(s_n)$, and $s_{n+1}=s_n^+$.
If $s_n=H$, set all subsequent stopping times equal to $H$.
\cref{eq:4.4_shif_const} proves inductively that every restart datum has finite cubic moment.
The preceding interval includes the jump at $s_n$, whereas the subsequent stochastic increment excludes it.
\cref{lem:stopping-time-restart} and the semigroup identity therefore combine the two representations into the original closed mild equation on $[0,s_{n+1}]$.
Finite sums of the interval suprema cubed and subadditivity of the $3/5$ power of the $L^5$ integrals give finite cubic critical moments at every stage.
Put $\zeta_H=\lim_n s_n\le H$.
The compatible sample paths define $u^H$ on $[0,\zeta_H)$, with announcing times
\begin{equation*}
\lambda_n^H:=s_n\wedge\frac{nH}{n+1},\qquad n\ge1.
\end{equation*}
They are positive, increasing, bounded, and converge to $\zeta_H$.
If $\zeta_H<H$, strict advance gives $s_n<\zeta_H$; if $\zeta_H=H$, the deterministic upper bound is strictly below $H$.
Thus $\lambda_n^H<\zeta_H$, and each closed restriction is obtained from a finite-stage solution.

On $\{\zeta_H<H\}$ there are infinitely many stages, each of which either has length $h_H$ or attains the cumulative $L^5$ integral threshold $\varepsilon^5$.
Only finitely many can have length $h_H$.
Additivity of the integral therefore gives
\begin{equation}
\label{eq:4.4_fin_horizon}
\int_0^{\zeta_H}\|u^H(s)\|_5^5\dd s=\infty \quad\text{on }\{\zeta_H<H\},\quad\text{almost surely}.
\end{equation}
Both $h_H$ and $\varepsilon$ remain fixed throughout the construction on $[0,H]$, independently of the restart state.

Construct $u^m$ for $H=m\in\N$.
\cref{prop:critical-uniqueness-stability}, applied to pairs of announcing times, identifies these solutions on their common intervals of existence.
In fact,
\begin{equation}
\label{eq:4.4_hor_lifetime}
\zeta_m=\zeta_{m+1}\wedge m \quad\text{almost surely}.
\end{equation}
If $\zeta_m<\min(\zeta_{m+1},m)$, the divergence at $\zeta_m$ contradicts the finite critical norm on an announcing interval of $u^{m+1}$ extending beyond it.
If $\zeta_{m+1}<\zeta_m$, interchange the two solutions.
Each event is covered by countably many such announcing intervals, which proves the identity.

Set $\tau_*=\lim_m\zeta_m$ and use the consistent sample paths to define $u$ on $[0,\tau_*)=\bigcup_m[0,\zeta_m)$.
Countability gives a common full-probability set; adaptedness follows from $\{t<\zeta_m\}\in\mathcal F_t$.
The first construction interval gives $\tau_*>0$.
Define
\begin{equation*}
\tau_n:=\max_{1\le m\le n}\lambda_n^m.
\end{equation*}
Then $0<\tau_n\le n$, $\tau_n$ is increasing, $\tau_n\uparrow\tau_*$, and $\tau_n<\tau_*$ whenever $\tau_*<\infty$.
The finitely many agreeing closed intervals $[0,\lambda_n^m]$, $m\le n$, cover $[0,\tau_n]$.
Their finite moment bounds apply to the union; the semigroup continuation identity for the restricted drift terms and stochastic integrands gives the closed equation there.
Thus these are valid localizations in \cref{def:critical-solution}.

On $\{\tau_*<m\}$, \eqref{eq:4.4_hor_lifetime} gives $\zeta_m=\tau_*<m$.
\cref{eq:4.4_fin_horizon}, followed by a union over $m$, proves \eqref{eq:4.4_crit_space} and hence \eqref{eq:4.4_comb_crit}.
Any other local solution $v$ with lifetime $\zeta$ agrees with $u$ before $\zeta\wedge\tau_*$.
On $\{\zeta>\tau_*\}$, one of its announcing intervals extends beyond $\tau_*$, contradicting this divergence.
Thus $\zeta\le\tau_*$ and $v$ is a restriction of $u$.
This proves maximality, uniqueness of the maximal lifetime, and independence of the construction choices.

For the local solution in the statement, the continuation criterion implies $\zeta<\tau_*$ on $B$.
The path of $u$ is c\`adl\`ag on a neighborhood of $\zeta$ there, so $\lim_{t\uparrow\zeta}v(t)=u(\zeta-)$ exists in $E$.
The stopping time $\zeta$, announced by $(\kappa_n)$, has a predictable graph: on $(0,\infty)$, $$ \mathds1_{\{s=\zeta<\infty\}} =\lim_n\mathds1_{\{\kappa_n<s\le\zeta\}}.
$$ Choose $Z_j\uparrow Z$ with $\mu(Z_j)<\infty$.
For every finite $T$, compensation gives $$ \E\int_0^T\!\int_{Z_j}\mathds1_{\{s=\zeta\}}\, \PN(\mathrm ds,\mathrm dz) =\E\int_0^T\!\int_{Z_j}\mathds1_{\{s=\zeta\}}\, \mu(\mathrm dz)\dd s=0.
$$ Taking a countable union of the exceptional null sets over $j\in\N$ and positive integer horizons $T$ shows that the Poisson random measure has no atom at $\zeta$.
The jump identity on closed localizations of $u$ yields $u(\zeta)=u(\zeta-)$ on $B$, proving \eqref{eq:4.4_crit_term}.
For the construction restarted at $\zeta$, set $$ B_{T,M}:=\left\{\zeta\le T,\quad \sup_{t<\zeta}\|v(t)\|_3\le M,\quad \int_0^\zeta\|v(t)\|_5^5\dd t\le M^5\right\}.
$$ These $\mathcal F_\zeta$-events exhaust $B$ over positive integers $T,M$.
The terminal state has norm at most $M$ on $B_{T,M}$.
Extend it by zero off that event to obtain $$ \xi_{T,M}:=\mathds1_{B_{T,M}}v_\zeta \in L^3(\Omega,\mathcal F_{\zeta\wedge T};E).
$$ The $\mathcal F_{\zeta\wedge T}$-measurability follows from $B_{T,M}\subseteq\{\zeta\le T\}$.
The first local construction, with initial time $\zeta\wedge T$ and deterministic terminal time $T+1$, gives a strictly positive continuation on $B_{T,M}$.
\cref{lem:stopping-time-restart} supplies the corresponding stopping-time restart identity, and stopping-time uniqueness identifies it with $u$ on the common interval.
The maximal process already gives the extension on the original stochastic basis.
At a general stopping time in the interval of existence, the restart state is the post-jump value $u(\sigma)$. 
Equality with the left limit above uses the fact that $\zeta$ is announced by $(\kappa_n)$.
\end{proof}

\begin{proof}
[Proof of \cref{thm:main-local}] For data satisfying \eqref{eq:initial}, the result follows from \cref{prop:critical-uniqueness-stability,thm:maximal-critical-solution}.
Indeed, the bounds for the cubic term in \eqref{eq:3.3_tame1} and \eqref{eq:3.3_tame2} follow from \cref{hyp:T1}.
After applying the cumulative cutoff, the Lipschitz constant of the deterministic part is bounded by $C_{\nu,\vartheta}(\varepsilon+C_g\varepsilon^2)$.
The stochastic contraction and the stopping-time restart use \cref{hyp:L}.
These facts also apply when $g_N\equiv0$.

Now let $u_0\in L^0(\Omega,\mathcal F_0;E)$.
For $k\ge1$ put $$ D_k=\{k-1\le\|u_0\|_3<k\},\qquad a_k=\mathds1_{D_k}u_0, \qquad \Sigma_k=\mathds1_{D_k}\Sigma, \qquad \mathcal G_k=\mathds1_{D_k}\mathcal G .
$$ The events $D_k$ form an $\mathcal F_0$-measurable partition up to a null set.
The data $a_k$ are essentially bounded in $E$, and $(\Sigma_k,\mathcal G_k)$ satisfy \cref{hyp:L} with the original bounds on each finite horizon.
On the given stochastic basis, let $(u^k,\{\tau_n^k\},\zeta_k)$ be the maximal solution for $(a_k,\Sigma_k,\mathcal G_k)$ supplied by \cref{thm:maximal-critical-solution}.
On $D_k$ set $$ \tau_*=\zeta_k,\qquad \tau_n=n\wedge\tau_n^k,\qquad u(t)=u^k(t)\quad(0\le t<\zeta_k).
$$ An $\mathcal F_0$-measurable choice among countably many stopping times is a stopping time.
Thus $0<\tau_n\le n$, $\tau_n\uparrow\tau_*$, and $\tau_n<\tau_*$ on $\{\tau_*<\infty\}$, almost surely.
Each stopped process is adapted, has c\`adl\`ag $E$-valued paths, and has finite space--time $L^5$ norm on its closed solution interval.

To verify the equation, fix $n$.
On $D_k$ the initial term and both deterministic convolutions in \eqref{eq:2.4_stop_Con} agree with those for $u^k$ stopped at $n\wedge\tau_n^k$.
The same holds for the stochastic convolutions by \eqref{eq:3.6_restart_locality} at initial time zero, applied after localization.
Indeed, on each bounded closed interval \eqref{eq:local-growth} and the pathwise $E$ bound make the two quadratic integrability quantities and the cubic Poisson quantity finite almost surely.
The $\mathcal F_0$ indicators therefore commute with both localized integrals.
Choosing the c\`adl\`ag convolution representative on each $D_k$ gives \eqref{eq:closed-mild-equality} simultaneously in time, including at $\tau_n$.
The exceptional sets may be combined over $k,n$ and integer horizons.

For comparison, let $v$ be any other critical solution with the same datum, coefficients and noises.
On each of its closed intervals, $\mathds1_{D_k}v$ solves the equation with data $(a_k,\Sigma_k,\mathcal G_k)$.
For example, $$ \Sigma_k(t,\mathds1_{D_k}v(t-)) =\mathds1_{D_k}\Sigma(t,v(t-)),\qquad \mathcal G_k(t,\mathds1_{D_k}v(t-),z) =\mathds1_{D_k}\mathcal G(t,v(t-),z).
$$ The tensor product satisfies the same indicator identity, and $\TN(0)=0$.
Thus \cref{prop:critical-uniqueness-stability} applied to two such restricted solutions proves uniqueness on every common closed interval, including its terminal value.
If $v$ has lifetime $\eta$, \cref{thm:maximal-critical-solution} applied to $\mathds1_{D_k}v$ gives $\eta\le\zeta_k$ and identifies it with $u^k$ up to $\eta$.
Restricting to $D_k$ and taking their union proves maximality and pathwise uniqueness for the original datum.
Finally, on each $D_k$ the lifetime and the space--time integral are those of $u^k$.
The finite-moment criterion \eqref{eq:4.4_crit_space} and countability give the asserted continuation criterion.

Let $v$ be any local solution with the same datum, coefficients and noises and lifetime $\zeta$.
On $\{\zeta<\infty,\ \int_0^\zeta\|v(t)\|_5^5\,\dd t<\infty\}$, the restriction property and the criterion imply $\zeta<\tau_*$ almost surely.
Hence $v(t)=u(t)$ for $t<\zeta$ and its $E$-valued left limit is $u(\zeta-)$.
The announcing times in \cref{def:critical-solution} make the graph of the finite lifetime predictable.
Compensation on an increasing exhaustion of $Z$ by finite-measure sets, as in the proof of \eqref{eq:4.4_crit_term}, shows that $\PN$ has no atom at $\zeta$.
The jump formula on closed restrictions of $u$ gives $u(\zeta-)=u(\zeta)$.
The continuation is the restriction of $u$ after $\zeta$.
\end{proof}
 \cref{Section 5} and \cref{Section 6} use \cref{hyp:T} and \cref{hyp:G} to derive the remainder estimates needed for global continuation.

\begin{remark}
[Stopped weak formulation]\label{rem:stopped-weak-formulation} Every closed local critical mild solution $(u,\tau)$ in the sense of \cref{def:critical-solution} satisfies the following weak formulation.
For every solenoidal $\varphi\in C_c^\infty(\R^3;\R^3)$, almost surely, simultaneously for all $t\ge0$,
\begin{equation*}
\begin{aligned}
\langle u(t\wedge\tau),\varphi\rangle =&\langle u_0,\varphi\rangle +\nu\int_0^t\chi_\tau(s) \langle u(s),\Delta\varphi\rangle\dd s+\int_0^t\chi_\tau(s) \int_{\R^3}(u(s)\otimes u(s)): \nabla\varphi\dd x\dd s \\
&-\int_0^t\chi_\tau(s) \langle\TN(u(s)),\varphi\rangle\dd s +\int_0^t\chi_\tau(s) \langle\Sigma(s,u(s-))\dd W(s),\varphi\rangle \\
&+\int_0^t\int_Z\chi_\tau(s) \langle\mathcal G(s,u(s-),z),\varphi\rangle \,\tPN(\mathrm ds,\mathrm dz).
\end{aligned}
\end{equation*}

To establish this identity, first restrict to $\{\|u_0\|_3\le k\}\in\mathcal F_0$ using \eqref{eq:3.6_restart_locality}, then localize the critical norms and the Wiener and Poisson integrability quantities.
The initial-data events exhaust $\Omega$ up to a null set.
The regularization argument of \cref{lem:stopping-time-restart}, together with the deterministic convolution estimates of \cref{lem:3.2_critical-convection} and \cref{lem:critical-taming}, gives the weak formulation for $Y_\tau$ at deterministic times.
H\"older's inequality justifies passage to the limit in the deterministic terms, while the scalar It\^o isometries and martingale maximal inequalities give convergence of the stochastic integrals.

Using the c\`adl\`ag modifications and a countable dense set of times, the identity extends to all times outside one null set.
Evaluating the weak formulation for $Y_\tau$ at $t\wedge\tau$ and applying \eqref{eq:closed-mild-equality} yields the asserted identity.
Since $\varphi$ is solenoidal, $\mathcal P\varphi=\varphi$; duality for the Leray projection and the tensor convention \eqref{eq:2.1_tensor-convention} give the stated sign of the convection term.

On every finite deterministic interval $[0,T]$, the critical norms on $[0,T\wedge\tau]$ and the corresponding restricted Wiener and Poisson integrability quantities are finite almost surely.
Consequently, the localization times equal $T$ for all sufficiently large indices, almost surely.
Removing the localization and intersecting the resulting full-probability events over integer $T\ge1$ establishes the identity simultaneously for all $t\ge0$.
\end{remark}
\section{Caloric decomposition and finite-energy remainder}
\label{Section 5} The continuation criterion in \cref{Section 4} reduces global existence to bounds for the critical norms on finite time intervals.
Since the initial velocity need not belong to $L^2$, we derive the additional a priori estimates for the remainder after subtracting the heat evolution of the initial datum.
We first establish the energy regularity of this remainder and then obtain localized bounds independent of the solution's critical norms.
Finally, we select deterministic positive times with measurable $H^1_\sigma$-valued remainder states, which will be used in the continuation argument of \cref{Section 6}.

Throughout this section, assume \cref{hyp:T}, \cref{hyp:L}, and \cref{hyp:G}.
Let $(u,\tau)$ be a closed local probabilistically strong critical solution, in the sense of \cref{def:critical-solution}, with $0<\tau\le T<\infty$.
On its closed interval of existence, define $$ h(t):=S(t)u_0,\qquad R(t):=u(t)-h(t),\qquad 0\le t\le\tau.
$$ For a process defined on a deterministic time interval, we use the stopped remainder $R^\tau(t):=u(t\wedge\tau)-h(t\wedge\tau)$.
Thus the heat evolution $h$ is stopped at $\tau$ as well.

The decomposition $u=h+R$ does not by itself imply that $R$ belongs to $L^2_\sigma$.
The following lemma identifies $R$ with deterministic and stochastic convolutions having finite-energy regularity and zero initial value.
This identification justifies the $L^2$ It\^o identity on closed solution intervals, including a possible jump at the terminal stopping time.
\begin{lemma}
[Energy regularity of the caloric remainder] \label{lem:remainder-energy-entrance} Under the assumptions above, the remainder admits an adapted $L^2_\sigma$-valued version, agreeing with $u-h$ simultaneously on $[0,\tau]$, such that, almost surely,
\begin{equation}
\label{eq:5.1_rem_eng}
R\in D([0,\tau];L^2_\sigma) \cap L^2(0,\tau;H^1_\sigma) \cap L^4(0,\tau;L^4_\sigma), \qquad R(0)=0.
\end{equation}
\begin{equation}
\label{eq:5.1_rem_zero}
\lim_{t\downarrow0}\|R(t)\|_2=0\quad\text{almost surely}.
\end{equation}

Define
\begin{align*}
\mathscr D_\tau &:=\int_0^\tau \bigl(\|u(s)\|_4^4+C_g^2\|u(s)\|_{18/5}^6\bigr)\dd s, \\
\mathscr Q_\tau &:=\int_0^\tau\left( \|\Sigma(s,u(s-))\|_{\HS(\U;L^2_\sigma)}^2 +\int_Z\|\mathcal G(s,u(s-),z)\|_2^2\,\mu(\dd z) \right)\dd s.
\end{align*}
Both quantities are finite almost surely.
If their expectations are finite, then
\begin{equation}
\label{eq:5.1_rem_est}
\E\sup_{0\le t\le\tau}\|R(t)\|_2^2 +\nu\E\int_0^\tau\|\nabla R(s)\|_2^2\dd s \le C\bigl(\nu^{-1}\E\mathscr D_\tau +\E\mathscr Q_\tau\bigr),
\end{equation}
where $C$ is independent of $T$, $\tau$, and the solution.
The pathwise assertion \eqref{eq:5.1_rem_eng} holds for every closed critical solution with $u_0\in L^0(\Omega,\mathcal F_0;E)$.
The following identity holds almost surely for every $0\le t\le T$:
\begin{equation}
\label{eq:5.1_rem_idnt}
\begin{aligned}
&\|R^\tau(t)\|_2^2 +2\nu\int_0^{t\wedge\tau}\|\nabla R(s)\|_2^2\dd s =-2\int_0^{t\wedge\tau} \langle\LP\diver(u\otimes u),R\rangle\dd s -2\int_0^{t\wedge\tau} \langle\LP\TN(u),R\rangle\dd s \\
&+\int_0^{t\wedge\tau} \|\Sigma(s,u(s-))\|_{\HS(\U;L^2_\sigma)}^2\dd s +M_R(t) +\int_{(0,t]\times Z}\chi_\tau(s) \|\mathcal G(s,u(s-),z)\|_2^2\,\PN(\dd s,\dd z),
\end{aligned}
\end{equation}
where $\chi_\tau(s)=\mathds1_{\{0<s\le\tau\}}$ and
\begin{equation*}
\begin{aligned}
M_R(t):=&2\int_0^t\chi_\tau(s) (R^\tau(s-),\Sigma(s,u(s-))\dd W(s))_{L^2} \\
&+2\int_0^t\int_Z\chi_\tau(s) (R^\tau(s-),\mathcal G(s,u(s-),z))_{L^2} \,\tPN(\dd s,\dd z).
\end{aligned}
\end{equation*}
Here $M_R$ is a real-valued local martingale.
The convection term is interpreted through $H^{-1}$--$H^1$ duality.
The taming term is interpreted through either $L^{6/5}$--$L^6$ or $L^{4/3}$--$L^4$ duality.
The identity includes the possible jump at $\tau$.
Moreover,
\begin{equation}
\label{eq:5.1_rem_jump}
\tilde\Delta R^\tau(t)=\chi_\tau(t) \int_Z\mathcal G(t,u(t-),z)\,\PN(\{t\},\dd z) \quad\text{in }L^2_\sigma,
\end{equation}
simultaneously for $0<t\le T$.
\end{lemma}

\begin{proof}

For almost every sample path on which the critical norms are finite, define $$ U_\tau:=\sup_{s\le\tau}\|u(s)\|_3, \qquad B_\tau:=\int_0^\tau\|u(s)\|_5^5\dd s.
$$ Both quantities are finite.
Spatial interpolation gives
\begin{equation*}
\|u(s)\|_4 \le\|u(s)\|_3^{3/8}\|u(s)\|_5^{5/8},\qquad \|u(s)\|_{18/5} \le\|u(s)\|_3^{7/12}\|u(s)\|_5^{5/12}.
\end{equation*}
H\"older's inequality in time yields
\begin{align}
\int_0^\tau\|u(s)\|_4^4\dd s &\le T^{1/2}U_\tau^{3/2}B_\tau^{1/2},
\label{eq:5.1_rem_1}\\
\int_0^\tau\|u(s)\|_{18/5}^6\dd s &\le T^{1/2}U_\tau^{7/2}B_\tau^{1/2}.
\end{align}
Extend the deterministic drift terms by zero outside $(0,\tau]$ and define $$ F:=-\chi_\tau u\otimes u\in L^2(0,T;L^2), \qquad f:=-\chi_\tau\TN(u)\in L^2(0,T;L^{6/5}).
$$ \cref{lem:taming} gives
\begin{equation*}
\|F\|_{L^2(0,T;L^2)}^2 +\|f\|_{L^2(0,T;L^{6/5})}^2 \le\mathscr D_\tau<\infty.
\end{equation*}
For $v\in H^1_\sigma$, Leray duality and Sobolev's inequality give
\begin{equation}
\bigl|\langle\LP\diver F+\LP f,v\rangle\bigr| \le\bigl(\|F\|_2+C_S\|f\|_{6/5}\bigr)\|\nabla v\|_2.
\end{equation}

For drift terms $F$ and $f$ belonging to the preceding integrability spaces, the distributional convolution $$ V(t):=\int_0^t S(t-s) [\LP\diver F(s)+\LP f(s)]\dd s $$ has a unique representative in $C([0,T];L^2_\sigma)\cap L^2(0,T;H^1_\sigma)$, with $V(0)=0$, and
\begin{equation}
\label{eq:5.1_det_con}
\sup_{t\le T}\|V(t)\|_2^2 +\nu\int_0^T\|\nabla V(s)\|_2^2\dd s \le\frac C\nu\int_0^T (\|F(s)\|_2^2+\|f(s)\|_{6/5}^2)\dd s.
\end{equation}
Indeed, approximation by smooth functions, testing with $V$, and Young's inequality yield this estimate.
Applying the same estimate to differences gives convergence in the asserted spaces, since the $L^2_tL^2_x$ norm is bounded by $T^{1/2}$ times the uniform-in-time $L^2$ norm.
The limit agrees with the distributional convolution.
Continuity of this construction and its dependence only on the restrictions of $F$ and $f$ to preceding times yield a measurable adapted process for the progressively measurable drift terms above, without requiring $\E\mathscr D_\tau<\infty$.

Define the predictable stochastic integrands $$ \Phi:=\chi_\tau\Sigma(\cdot,u_-), \qquad \Psi:=\chi_\tau\mathcal G(\cdot,u_-,\cdot).
$$ These integrands are predictable, with compatible $L^2_\sigma$-valued versions by \cref{hyp:G}, and
\begin{equation}
\label{eq:5.1_fin_eng}
\mathscr Q_\tau \le C_T^{(1)}T(1+U_\tau^2)<\infty \quad\text{almost surely}.
\end{equation}
Only the $L^2$-valued growth bound in \cref{hyp:G} is used here.
For square-integrable stochastic integrands, the convolution
\begin{equation}
\label{eq:5.1_stoch_con1}
\tilde Z(t):=\int_0^tS(t-s)\Phi(s)\dd W(s) +\int_0^t\int_ZS(t-s)\Psi(s,z)\, \tPN(\dd s,\dd z)
\end{equation}
has an adapted version in $D([0,T];L^2_\sigma)\cap L^2(0,T;H^1_\sigma)$, with $\tilde Z(0)=0$, satisfying
\begin{equation}
\label{eq:5.1_hilbert_eng}
\begin{aligned}
\E\sup_{t\le T}\|\tilde Z(t)\|_2^2 +\nu\E\int_0^T\|\nabla \tilde Z(s)\|_2^2\dd s \le C\E\int_0^T\left( \|\Phi(s)\|_{\HS(\U;L^2_\sigma)}^2 +\int_Z\|\Psi(s,z)\|_2^2\,\mu(\dd z) \right)\dd s.
\end{aligned}
\end{equation}
To establish this estimate, first consider bounded elementary predictable stochastic integrands with smooth compactly supported solenoidal values, finitely many Wiener coordinates, and finitely many measurable subsets of $Z$ having finite $\mu$-measure.
On the corresponding finite mark space, apply \cite[Theorem 2.1, $p=2$, $M=3$]{MR4165650} to obtain

$$\|\tilde Z(t)\|_2^2+2\nu\int_0^t\|\nabla \tilde Z\|_2^2\dd s =M_{\tilde Z}(t)+\int_0^t\|\Phi\|_{\HS(\U;L^2_\sigma)}^2\dd s+\int_{(0,t]\times Z}\|\Psi\|_2^2\,\PN(\dd s,\dd z),$$

where the stochastic integrands of $M_{\tilde Z}$ are the linear expressions $2(\tilde Z_-,\Phi\,\cdot)_{L^2}$ and $2(\tilde Z_-,\Psi)_{L^2}$.
The Burkholder--Davis--Gundy inequality for optional quadratic variation, followed by Young's inequality and compensation, gives $$ \E\sup_{t\le T}|M_{\tilde Z}(t)| \le\tfrac12\E\sup_{t\le T}\|\tilde Z(t)\|_2^2 +C\E\int_0^T\left(\|\Phi\|_{\HS(\U;L^2_\sigma)}^2 +\int_Z\|\Psi\|_2^2\,\mu(\dd z)\right)\dd s.
$$ This proves \eqref{eq:5.1_hilbert_eng}.
The corresponding estimate for differences extends the construction by density to all square-integrable stochastic integrands: a subsequence converges almost surely uniformly in $L^2$ and in $L^2(0,T;H^1)$, and the It\^o isometries identify its limit with \eqref{eq:5.1_stoch_con1} at each time.
The nonnegative jump-square contribution is increasing, so no fourth moment of the Poisson integrand is required.

For the stochastic integrands satisfying the pathwise integrability condition \eqref{eq:5.1_fin_eng}, introduce the stopping times $$ \rho_m:=T\wedge\inf\{t:A(t)\ge m\},\qquad A(t):=\int_0^t\left(\|\Phi\|_{\HS(\U;L^2_\sigma)}^2 +\int_Z\|\Psi\|_2^2\,\mu(\dd z)\right)\dd s.
$$ The stochastic integrands restricted to $(0,\rho_m]$ have expected quadratic integrability quantities bounded by $m$.
The stopped difference estimate gives equality of their stochastic convolutions up to the smaller stopping time.
Since $A$ is continuous and $A(T)=\mathscr Q_\tau<\infty$, one has $\rho_m=T$ for all sufficiently large $m$ on each sample path satisfying this integrability condition.
This defines $\tilde Z$ with the asserted energy regularity.

The adapted process $\widehat R:=V+\tilde Z$ satisfies $$ \widehat R\in D([0,T];L^2_\sigma)\cap L^2(0,T;H^1_\sigma), \qquad \widehat R(0)=0.
$$ Its deterministic terms agree in distributions with the corresponding terms of $Y_\tau-h$, where $Y_\tau$ is defined in \eqref{eq:2.4_stop_Con}.
After pairing with a smooth compactly supported test function, the compatibility condition in \cref{hyp:G} identifies the $L^2$- and $L^3$-valued stochastic integrands and, consequently, their scalar stochastic integrals, initially under localization and subsequently without localization.
Thus the two processes agree in distributions at each deterministic time, without requiring an embedding between $L^2$ and $L^3$.
A countable separating family of test functions, rational times, and the endpoint $T$, together with the c\`adl\`ag property, gives simultaneous equality on $[0,T]$.
By \eqref{eq:closed-mild-equality}, $\widehat R=u-h$ on $[0,\tau]$.
Restriction to $[0,\tau]$ and stopping at $\tau$ give the required version; right continuity at zero proves \eqref{eq:5.1_rem_zero}.

The remaining $L^4$ regularity follows from \eqref{eq:5.1_rem_1} and \cref{lem:heat-extension}: $$ \int_0^\tau\|R\|_4^4\dd s \le8\int_0^\tau\|u\|_4^4\dd s +8T^{1/2}\|h\|_{L^8(0,T;L^4)}^4<\infty.
$$ When $\E\mathscr D_\tau$ and $\E\mathscr Q_\tau$ are finite, \eqref{eq:5.1_det_con} and \eqref{eq:5.1_hilbert_eng} imply \eqref{eq:5.1_rem_est}.

Subtracting the stopped heat equation for $h$ from \cref{rem:stopped-weak-formulation}, and using the energy regularity established above, gives in $H^{-1}$ $$ \dd R^\tau=\chi_\tau\bigl[\nu\Delta R-\LP\diver(u\otimes u) -\LP\TN(u)\bigr]\dd t +\Phi\dd W+\int_Z\Psi(t,z)\,\tPN(\dd t,\dd z).
$$ Apply $J_\delta=e^{\delta\Delta}$ and the Hilbert-space It\^o formula.
The regularized drift belongs locally to $L^2$.
As $\delta\downarrow0$, $J_\delta R^\tau\to R^\tau$ uniformly in time in $L^2$, since the range of an $L^2$-valued c\`adl\`ag path is relatively compact.
Convergence of the gradients in $L^2(0,\tau;L^2)$, of $J_\delta F$ in $L^2_tL^2_x$, and of $J_\delta f$ in $L^2_tL^{6/5}_x$ permits passage to the limit in the drift and viscous terms.
Localize the energy and quadratic noise integrals, together with the running $L^2$ norm.
At an exit caused by a jump, the norm bound is used for the left limits appearing in the stochastic integrands; the terminal value of $R^\tau$ and the jump at that time remain included.
Quadratic-variation estimates then yield uniform convergence in probability of the linear stochastic integrals; the It\^o correction and nonnegative jump terms converge by $L^2$ convergence of the stochastic integrands and compensation.
Removing the localization proves \eqref{eq:5.1_rem_idnt} simultaneously in time.
The scalar stochastic integrands are locally square integrable, with predictable quadratic variation bounded by $$ 4\sup_{t\le\tau}\|R(t)\|_2^2\,\mathscr Q_\tau<\infty.
$$ The alternative duality interpretation of the taming term follows from $u,R\in L^4((0,\tau)\times\R^3)$.
Finally, testing the stopped weak equation against a countable separating family gives \eqref{eq:5.1_rem_jump}.
The indicator $\chi_\tau$ includes the terminal jump in both identities.
\end{proof}

To estimate the convective contribution to the remainder energy identity, we use the divergence-free cancellations to express it in terms of $h$ and spatial derivatives of $R$.
The following bound separates a term that can be absorbed by the viscous dissipation from terms controlled by $\|h\|_5^5\|R\|_2^2$ and $\|h\|_4^4$.

\begin{lemma}
[Mixed convection estimates] \label{lem:mixed-remainder-convection} Let $h\in L^3_\sigma\cap L^4\cap L^5$ and $R\in H^1_\sigma$, and set $u=h+R$.
Then $\LP\diver(u\otimes u)\in H^{-1}$ and
\begin{equation}
\label{eq:5.2_rem}
\langle \LP\diver(u\otimes u),R \rangle_{H^{-1},H^1} = -\int_{\R^3} ((R\cdot\nabla)R)\cdot h\dd x-\int_{\R^3} (h\otimes h):\nabla R\dd x.
\end{equation}
In particular,
\begin{equation}
\label{eq:5.2_rem-bound}
\left| \langle \LP\diver(u\otimes u),R \rangle_{H^{-1},H^1} \right| \le \|h\|_3\|R\|_6\|\nabla R\|_2+\|h\|_4^2\|\nabla R\|_2.
\end{equation}
Furthermore, for every $\eta>0$,
\begin{equation}
\label{eq:5.2_mixed_remain}
\left| \langle \LP\diver(u\otimes u),R \rangle_{H^{-1},H^1} \right| \le \eta\|\nabla R\|_2^2 +C\eta^{-4}\|h\|_5^5\|R\|_2^2+C\eta^{-1}\|h\|_4^4.
\end{equation}
The constant $C$ depends only on the spatial Sobolev constant.
For the solution of \cref{lem:remainder-energy-entrance}, the convection pairing is integrable in time, and the corresponding time-integrated estimate holds.
The estimates require neither an $L^2$ norm nor a spatial derivative of $h$.
\end{lemma}

\begin{proof}
Since $R\in H^1_\sigma\subset L^4$ and $h\in L^4$, we have $u\otimes u\in L^2$.
The boundedness of $\LP$ on $H^{-1}$ therefore makes the pairing with $R$ well defined.

Using \eqref{eq:negative-duality} and distributional integration by parts, we obtain $$
\begin{aligned}
\langle \LP\diver(u\otimes u),R \rangle_{H^{-1},H^1} &= -\int_{\R^3}(u\otimes u):\nabla R\dd x \\
&= -\int_{\R^3} (R\otimes R+h\otimes R +R\otimes h+h\otimes h) :\nabla R\dd x.
\end{aligned}
$$ By \eqref{eq:tensor-product} and \eqref{eq:skew}, $$ \int_{\R^3} (R\otimes R):\nabla R\dd x = B(R,R,R)=0 ; \qquad \int_{\R^3} (h\otimes R):\nabla R\dd x = B(h,R,R)=0.
$$ Moreover, the tensor convention \eqref{eq:2.1_tensor-convention} gives $$ (R\otimes h):\nabla R = ((R\cdot\nabla)R)\cdot h.
$$ Consequently, \eqref{eq:5.2_rem} follows.
H\"older's inequality immediately gives \eqref{eq:5.2_rem-bound}.
For the stronger estimate, H\"older's inequality, interpolation between $L^2$ and $L^6$, and Young's inequality give $$ \left| \int_{\R^3} ((R\cdot\nabla)R)\cdot h\dd x \right|\le \|h\|_5\|R\|_{10/3}\|\nabla R\|_2\le C\|h\|_5\|R\|_2^{2/5} \|\nabla R\|_2^{8/5}\le \frac{\eta}{2}\|\nabla R\|_2^2 +C\eta^{-4}\|h\|_5^5\|R\|_2^2.
$$ Similarly, $$ \left| \int_{\R^3} (h\otimes h):\nabla R\dd x \right| \le \|h\|_4^2\|\nabla R\|_2\le \frac{\eta}{2}\|\nabla R\|_2^2 +C\eta^{-1}\|h\|_4^4.
$$ Combining these inequalities proves \eqref{eq:5.2_mixed_remain}.

For the time-dependent assertion, let $(u,\tau)$ be the closed local critical solution considered in \cref{lem:remainder-energy-entrance}, with $\tau\le T<\infty$.
By \cref{lem:heat-extension}, almost surely, $$ \int_0^T\|h(s)\|_5^5\dd s<\infty, \qquad \int_0^T\|h(s)\|_4^4\dd s \le T^{1/2} \|h\|_{L^8(0,T;L^4)}^4<\infty.
$$ Furthermore, \cref{lem:remainder-energy-entrance} gives $$ \sup_{0\le s\le\tau}\|R(s)\|_2^2 + \int_0^\tau\|\nabla R(s)\|_2^2\dd s <\infty \qquad\text{almost surely}.
$$ Thus the right-hand side of \eqref{eq:5.2_mixed_remain} is integrable on $(0,\tau)$ almost surely.
The asserted time integrability and integrated estimate follow.
\end{proof}

\begin{lemma}
[Shifted taming coercivity] \label{lem:shifted-taming-coercivity} Assume \cref{hyp:T}.
Let $h\in L^4$, let $R\in L^2_\sigma\cap L^4$, and set $u=h+R$.
The $L^{4/3}$--$L^4$ duality pairing satisfies
\begin{equation}
\label{eq:5.3_shifted_tame}
\langle\LP\TN(u),R\rangle =\int_{\R^3}\TN(h+R)\cdot R\dd x.
\end{equation}
For every $0<\alpha<(8\nu)^{-1}$,
\begin{equation}
\label{eq:5.3_shifted-tame2}
\begin{aligned}
\langle\LP\TN(h+R),R\rangle \ge&\left(\frac1{8\nu}-\alpha\right)\|R\|_4^4 -\frac{c_{N,\nu}}2\|R\|_2^2-\frac{3}{4^{4/3}}\alpha^{-1/3} C_g^{4/3}\|h\|_4^4.
\end{aligned}
\end{equation}
In particular, with $C_*:=\frac{3}{2^{4/3}}$,
\begin{equation}
\label{eq:5.3_shifted-tame3}
\langle\LP\TN(h+R),R\rangle \ge\frac1{16\nu}\|R\|_4^4 -\frac{c_{N,\nu}}2\|R\|_2^2 -C_*\nu^{1/3}C_g^{4/3}\|h\|_4^4.
\end{equation}
For $h$ and the remainder $R$ associated with the solution considered in \cref{lem:remainder-energy-entrance}, these estimates hold almost surely for almost every $t\in(0,\tau)$ and may be integrated over every subinterval of $[0,\tau]$.
\end{lemma}

\begin{proof}
Since $u\in L^4$, \cref{lem:taming} gives $\TN(u)\in L^{4/3}$.
Since $R\in L^4_\sigma$, the Leray duality identity \eqref{eq:leray-duality} yields \eqref{eq:5.3_shifted_tame}.
Applying \eqref{eq:quantitative-monotonicity} pointwise with $a=h+R$ and $b=h$, and integrating over $\R^3$, gives $$ \int_{\R^3} (\TN(h+R)-\TN(h))\cdot R\dd x \ge\frac1{8\nu}\|R\|_4^4 -\frac{c_{N,\nu}}2\|R\|_2^2.
$$ The remaining contribution satisfies $$ \left|\int_{\R^3}\TN(h)\cdot R\dd x\right| \le C_g\int_{\R^3}|h|^3|R|\dd x.
$$ For $x,y\ge0$, the parameter-dependent Young inequality states that $$ xy\le\alpha x^4 +\frac{3}{4^{4/3}}\alpha^{-1/3}y^{4/3}.
$$ Applying this inequality with $x=|R|$ and $y=C_g|h|^3$, followed by integration, proves \eqref{eq:5.3_shifted-tame2}.
The choice $\alpha=(16\nu)^{-1}$ then yields \eqref{eq:5.3_shifted-tame3}.

For the solution of \cref{lem:remainder-energy-entrance}, one has, almost surely, $$ R\in L^\infty(0,\tau;L^2)\cap L^4(0,\tau;L^4).
$$ Moreover, since $\tau\le T<\infty$, \cref{lem:heat-extension} gives $$ \int_0^\tau\|h(s)\|_4^4\dd s \le T^{1/2}\|h\|_{L^8(0,T;L^4)}^4<\infty.
$$ The growth bound and H\"older's inequality imply $$ \int_0^\tau \left|\langle\LP\TN(h+R),R\rangle\right|\dd s \le C_g \|h+R\|_{L^4(0,\tau;L^4)}^3 \|R\|_{L^4(0,\tau;L^4)} <\infty.
$$ Thus all terms are integrable in time almost surely, and the estimates may be integrated over every subinterval of $[0,\tau]$.
\end{proof}

The preliminary estimate \eqref{eq:5.1_rem_est} depends on quantities controlled by the critical norms of $u$.
Combining the two preceding drift estimates with \cref{hyp:G}, the following proposition gives energy bounds localized on events where $\|u_0\|_3$ is bounded, with constants independent of the solution's critical norms and of the auxiliary stopping thresholds.
This uniformity permits passage from closed solution intervals to the maximal interval of existence.
\cref{Section 6} will use these bounds to obtain Sobolev regularity and control the critical continuation quantities.
\begin{proposition}
[Uniform localized remainder energy estimates] \label{prop:uniform-energy} Assume \cref{hyp:T}, \cref{hyp:L}, and \cref{hyp:G}.
Let $u_0$ be an $\mathcal F_0$-measurable $L^3_\sigma$-valued initial datum.
Fix $T>0$, and let $(u,\tau)$ be a closed local critical solution with $0<\tau\le T$.
Set $h(t)=S(t)u_0$ for $t\ge0$ and $R(t)=u(t)-h(t)$ for $0\le t\le\tau$.
Define
\begin{equation}
\label{eq:5.4_rem_eng}
\mathcal E_R(\tau) :=\sup_{0\le t\le\tau}\|R(t)\|_2^2 +\nu\int_0^\tau\|\nabla R(s)\|_2^2\dd s +\nu^{-1}\int_0^\tau\|R(s)\|_4^4\dd s .
\end{equation}
For $K\ge1$, put $\Omega_K=\{\|u_0\|_3\le K\}\in\mathcal F_0$.
There exists $C=C(\nu,C_g,c_{N,\nu},C_T^{(1)})\ge1$ such that, with
\begin{equation}
\label{eq:5.4_loc_eng1}
\Gamma_{T,K}=C e^{CT}\bigl[T(1+K^2)+T^{1/2}K^4\bigr],
\end{equation}
one has
\begin{equation}
\label{eq:5.4_loc_eng2}
\E[\mathds1_A\mathcal E_R(\tau)] \le\Gamma_{T,K}\P(A) \qquad(A\in\mathcal F_0,\ A\subseteq\Omega_K).
\end{equation}
The constant $C$ is independent of $K$, the solution's critical norms, and the choice of stopping times.
The supremum includes $R(\tau)$.

For the maximal critical solution $(u,\{\tau_n\},\tau_*)$, set $R=u-h$ on $[0,\tau_*)$ and write
\begin{equation}
\label{eq:5.4_max_eng}
\begin{aligned}
\mathcal E_R^*(T) := {}&\sup_{\substack{0\le t\le T \\
t<\tau_*}}\|R(t)\|_2^2 +\nu\int_0^T\mathds1_{\{s<\tau_*\}}\|\nabla R(s)\|_2^2\dd s +\nu^{-1}\int_0^T\mathds1_{\{s<\tau_*\}}\|R(s)\|_4^4\dd s .
\end{aligned}
\end{equation}
The integrands are extended by zero outside the interval of existence.
Then
\begin{equation}
\label{eq:5.4_max_eng2}
\E[\mathds1_A\mathcal E_R^*(T)] \le\Gamma_{T,K}\P(A) \qquad(A\in\mathcal F_0,\ A\subseteq\Omega_K).
\end{equation}
If, in addition, \eqref{eq:initial} holds, then for every $T>0$ there is a deterministic constant $C_T$, depending only on $T$ and the structural and coefficient bounds above, such that
\begin{equation}
\label{eq:5.4_fractional_energy}
\E\bigl[(\mathcal E_R^*(T))^{3/4}\bigr] \le C_T\bigl(1+\E\|u_0\|_3^3\bigr).
\end{equation}
Almost surely, $\mathcal E_R^*(T)<\infty$ simultaneously for all finite $T$.
\end{proposition}

\begin{proof}
Write $X=\|R\|_2^2$, $Y=\|\nabla R\|_2^2$, and $\tilde Y=\|R\|_4^4$ on the solution interval.
The quantities $Y,\tilde Y$ are defined almost everywhere in time.
Write $$ q_W=\|\Sigma(\cdot,u_-)\|_{\HS(\U;L^2_\sigma)}^2,\qquad q_J=\int_Z\|\mathcal G(\cdot,u_-,z)\|_2^2\,\mu(\dd z),\qquad q=q_W+q_J .
$$ At a time where the spatial norms are finite, set $H=\|h\|_4$, $r=\|R\|_4$, and $d=\|\nabla R\|_2$.
By \eqref{eq:5.2_rem} and H\"older's inequality, $$ 2|\langle\LP\diver(u\otimes u),R\rangle| \le 2d(Hr+H^2) \le \nu d^2+\frac{r^4}{16\nu}+\frac{18}{\nu}H^4 .
$$ Indeed, $(Hr+H^2)^2\le2H^2r^2+2H^4$ and $2H^2r^2\le r^4/16+16H^4$.
With $\alpha=(32\nu)^{-1}$, \eqref{eq:5.3_shifted-tame2} gives $$ 2\langle\LP\TN(h+R),R\rangle \ge\frac{3}{16\nu}r^4-c_{N,\nu}X-3\nu^{1/3}C_g^{4/3}H^4 .
$$ Inserting both estimates into \eqref{eq:5.1_rem_idnt} gives, in integrated form on $[0,\tau]$,
\begin{equation}
\label{eq:5.4_mart_ineq}
\dd X+\left(\nu Y+\frac{\tilde Y}{8\nu}\right)\dd t \le (aX+b+q_W)\dd t+\dd M_R+\dd J_R,
\end{equation}
where $M_R$ is the local martingale in \eqref{eq:5.1_rem_idnt}, $a=c_{N,\nu}$, and $$
\begin{aligned}
b(t)&=(18/\nu+3\nu^{1/3}C_g^{4/3})\|h(t)\|_4^4,& J_R(t)&=\int_{(0,t]\times Z}\chi_\tau(s) \|\mathcal G(s,u(s-),z)\|_2^2\,\PN(\dd s,\dd z).
\end{aligned}
$$ The heat estimate in \cref{lem:heat-extension} gives $b\in L^1(0,T)$ almost surely.
Since $u$ is c\`adl\`ag, $u(s-)=u(s)$ almost everywhere in time.
By \cref{hyp:G} and $\|R\|_3^2\le C_S\|R\|_2\|\nabla R\|_2$, we have $q\le c_T(1+\|h\|_3^2+X^{1/2}Y^{1/2})$ almost everywhere on $(0,\tau)$.
In particular, for every $\delta>0$,
\begin{equation}
\label{eq:5.4_quad_noise}
q\le\delta Y+c_\delta X+c_T(1+\|h\|_3^2),
\end{equation}
where $c_T,c_\delta$ depend only on the coefficient bound on $[0,T]$, the Sobolev constant, and, for $c_\delta$, on $\delta$.

Set $$
\begin{aligned}
L(t)&=\int_0^{t\wedge\tau}(Y+\tilde Y+q+aX+b)(s)\dd s, \\
\theta_m&=\tau\wedge\inf\{t\in[0,T]: \|R^\tau(t)\|_2^2>m\ \text{or}\ L(t)\ge m\},
\end{aligned}
$$ with $\inf\varnothing=\infty$.
By \cref{lem:remainder-energy-entrance}, $\theta_m$ increases to $\tau$ and equals $\tau$ eventually almost surely.
For $0<s\le\theta_m$, $X(s-)\le m$, and $L(\theta_m)\le m$ by continuity.
For the stopped martingale $M_R^{\theta_m}(t)=M_R(t\wedge\theta_m)$ this gives $$ \E\langle M_R^{\theta_m}\rangle_T \le4m\,\E\int_0^{\theta_m}q(s)\dd s\le4m^2,\qquad \E J_R(\theta_m)=\E\int_0^{\theta_m}q_J(s)\dd s\le m .
$$ Thus $M_R^{\theta_m}$ is square integrable.
Since $X(0)=0$, \eqref{eq:5.4_mart_ineq} also gives $$ \E\sup_{t\le\theta_m}X(t)<\infty .
$$ This supremum includes a possible jump at $\theta_m$.

Fix $\theta=\theta_m$.
For a deterministic $\lambda>0$ put $w(t)=e^{-(a+\lambda)t}$.
since $w$ is deterministic and $C^1$ with $w'(t)=-(a+\lambda)w(t)$, applying It\^o's product rule to $w(t)X(t)$ via the exact identity \eqref{eq:5.1_rem_idnt} on $[0,\theta]$ and substituting the two drift bounds above (using $w>0$) gives, 
$$
\begin{aligned}
w(t)X(t)&+\int_0^t w(s)\left(\nu Y(s)+\frac{\tilde Y(s)}{8\nu} +\lambda X(s)\right)\dd s \\
&\le\int_0^t w(s)(b(s)+q_W(s))\dd s +\int_0^t w(s)\dd M_R(s)+\int_0^t w(s)\dd J_R(s).
\end{aligned}
$$ Write $$
\begin{aligned}
U&=\sup_{t\le\theta}w(t)X(t),& V&=\int_0^\theta w\left(\nu Y+\frac{\tilde Y}{8\nu} +\lambda X\right)\dd s, \\
B&=\int_0^\theta wb\dd s,& Q_W&=\int_0^\theta wq_W\dd s,\qquad J_w=\int_0^\theta w\dd J_R,
\end{aligned}
$$ and $\mathcal M(t)=\int_0^{t\wedge\theta}w\dd M_R$.
Taking the supremum of the energy and evaluating the integral at $\theta$ gives $U+V\le2(B+Q_W+J_w+\sup_{t\le T}|\mathcal M(t)|)$.
For $0<s\le\theta$, continuity of $w$ yields $w(s)X(s-)\le U$.
The continuous and purely discontinuous martingale parts are orthogonal, so their optional quadratic variations give $[\mathcal M]_T\le4U(Q_W+J_w)$.

For $A\in\mathcal F_0$, the Burkholder--Davis--Gundy and Young inequalities applied to $\mathds1_A\mathcal M$, followed by compensation of the predictable integrand in $J_w$, give $$ \E\!\left[\mathds1_A\sup_{t\le T}|\mathcal M(t)|\right] \le\varepsilon\,\E[\mathds1_A U] +C_\varepsilon\,\E\!\left[\mathds1_A\int_0^\theta wq\dd s\right].
$$ Absorbing the first term for fixed sufficiently small $\varepsilon$ yields a numerical constant $C_0\ge1$ such that $$ \E[\mathds1_A(U+V)] \le C_0\,\E\!\left[\mathds1_A\int_0^\theta w(b+q)\dd s\right].
$$ Take $\delta=\nu/(2C_0)$ in \eqref{eq:5.4_quad_noise} and fix $\lambda=1+2C_0c_\delta$.
Absorption of the $Y$ and $X$ terms gives
\begin{equation}
\label{eq:5.4_weighted_rem}
\begin{aligned}
&\E\!\left[\mathds1_A\left\{U+\int_0^\theta w \left(\frac{\nu}{2}Y+\frac{\tilde Y}{8\nu} +\frac{\lambda}{2}X\right)\dd s\right\}\right] \le C\,\E\!\left[\mathds1_A\int_0^\theta w(1+\|h\|_3^2+\|h\|_4^4)\dd s\right].
\end{aligned}
\end{equation}

On $\Omega_K$, heat contraction and \eqref{eq:3.1_heat-critical} give $$ \sup_{s\le T}\|h(s)\|_3\le K,\qquad \int_0^T\|h(s)\|_4^4\dd s \le T^{1/2}\|h\|_{L^8(0,T;L^4)}^4 \le C\nu^{-1/2}T^{1/2}K^4 .
$$ The deterministic bound $w(t)\ge e^{-(a+\lambda)T}$ and \eqref{eq:5.4_weighted_rem} imply $\E[\mathds1_A\mathcal E_R(\theta_m)]\le\Gamma_{T,K}\P(A)$ for all $A\in\mathcal F_0$ with $A\subseteq\Omega_K$, after increasing the structural constant in \eqref{eq:5.4_loc_eng1}.
Monotone convergence as $m\to\infty$ proves \eqref{eq:5.4_loc_eng2}.

Apply \eqref{eq:5.4_loc_eng2} to the closed solutions stopped at $T\wedge\tau_n$.
Their $L^2$-valued remainders agree on overlaps because each agrees with $u-h$ in distributions.
Their energies increase to $\mathcal E_R^*(T)$.
On $\{T<\tau_*\}$ the intervals eventually include $T$.
On $\{\tau_*\le T\}$ they exhaust $[0,\tau_*)$ by strict announcement at a finite lifetime.
Monotone convergence gives \eqref{eq:5.4_max_eng2}.

Finally, put $Z_T=\mathcal E_R^*(T)$ and decompose $\Omega$, up to a null set, into the $\mathcal F_0$-events $$ D_0=\{\|u_0\|_3\le1\},\qquad D_j=\{2^{j-1}<\|u_0\|_3\le2^j\}\quad(j\ge1).$$ 
For every $j\ge0$, $D_j\subseteq\Omega_{2^j}$.
H\"older's inequality and \eqref{eq:5.4_max_eng2} yield $$ \E[\mathds1_{D_j} Z_T^{3/4}] \le \bigl(\E[\mathds1_{D_j} Z_T]\bigr)^{3/4}\P(D_j)^{1/4} \le \Gamma_{T,2^j}^{3/4}\P(D_j) \le C_T(1+2^{3j})\P(D_j).
$$ For $j\ge1$, $2^{3j}\le 8\|u_0\|_3^3$ on $D_j$.
Under \eqref{eq:initial} the upper bounds are summable in $j$, which proves \eqref{eq:5.4_fractional_energy}.
For each fixed $T$, \eqref{eq:5.4_max_eng2} also gives $Z_T<\infty$ almost surely on each $D_j$.
Applying this observation at positive integer horizons and using monotonicity of $T\mapsto\mathcal E_R^*(T)$ gives almost sure finiteness simultaneously for all finite $T$.
\end{proof}

The energy estimate yields $H^1_\sigma$ regularity for almost every time, whereas the Sobolev equation in \cref{Section 6} must be initialized with an $\mathcal F_s$-measurable state.
The next lemma supplies a deterministic sequence of admissible starting times decreasing to zero, together with localized second-moment bounds.

\begin{lemma}
[Deterministic positive-time $H^1$ regularity of the remainder] \label{lem:admissible-h1} Under the assumptions of \cref{prop:uniform-energy}, let $(u,\{\tau_n\},\tau_*)$ be the maximal critical solution, and let $R=u-S(\cdot)u_0$ on its interval of existence.
There exists a deterministic set $G\subset(0,\infty)$ of full Lebesgue measure such that, for every $s\in G$, the random variable
\begin{equation*}
\xi_s:=
\begin{cases}
R(s),&s<\tau_*, \\
0,&s\ge\tau_*
\end{cases}
\end{equation*}
admits a strongly $\mathcal F_s$-measurable $H^1_\sigma$-valued version.
For every $K\in\N$,
\begin{equation}
\label{eq:5.5_adm_h1}
\E[\mathds1_{\Omega_K}\|\xi_s\|_{H^1}^2]<\infty.
\end{equation}

There exist deterministic times
\begin{equation}
\label{eq:5.5_det_h1}
s_j\in G\cap(2^{-j-1},2^{-j}),\qquad j\ge1.
\end{equation}
They satisfy $s_j\downarrow0$, and almost surely $s_j<\tau_*$ for all sufficiently large $j$.
For $j,K,M\in\N$, define
\begin{equation}
\label{eq:5.5_h1_evts}
A_{j,K,M}:=\{s_j<\tau_*\}\cap\Omega_K \cap\{\|\xi_{s_j}\|_{H^1}\le M\}.
\end{equation}
Then $A_{j,K,M}\in\mathcal F_{s_j}$ and $$ \mathds1_{A_{j,K,M}}\xi_{s_j} \in L^\infty(\Omega,\mathcal F_{s_j};H^1_\sigma), \qquad \bigcup_{K,M\ge1}A_{j,K,M} =\{s_j<\tau_*\} \quad\text{almost surely}.
$$ Since the times $s_j$ are deterministic, they are stopping times, and the corresponding bounded $H^1_\sigma$-valued initial states are $\mathcal F_{s_j}$-measurable.
Moreover, for $0<a<b\le T$ and $K\in\N$, there exists a deterministic $s\in G\cap(a,b)$ such that
\begin{equation}
\label{eq:5.5_h1_2}
\E[\mathds1_{\Omega_K}\|\xi_s\|_{H^1}^2] \le\frac{2(T+\nu^{-1})\Gamma_{T,K}}{b-a} \P(\Omega_K).
\end{equation}
\end{lemma}

\begin{proof}
The consistent remainders associated with the closed stopping times define the progressively measurable $L^2_\sigma$-valued process $$ \widetilde R(t):=\lim_{n\to\infty} \mathds1_{\{t\le\tau_n\}}R^{\tau_n}(t) =\mathds1_{\{t<\tau_*\}}R(t), $$ where the last expression is defined to be zero outside the interval of existence.
Each process in the approximating sequence is progressively measurable.
The inequalities $\tau_n<\tau_*$ on $\{\tau_*<\infty\}$ imply that $\widetilde R$ vanishes at and after $\tau_*$ without requiring the existence of a terminal limit of $R$.

Define $\tilde h(v)=\|v\|_{H^1}^2$ for $v\in H^1_\sigma$ and $\tilde h(v)=\infty$ otherwise.
Monotone convergence in the Fourier representation gives
\begin{equation*}
\tilde h(v)=\|v\|_2^2 +\sup_{\ell\ge1}\|\nabla e^{\Delta/\ell}v\|_2^2,
\end{equation*}
where each term is continuous on $L^2_\sigma$.
Consequently, $\tilde h$ is Borel measurable, and $H^1_\sigma$ is a Borel subset of $L^2_\sigma$.
\cref{prop:uniform-energy} implies
\begin{equation}
\label{eq:5.5_h1_avg}
\int_0^T\E[\mathds1_{\Omega_K}\tilde h(\widetilde R(s))]\dd s \le(T+\nu^{-1})\Gamma_{T,K}\P(\Omega_K)<\infty.
\end{equation}
Indeed, the two time integrals defining the $H^1$ norm are bounded by $T\mathcal E_R^*(T)$ and $\nu^{-1}\mathcal E_R^*(T)$, respectively, and $\tilde h(0)=0$.

Tonelli's theorem, followed by exclusion of the Lebesgue-null sets corresponding to integer $T,K$, gives a deterministic set $G$ of full Lebesgue measure on which the localized expectations are finite for every $K$.
Since $\bigcup_K\Omega_K$ has probability one, $\widetilde R(s)\in H^1_\sigma$ almost surely for each fixed $s\in G$.
For any countable selection of such times, the corresponding probability-one events can be intersected to obtain a common probability-one event.

The operators $e^{\Delta/\ell}:L^2_\sigma\to H^1_\sigma$ are continuous and converge in $H^1$ to the identity on the Borel subset $H^1_\sigma$.
Consequently, the map taking the value $v\in H^1_\sigma$ on this subset and zero on its complement is Borel measurable.
Composing this map with $\widetilde R(s)$ gives the strongly $\mathcal F_s$-measurable version of $\xi_s$ and establishes \eqref{eq:5.5_adm_h1}.

Choose $s_j$ as in \eqref{eq:5.5_det_h1}.
The disjoint, consecutively ordered intervals imply that the sequence decreases strictly to zero.
Since $\tau_*>0$ almost surely, $s_j<\tau_*$ for all sufficiently large $j$ almost surely.
The events defined in \eqref{eq:5.5_h1_evts} belong to $\mathcal F_{s_j}$.
The asserted boundedness and the identity for their union follow from the almost sure finiteness of $\|\xi_{s_j}\|_{H^1}$ and $\|u_0\|_3$.

Finally, restrict \eqref{eq:5.5_h1_avg} to $(a,b)$.
If its upper bound is strictly positive, the integrand does not exceed twice that bound divided by $b-a$ on a set of positive Lebesgue measure.
If the upper bound is zero, the integrand vanishes almost everywhere.
In either case, the corresponding set intersects $G$, proving \eqref{eq:5.5_h1_2}.
\end{proof}

For each deterministic $\delta>0$, choose a time $s_j<\delta$ from \cref{lem:admissible-h1}.
The events $A_{j,K,M}$ exhaust $\{s_j<\tau_*\}$, which contains $\{\delta<\tau_*\}$, and yield bounded, $\mathcal F_{s_j}$-measurable $H^1_\sigma$-valued initial states for the shifted equation.
\section{Positive-time regularity and global continuation}
\label{Section 6} To apply the continuation criterion of \cref{thm:main-local}, we strengthen the finite-energy estimates of \cref{Section 5} to positive-time Sobolev bounds for the remainder.
We construct the shifted evolution from $H^1_\sigma$-valued states and derive estimates that exclude finite-time explosion in that class.
The deterministic starting times of \cref{lem:admissible-h1} and critical pathwise uniqueness identify this evolution with the remainder of the maximal critical solution.
Sobolev embedding and the heat estimates then yield the critical bounds required for global continuation.

Throughout this section, we assume \cref{hyp:T}, \cref{hyp:L}, and \cref{hyp:G}.
For the global theorem we first prove the conclusions under \eqref{eq:initial} and then pass to general measurable data by localization on $\mathcal F_0$.
For a critical solution $u$ on its interval of existence, define $$ h(t):=S(t)u_0,\qquad R(t):=u(t)-h(t),\qquad \Omega_K:=\{\|u_0\|_3\le K\}.
$$ The $L^2$-based Sobolev regularity considered below concerns $R$ and no $L^2$ integrability assumption is imposed on $h$.
For a fixed solenoidal vector field $h$, define
\begin{equation*}
\begin{aligned}
F_h(v)&:=-\LP\bigl[((h+v)\cdot\nabla)(h+v)+\TN(h+v)\bigr], \\
\mathcal A_h(v)&:=\nu\Delta v+F_h(v), \qquad Q_h(v):=\int_{\R^3}|h+v|^2|\nabla v|^2\dd x.
\end{aligned}
\end{equation*}
For $v\in H^1_\sigma$, the convection term is understood in the distributional sense.
For $v\in H^2_\sigma$ and the positive-time caloric fields considered below, $\mathcal A_h(v)\in L^2_\sigma$.

For $v\in H^2_\sigma$, define
\begin{equation}
\label{eq:6.1_h1_diss}
\mathcal D_h(v):= \nu\|\nabla v\|_2^2+\frac{\nu}{2}\|\Delta v\|_2^2 +\frac1{8\nu}\bigl(\|v\|_4^4+Q_h(v)\bigr).
\end{equation}
Unless specified otherwise, constants in this section may depend on $\nu,C_g,c_{N,\nu}$ and the coefficient bounds on the deterministic time interval under consideration.

Fix a deterministic $s>0$ and a strongly $\mathcal F_s$-measurable $H^1_\sigma$-valued random variable $\xi$, finite almost surely.
Consider the shifted equation for $ t\ge s,$
\begin{equation}
\label{eq:6.1_shifted}
\left\{
\begin{aligned}
\dd V(t)=&\mathcal A_{h(t)}(V(t))\dd t +\Sigma(t,h(t)+V(t-))\dd W(t)+\int_Z\mathcal G(t,h(t)+V(t-),z)\, \tPN(\dd t,\dd z), \\
V(s)=&\xi .
\end{aligned}
\right.
\end{equation}
Set $\eta=h(s)+\xi$.
For the localized assertion below, let $A\in\mathcal F_s$ satisfy $$ A\subset\{\|u_0\|_3\le K,\ \|\xi\|_{H^1}\le M\} $$ for some finite $K,M$.
Define, for $t>s$, $$ \eta_A:=\mathds1_A\eta,\qquad \Sigma_A(t,v):=\mathds1_A\Sigma(t,v),\qquad \mathcal G_A(t,v,z):=\mathds1_A\mathcal G(t,v,z).
$$

We first start with the proof of the local $H^1_\sigma$ solvability of the shifted \cref{eq:6.1_shifted}, together with its explosion alternative and jump identity.

\begin{lemma}
[Local $H^1$ well-posedness for the shifted equation] \label{lem:shifted-local-h1} Under \cref{hyp:T}, \cref{hyp:L}, and \cref{hyp:G}, the \cref{eq:6.1_shifted} admits a pathwise unique maximal adapted solution with lifetime $\rho>s$ almost surely.
On every bounded closed localizing interval $[s,\theta]$ with $s<\theta<\rho$ almost surely, it satisfies
\begin{equation*}
V\in D([s,\theta];H^1_\sigma) \cap L^2(s,\theta;H^2_\sigma) \quad\text{almost surely}.
\end{equation*}
The solution satisfies the heat-semigroup mild formulation in $H^1_\sigma$ and the time-integrated formulation in $L^2_\sigma$, with the stochastic integrals defined by localization.
Moreover,
\begin{equation}
\label{eq:6.1_shifted_h1-1}
\{\rho<\infty\} \subset \left\{\sup_{s\le t<\rho}\|V(t)\|_{H^1}=\infty\right\} \quad\text{up to a null set},
\end{equation}
and the jump identity
\begin{equation}
\label{eq:6.1_h1_jump}
\tilde\Delta V(t) =\int_Z\mathcal G(t,h(t)+V(t-),z)\, \PN(\{t\},\dd z)
\end{equation}
holds in $H^1_\sigma$, simultaneously for $s<t\le\theta$, including the terminal value.

Set $U=h+V$.
On every such interval, $U$ belongs to $$ U\in D([s,\theta];E)\cap L^5(s,\theta;L^5_\sigma) $$ and satisfies the critical mild equation with initial time $s$ and datum $\eta$.
For every bounded stopping time $\theta$ with $s<\theta<\rho$ almost surely, define $$ U^\theta(t):=U(t\wedge\theta),\qquad t\ge s.
$$ Then $(U^\theta,\theta)$ is a closed local critical mild solution in the sense of \cref{def:critical-solution}, with initial time $s$ and datum $\eta\in L^0(\Omega,\mathcal F_s;E)$.
For every event $A$ specified above, $\eta_A\in L^3(\Omega,\mathcal F_s;E)$, and the process $$ U_A^\theta(t):=\mathds1_AU^\theta(t),\qquad t\ge s, $$ together with $\theta$, is a closed local critical mild solution with initial time $s$, datum $\eta_A$, and coefficients $\Sigma_A,\mathcal G_A$.
This process is stopped at $\theta$ and vanishes on $A^c$.
The $H^1_\sigma$-valued construction and pathwise uniqueness remain valid under localization on any event in $\mathcal F_s$.
\end{lemma}

\begin{proof}
\emph{Step 1:} First assume $\|u_0\|_3\le K$ almost surely and $\xi\in L^2(\Omega,\mathcal F_s;H^1_\sigma)$, and fix $T>s$.
\cref{lem:heat-extension} gives a bound for
\begin{equation}
\label{eq:6.1_h_bds}
\sup_{s\le t\le T} \left( \|h(t)\|_3+\|h(t)\|_6+\|h(t)\|_\infty +\|\nabla h(t)\|_\infty +\bigl\||h(t)|\,|\nabla h(t)|\bigr\|_2 \right)
\end{equation}
in terms of a deterministic constant depending on $s,T,K,\nu$.
Denote by $H^\alpha_\sigma$ the solenoidal inhomogeneous Bessel-potential space.
The embedding $H^{1/2}\hookrightarrow L^3$ gives $L^{3/2}\hookrightarrow H^{-1/2}$ by duality.
For $v\in H^1_\sigma$, $$ \|(v\cdot\nabla)v\|_{3/2} \le\|v\|_6\|\nabla v\|_2\le C\|v\|_{H^1}^2.
$$ The remaining convective terms, $(h\cdot\nabla)v$, $(v\cdot\nabla)h$, and $(h\cdot\nabla)h$, belong to $L^2$ by \eqref{eq:6.1_h_bds}.
\cref{lem:taming} further yields $$
\begin{aligned}
\|\TN(h+v)\|_2 &\le C_g(\|h\|_6+C\|v\|_{H^1})^3, \\
\|\TN(h+v)-\TN(h+w)\|_2 &\le CC_g(\|h\|_6^2+\|v\|_{H^1}^2+\|w\|_{H^1}^2) \|v-w\|_{H^1}.
\end{aligned}
$$ Using $$ (v\cdot\nabla)v-(w\cdot\nabla)w =((v-w)\cdot\nabla)v+(w\cdot\nabla)(v-w) $$ and boundedness of $\LP$ on the corresponding spaces, we obtain
\begin{equation*}
\|F_{h(t)}(v)-F_{h(t)}(w)\|_{H^{-1/2}} \le L_{s,T,K,m}\|v-w\|_{H^1}
\end{equation*}
for $\|v\|_{H^1},\|w\|_{H^1}\le m$, uniformly on $[s,T]$.
The same estimates provide a bound for $F_{h(t)}(0)$.

\emph{Step 2:} Choose a Lipschitz $\psi:[0,\infty)\to[0,1]$ equal to one on $[0,1]$ and zero on $[2,\infty)$, and define $$ F_m(t,v):=\psi(\|v\|_{H^1}/m)F_{h(t)}(v).
$$ This map is bounded and globally Lipschitz from $H^1_\sigma$ to $H^{-1/2}_\sigma$.
By \cref{hyp:G} and $H^1\hookrightarrow L^3$, the shifted stochastic coefficients are globally Lipschitz in $H^1$, with quadratic growth bounded by $C_{T,K}(1+\|v\|_{H^1}^2)$.
For an interval $I=[a,a+\ell]\subset[s,T]$, $\ell\le1$, define $$ \|v\|_{\mathscr S^2(I;H^1)} :=\left(\E\sup_{t\in I}\|v(t)\|_{H^1}^2\right)^{1/2}.
$$ The heat-semigroup estimate $$ \|S(r)\|_{\mathcal L(H^{-1/2},H^1)} \le C_\nu r^{-3/4},\qquad 0<r\le1, $$ shows that the drift convolution has Lipschitz constant $C_\nu L_{s,T,K,m}\ell^{1/4}$ in this norm.

Let $[a,b]\subset[s,T]$ be a deterministic interval.
Let $\Phi$ be a predictable $\HS(\U;H^1_\sigma)$-valued process, and let $\Psi$ be a $\mathscr P\otimes\mathcal Z$-measurable $H^1_\sigma$-valued function on $\Omega\times[a,b]\times Z$.
Assume $$ \E\int_a^b \left( \|\Phi(t)\|_{\HS(\U;H^1_\sigma)}^2 +\int_Z\|\Psi(t,z)\|_{H^1}^2\,\mu(\dd z) \right)\dd t<\infty.
$$ Define the stochastic convolution with zero initial datum by $$ \tilde Z(t) :=\int_a^tS(t-r)\Phi(r)\dd W(r) +\int_{(a,t]\times Z}S(t-r)\Psi(r,z)\, \tPN(\dd r,\dd z),\qquad a\le t\le b.
$$ Then
\begin{equation}
\label{eq:6.1_h1_hilbert}
\begin{aligned}
&\E\sup_{a\le t\le b}\|\tilde Z(t)\|_{H^1}^2 +\nu\E\int_a^b \bigl(\|\nabla\tilde Z(t)\|_2^2 +\|\Delta\tilde Z(t)\|_2^2\bigr)\dd t \\
&\qquad\le C\E\int_a^b \left( \|\Phi(t)\|_{\HS(\U;H^1_\sigma)}^2 +\int_Z\|\Psi(t,z)\|_{H^1}^2\,\mu(\dd z) \right)\dd t.
\end{aligned}
\end{equation}
Apply \eqref{eq:5.1_hilbert_eng}, after translating the time origin to $a$, to the stochastic integrands transformed by $\mathcal J:=(I-\Delta)^{1/2}$.
This operator commutes with $S(t)$ and is an isometry $H^1_\sigma\to L^2_\sigma$.
The corresponding dissipation satisfies $$ \|\nabla\mathcal J\tilde Z\|_2^2 =\|\nabla\tilde Z\|_2^2+\|\Delta\tilde Z\|_2^2.
$$ The construction and difference estimate established in \cref{lem:remainder-energy-entrance} also give the continuous Wiener and c\`adl\`ag Poisson versions, using only second moments.
The stochastic convolution map consequently has Lipschitz constant $CL_T^{(1)}\ell^{1/2}$.
Choose $\ell$ so that $$ C_\nu L_{s,T,K,m}\ell^{1/4} +CL_T^{(1)}\ell^{1/2}<1.
$$ The adapted c\`adl\`ag processes with finite $\mathscr S^2$ norm form a complete space modulo indistinguishability.
The drift convolution is continuous by the integrable kernel bound.
The Banach fixed-point theorem, applied successively on finitely many such intervals, therefore yields $V_m\in\mathscr S^2([s,T];H^1)$.
The truncation is applied only to the deterministic drift.

\emph{Step 3:} We establish
\begin{equation}
\label{eq:6.1_h2_1}
\E\int_s^T\|V_m(t)\|_{H^2}^2\dd t<\infty .
\end{equation}
The heat-semigroup evolution of the initial datum has this regularity by Fourier integration and $\xi\in L^2(\Omega;H^1)$; the stochastic convolution does by \eqref{eq:6.1_h1_hilbert} and its maximal estimate.
For $w(t)=\int_s^tS(t-r)f(r)\dd r$, Plancherel's theorem and Young's inequality in time give
\begin{equation}
\label{eq:6.1_sob_bds}
\|w\|_{L^2(s,T;H^{\alpha+2})} \le C_{\nu,T-s}\|f\|_{L^2(s,T;H^\alpha)}, \qquad \alpha=-\tfrac12,\ 0.
\end{equation}
Indeed, $$ \sup_{\zeta\in\R^3}\int_0^{T-s} (1+|\zeta|^2)e^{-\nu|\zeta|^2r}\dd r \le T-s+\nu^{-1}.
$$ First apply this with $\alpha=-1/2$ to the bounded truncated drift.
It follows that $$ V_m\in L^2(\Omega;L^2(s,T;H^{3/2})).
$$ Since $H^{3/2}\hookrightarrow W^{1,3}$, $$ \psi(\|v\|_{H^1}/m)\|(v\cdot\nabla)v\|_2 \le Cm\|\nabla v\|_3\le Cm\|v\|_{H^{3/2}}.
$$ The mixed convection and taming terms satisfy $L^2$ bounds where the cutoff is nonzero.
Consequently,
\begin{equation}
\|F_m(t,V_m(t))\|_2 \le C_{s,T,K,m} \bigl(1+\|V_m(t)\|_{H^{3/2}}\bigr) \quad\text{for almost every }t.
\end{equation}
The second application of \eqref{eq:6.1_sob_bds}, with $\alpha=0$, proves \eqref{eq:6.1_h2_1}.
In particular, $\nu\Delta V_m+F_m(\cdot,V_m)$ belongs to $L^2(\Omega;L^2(s,T;L^2_\sigma))$.
The stochastic Fubini theorem, applied first after spatial spectral truncation and then using these estimates, yields the time-integrated equation in $L^2_\sigma$.

\emph{Step 4:} Construct the truncated solutions consistently on bounded deterministic time intervals and define $$ \sigma_m:=\inf\{t\ge s:\|V_m(t)\|_{H^1}\ge m\}, \qquad \rho:=\lim_m\sigma_m, $$ with $\inf\varnothing=\infty$ and $\sigma_m=s$ when $\|\xi\|_{H^1}\ge m$.
The difference estimates on sufficiently short intervals also apply to deterministic drift terms and stochastic integrands restricted to stopping-time intervals: multiply them by $\mathds1_{\{a<r\le\theta\}}$, form convolutions at deterministic times, and then evaluate at $t\wedge\theta$, as in \cref{lem:stopping-time-restart}.
Before the first exit from an $H^1$ ball, the deterministic drift and quadratic noise integrals  are bounded.
The corresponding second moments of the stopped suprema are finite after initial-data localization, including the exit jump.
Absorption and iteration therefore establish uniqueness on the common closed stochastic interval.

For $n\ge m$, both cutoffs are inactive before $\sigma_m\wedge\sigma_n$, and the stochastic coefficients agree at the terminal time.
Hence $V_m=V_n$ there, including the endpoint, and $\sigma_m\le\sigma_n$.
Modifying the drift at the exit time does not change its Lebesgue time integral.
These solutions define $V$ on $[s,\rho)$; right continuity and the almost sure finiteness of $\xi$ imply $\rho>s$.
If $\rho<\infty$ and the norm before $\rho$ were bounded, then $\sigma_m=\rho$ for every sufficiently large $m$.
Consistency on closed stochastic intervals would make their finite values $V_m(\rho)$ equal, although their norms must be at least $m$, a contradiction.
This proves \eqref{eq:6.1_shifted_h1-1}.
Moreover, $\sigma_m<\rho$ on $\{\rho<\infty\}$ for every $m$.
Indeed, if $\sigma_m=\rho<\infty$, then $V=V_m$ on $[s,\rho)$, and the c\`adl\`ag regularity of $V_m$ on bounded deterministic intervals contradicts the explosion alternative.

Set $k_0:=\lceil\|\xi\|_{H^1}\rceil$ and $$ \theta_n:=\sigma_{k_0+n}\wedge(s+n),\qquad n\ge1.
$$ Since $k_0$ is $\mathcal F_s$-measurable, each $\theta_n$ is a bounded stopping time.
These stopping times satisfy $s<\theta_n\le\theta_{n+1}$, $\theta_n\uparrow\rho$, and $\theta_n<\rho$ on $\{\rho<\infty\}$, almost surely.
The regularity of the truncated solutions and their time-integrated equations therefore hold on these closed intervals.
The explosion alternative also excludes any extension in the stated solution class beyond a finite $\rho$, proving maximality.
The deterministic and Wiener convolutions are continuous in $H^1$; the Poisson jump identity follows by approximation with elementary predictable integrands and \eqref{eq:6.1_h1_hilbert}.
The identity includes the jump at the terminal time and requires no integrability assumption on the first absolute moment of the small-jump coefficient.

\emph{Step 5:} For $A\in\mathcal F_s$, $A\subset\Omega_K$, with $\mathds1_A\xi\in L^2(\Omega;H^1)$, replace $h,\xi$ by $h_A=\mathds1_Ah$, $\xi_A=\mathds1_A\xi$, and use $\Sigma_A=\mathds1_A\Sigma$, $\mathcal G_A=\mathds1_A\mathcal G$ on $(s,\infty)$.
Then $h_A(t)=S(t-s)h_A(s)$, and the preceding estimates apply.
These coefficients are predictable on the shifted interval and retain \cref{hyp:L} and \cref{hyp:G}.
The resulting solution vanishes on $A^c$ and satisfies the original equation on $A$.
Apply this construction on $$ A_n:=\{\|u_0\|_3\le n,\ \|\xi\|_{H^1}\le n\}\in\mathcal F_s.
$$ Uniqueness on stopped stochastic intervals identifies the solutions on their common intervals of existence; the explosion alternative implies equality of their maximal lifetimes there.
With $A_0=\varnothing$, the solutions defined on $A_n\setminus A_{n-1}$ determine the solution for general almost surely finite initial data.
The resulting lifetime is a stopping time because these events belong to $\mathcal F_s$.
The same argument establishes uniqueness after localization on $\mathcal F_s$-measurable events and preserves the stated regularity and jump identities.

\emph{Step 6:} On each bounded closed stochastic interval of localization, the embedding $H^1\hookrightarrow L^3\cap L^5$, boundedness of c\`adl\`ag paths, and \cref{lem:heat-extension} give $$ U=h+V\in D([s,\theta];E)\cap L^5(s,\theta;L^5_\sigma).
$$ Adding $h(t)=S(t-s)h(s)$ to the mild equation yields the critical mild identity with initial time $s$ and datum $\eta=h(s)+\xi$.
The stochastic integrals agree under $H^1_\sigma\hookrightarrow E$ by \cref{hyp:G} and localization; the drift convolutions agree in distributions with the representatives in \cref{lem:3.2_critical-convection} and \cref{lem:critical-taming}.
The stochastic integrals over the corresponding closed stochastic interval include the jump at its terminal time.
For a general measurable $\eta$, the convolution with forcing restricted to $(s,\theta]$ equals $U$ on $[s,\theta]$ and $S(t-\theta)U(\theta)$ for $t>\theta$, by the localized restart identity and the semigroup property.
This gives the c\`adl\`ag convolution representative in \cref{def:critical-solution}.
The solution process $U^\theta$ is frozen at $U(\theta)$ after $\theta$.

For an event $A$ specified before the lemma, $\|\eta_A\|_3\le(K+CM)\mathds1_A$, so the required cubic moment holds.
Defining $U_A:=\mathds1_AU$, $$
\begin{gathered}
U_A\otimes U_A=\mathds1_A(U\otimes U),\qquad \TN(U_A)=\mathds1_A\TN(U), \\
\Sigma_A(t,U_A(t-))=\mathds1_A\Sigma(t,U(t-)),\qquad \mathcal G_A(t,U_A(t-),z)=\mathds1_A\mathcal G(t,U(t-),z).
\end{gathered}
$$ Multiply the critical mild identity with initial time $s$ by $\mathds1_A$ and apply \eqref{eq:3.6_restart_locality} for the future stochastic integrals.
This proves compatibility with the $\mathcal F_s$-localized coefficients, including when the original stochastic coefficients do not vanish at zero.
\end{proof}

The global $H^1$ argument requires a coercive estimate for the drift paired with $(I-\Delta)v$.
The following lemma retains positive coefficients of $\|\Delta v\|_2^2$ and $Q_h(v)$ under the original taming hypothesis \cref{hyp:T}.

\begin{lemma}
[Coercivity of the shifted Sobolev drift] \label{lem:shifted-h1-drift} Assume \cref{hyp:T}.
Let $h\in L^3_\sigma\cap W^{1,\infty}$ satisfy $|h||\nabla h|\in L^2$, and let $v\in H^2_\sigma$.
Set $u=h+v$ and define $$ B_h(v):=\int_{\R^3}|h+v|^2|\nabla h|^2\dd x, \qquad \kappa_{\nu,g}:=\frac6\nu+72\nu C_g^2.
$$ Then all the following pairings and integrals are finite, and
\begin{equation}
\label{eq:6.2_shift_coercivity}
\begin{aligned}
2\bigl(\mathcal A_h(v),-\Delta v\bigr)_{L^2} &+\frac{\nu}{2}\|\Delta v\|_2^2+\frac1{8\nu}Q_h(v)\le c_{N,\nu}\|\nabla v\|_2^2+\kappa_{\nu,g}B_h(v).
\end{aligned}
\end{equation}
Moreover,
\begin{equation}
\label{eq:6.2_shifted_h1}
2\bigl(\mathcal A_h(v),(I-\Delta)v\bigr)_{L^2} +\mathcal D_h(v) \le C_h\bigl(1+\|v\|_{H^1}^2\bigr),
\end{equation}
where one may take
\begin{equation}
\label{eq:6.2_shifted_h1_2}
C_h=C\left( 1+\|h\|_5^5+\|h\|_4^4+\|\nabla h\|_\infty^2 +\bigl\||h||\nabla h|\bigr\|_2^2\right)
\end{equation}
with $C$ depending only on $\nu,C_g,c_{N,\nu}$.

In particular, if $h(t)=S(t)u_0$, $s>0$, and $\|u_0\|_3\le K$, then $C_{h(t)}$ is bounded uniformly for $s\le t\le T$ by a deterministic constant.
\end{lemma}

\begin{proof}
The assumptions imply that $h\in L^q$ for every $3\le q\le\infty$.
Furthermore, $v\in L^\infty\cap L^q$ for $2\le q\le\infty$.
Consequently, $u\in L^6\cap L^\infty$ and $\TN(u)\in L^2$.
The convection term belongs to $L^2$, since $$ (u\cdot\nabla)u=(u\cdot\nabla)v+(u\cdot\nabla)h, \qquad \|(u\cdot\nabla)h\|_2 \le\|\nabla h\|_\infty\|v\|_2 +\bigl\||h||\nabla h|\bigr\|_2 .
$$ Since $u\in W^{1,2}_{\mathrm{loc}}(\R^3;\R^3)\cap L^\infty$, the local Sobolev chain rule gives $$ \partial_j\TN(u) =D\TN(u)(\partial_jv+\partial_jh), \qquad j=1,2,3.
$$ By \cref{lem:taming}, $$ \|D\TN(u)\nabla v\|_2 \le3C_g\|u\|_\infty^2\|\nabla v\|_2.
$$ For the terms involving $\nabla h$, one has $$ \||u|^2|\nabla h|\|_2 \le C\left( \|\nabla h\|_\infty\|v\|_4^2 +\|h\|_\infty\|\nabla h\|_\infty\|v\|_2 +\|h\|_\infty\bigl\||h||\nabla h|\bigr\|_2 \right).
$$ Together with $\TN(u)\in L^2$ and $\|D\TN(u)\|_{\mathrm{op}}\le3C_g|u|^2$, these estimates imply $\TN(u)\in H^1$.
Since $v\in H^2$, the subsequent integration by parts is therefore valid in the weak Sobolev sense.
The orthogonality of the Leray projection in $L^2$ permits its removal in pairings with $v$ and $\Delta v$.

Set $D=\|\Delta v\|_2^2$, $Q=Q_h(v)$, and $B=B_h(v)$.
For the convection term involving $\nabla v$, Cauchy--Schwarz and Young's inequalities give
\begin{equation}
\label{eq:6.2_shifted_conv}
2\bigl|((u\cdot\nabla)v,\Delta v)_{L^2}\bigr| \le2Q^{1/2}D^{1/2} \le\frac{4\nu}{3}D+\frac3{4\nu}Q .
\end{equation}
The convection term involving $\nabla h$ satisfies
\begin{equation}
2\bigl|((u\cdot\nabla)h,\Delta v)_{L^2}\bigr| \le2B^{1/2}D^{1/2} \le\frac{\nu}{6}D+\frac6\nu B .
\end{equation}
For the taming term, integration by parts yields
\begin{equation}
\begin{aligned}
2(\TN(u),\Delta v)_{L^2} =&-2\sum_{j=1}^3\int_{\R^3} \partial_jv\cdot D\TN(u)\partial_jv\dd x -2\sum_{j=1}^3\int_{\R^3} \partial_jv\cdot D\TN(u)\partial_jh\dd x .
\end{aligned}
\end{equation}
By part \textup{(i)} of \cref{lem:taming}, $$ -2\sum_{j=1}^3\int_{\R^3} \partial_jv\cdot D\TN(u)\partial_jv\dd x \le-\frac1\nu Q+c_{N,\nu}\|\nabla v\|_2^2.
$$ For the remaining term, the derivative bound in the same lemma and Cauchy--Schwarz give an upper bound $6C_gQ^{1/2}B^{1/2}$ for its absolute value.
Young's inequality then gives
\begin{equation}
\label{eq:6.2_tame}
6C_g Q^{1/2}B^{1/2} \le\frac1{8\nu}Q+72\nu C_g^2B.
\end{equation}
Combining \eqref{eq:6.2_shifted_conv}-- \eqref{eq:6.2_tame} with $2(\nu\Delta v,-\Delta v)_{L^2}=-2\nu D$, the coefficients of $D$ and $Q$ satisfy $$ 2\nu-\frac{4\nu}{3}-\frac{\nu}{6} =\frac{\nu}{2}, \qquad \frac1\nu-\frac3{4\nu}-\frac1{8\nu} =\frac1{8\nu}.
$$ This proves \eqref{eq:6.2_shift_coercivity}.

For the $L^2$ component of the drift, \cref{lem:mixed-remainder-convection} and \cref{lem:shifted-taming-coercivity}, with $\eta=\nu/2$ in the former, give
\begin{equation}
\label{eq:6.2_L2}
\begin{aligned}
2(\mathcal A_h(v),v)_{L^2} +\nu\|\nabla v\|_2^2 +\frac1{8\nu}\|v\|_4^4 \le \bigl(c_{N,\nu}+c_1\|h\|_5^5\bigr)\|v\|_2^2 +c_2\|h\|_4^4 .
\end{aligned}
\end{equation}
Furthermore,
\begin{equation}
\label{eq:6.2_c_error}
B_h(v) \le 2\|\nabla h\|_\infty^2\|v\|_2^2 +2\bigl\||h||\nabla h|\bigr\|_2^2.
\end{equation}
Adding \eqref{eq:6.2_shift_coercivity} and \eqref{eq:6.2_L2}, and applying \eqref{eq:6.2_c_error}, proves \eqref{eq:6.2_shifted_h1}.

For the final assertion, \cref{lem:heat-extension} gives, on $\Omega_K$ and for $t\ge s>0$,
\begin{equation}
\label{eq:6.2_c_bnds}
\begin{aligned}
\|h(t)\|_5^5&\le C(\nu s)^{-1}K^5,& \|h(t)\|_4^4&\le C(\nu s)^{-1/2}K^4, \\
\|\nabla h(t)\|_\infty^2&\le C(\nu s)^{-2}K^2,& \bigl\||h(t)||\nabla h(t)|\bigr\|_2^2 &\le C(\nu s)^{-3/2}K^4.
\end{aligned}
\end{equation}
Consequently, \eqref{eq:6.2_shifted_h1_2} is uniformly bounded on every deterministic interval $[s,T]$ with $s>0$.
\end{proof}
The preceding coercivity estimate and \cref{hyp:G} yield localized $H^1$ bounds independent of the stopping thresholds, and hence global existence for the shifted equation.
At the deterministic times of \cref{lem:admissible-h1}, critical pathwise uniqueness identifies these solutions with the remainder.
The resulting estimates hold on every finite horizon, with all quantities restricted to times before the critical lifetime, and remain finite when that lifetime is finite.

\begin{proposition}
[Positive-time Sobolev persistence] \label{prop:positive-time-h1-persistence} Assume \cref{hyp:T}, \cref{hyp:L}, and \cref{hyp:G}.
\begin{enumerate}
[label=\textup{(\roman*)},leftmargin=*] \item Every solution constructed in \cref{lem:shifted-local-h1} exists globally on $[s,\infty)$.
For every $T>s$,
\begin{equation}
\label{eq:6.3_global_h1}
V\in D([s,T];H^1_\sigma)\cap L^2(s,T;H^2_\sigma) \quad\text{almost surely}.
\end{equation}
For $K\ge1$, there exists a finite deterministic constant $\Lambda_{s,T,K}$, depending only on $s,T,K,\nu,C_g,c_{N,\nu}$ and $C_T^{(1)}$, such that
\begin{equation}
\label{eq:6.3_h1_est_1}
\begin{aligned}
&\E\left[\mathds1_A\left\{ \sup_{s\le t\le T}\|V(t)\|_{H^1}^2 +\int_s^T\mathcal D_{h(r)}(V(r))\dd r \right\}\right] \le\Lambda_{s,T,K} \left(\E[\mathds1_A\|\xi\|_{H^1}^2] +(T-s)\P(A)\right)
\end{aligned}
\end{equation}
for every $A\in\mathcal F_s$ with $A\subset\Omega_K$ and $\E[\mathds1_A\|\xi\|_{H^1}^2]<\infty$.
The constant is independent of the stopping thresholds.

For $t>s$, set $U_V(t):=h(t)+V(t)$ and define $$
\begin{aligned}
q_W^{(1)}(t):= \|\Sigma(t,U_V(t-))\|_{\HS(\U;H^1_\sigma)}^2,&\qquad q_J^{(1)}(t):= \int_Z\|\mathcal G(t,U_V(t-),z)\|_{H^1}^2\,\mu(\dd z), \\
q^{(1)}(t):&=q_W^{(1)}(t)+q_J^{(1)}(t).
\end{aligned}
$$ After increasing $\Lambda_{s,T,K}$ if necessary,
\begin{equation}
\label{eq:6.3_h1_jump}
\begin{aligned}
&\E\left[\mathds1_A\left\{ \int_s^Tq^{(1)}(r)\dd r +\sum_{s<t\le T}\|\tilde\Delta V(t)\|_{H^1}^2 \right\}\right]\le\Lambda_{s,T,K} \left(\E[\mathds1_A\|\xi\|_{H^1}^2] +(T-s)\P(A)\right).
\end{aligned}
\end{equation}
Apart from the Lipschitz assumptions required for existence, these estimates use only the quadratic $H^1$ growth bound on the stochastic coefficients.
No fourth-moment assumption on the $H^1$-valued Poisson coefficient is required.

\item Let $(u,\{\tau_n\},\tau_*)$ be the maximal critical solution, and write $R=u-h$ on $[0,\tau_*)$.
Use the times $s_j$, the random variables $\xi_{s_j}$, and the events $A_{j,K,M}$ from \cref{lem:admissible-h1}.
For each $j$, let $V^{(j)}$ denote the global solution of \eqref{eq:6.1_shifted} with initial time $s_j$ and initial state $\xi_{s_j}$.
Then, on $\{s_j<\tau_*\}$,
\begin{equation}
\label{eq:6.3_sob_iden}
R(t)=V^{(j)}(t),\qquad s_j\le t<\tau_*,
\end{equation}
simultaneously at all such times, almost surely.

For $T>s_j$, define, on $\{s_j<\tau_*\}$, $$
\begin{aligned}
\mathcal H_{s_j}^*(T):= \sup_{\substack{s_j\le t\le T \\
t<\tau_*}} \|R(t)\|_{H^1}^2 +\int_{s_j}^T\mathds1_{\{r<\tau_*\}} \mathcal D_{h(r)}(R(r))\dd r,
\end{aligned}
$$ and set this functional to zero on $\{s_j\ge\tau_*\}$.
All time integrands are extended by zero outside the interval of existence.
Then
\begin{equation}
\label{eq:6.3_maximal_h1}
\E[\mathds1_{A_{j,K,M}}\mathcal H_{s_j}^*(T)] \le\Lambda_{s_j,T,K}(M^2+T-s_j)\P(A_{j,K,M}).
\end{equation}
Consequently, outside one null set, $R$ has $H^1_\sigma$-valued c\`adl\`ag paths on every compact subinterval of $(0,\tau_*)$ and, for every $0<\delta<T<\infty$, on $\{\delta<\tau_*\}$,
\begin{equation}
\label{eq:6.3_time_max}
\begin{aligned}
\sup_{\substack{\delta\le t\le T \\
t<\tau_*}} \|R(t)\|_{H^1}^2+\int_\delta^T\mathds1_{\{r<\tau_*\}} \left(\|R(r)\|_{H^2}^2 +\int_{\R^3}|u(r)|^2|\nabla R(r)|^2\dd x\right) \dd r<\infty.
\end{aligned}
\end{equation}
In particular, these quantities remain finite when $\tau_*\le T$, with the supremum and integrals taken strictly before $\tau_*$.
\end{enumerate}
\end{proposition}

\begin{proof}
\emph{Step 1:} Fix $T>s$ and an event $A\in\mathcal F_s$ such that $A\subset\Omega_K$ and $\E[\mathds1_A\|\xi\|_{H^1}^2]<\infty$.
The event localization of \cref{lem:shifted-local-h1} allows us to work with $h_A=\mathds1_Ah$, $\xi_A=\mathds1_A\xi$, and the coefficients $\Sigma_A=\mathds1_A\Sigma$ and $\mathcal G_A=\mathds1_A\mathcal G$ on $(s,\infty)$.
The resulting solution vanishes on $A^c$ and agrees with the original solution on $A$ throughout their common interval of existence.
We write the estimates on $A$ using the original notation.

Put $$ \Phi:=\Sigma(\cdot,h+V_-),\qquad \Psi:=\mathcal G(\cdot,h+V_-,\cdot),\qquad X:=\|V\|_{H^1}^2.
$$ For a bounded closed localizing time $\theta\le T$, with $s\le\theta<\rho$ whenever $\rho<\infty$, the energy identity is
\begin{equation*}
\begin{aligned}
X(t\wedge\theta) &=\|\xi\|_{H^1}^2 +2\int_s^{t\wedge\theta} \bigl(\mathcal A_{h(r)}(V(r)), (I-\Delta)V(r)\bigr)_{L^2}\dd r \\
&+\int_s^{t\wedge\theta}q_W^{(1)}(r)\dd r +M_\theta^{(1)}(t)+J_\theta^{(1)}(t), \qquad s\le t\le T,
\end{aligned}
\end{equation*}
where $$
\begin{aligned}
M_\theta^{(1)}(t) :&=2\int_s^t\mathds1_{\{r\le\theta\}} (V(r-),\Phi(r)\dd W(r))_{H^1}+2\int_{(s,t]\times Z}\mathds1_{\{r\le\theta\}} (V(r-),\Psi(r,z))_{H^1}\,\tPN(\dd r,\dd z), \\
J_\theta^{(1)}(t) :&=\int_{(s,t]\times Z}\mathds1_{\{r\le\theta\}} \|\Psi(r,z)\|_{H^1}^2\,\PN(\dd r,\dd z).
\end{aligned}
$$ The stochastic integrands are extended by zero outside $(s,\theta]$; in particular, the jump at $\theta$ is included.
To justify the identity, let $P_n$ be the spatial Fourier projection onto $\{|\zeta|\le n\}$ and put $\mathcal J=(I-\Delta)^{1/2}$.
Apply the $L^2$ It\^o formula to $\mathcal JP_nV$, using the time-integrated $L^2$ equation of \cref{lem:shifted-local-h1}.
For elementary predictable stochastic integrands on finite mark partitions, this is \cite[Theorem 2.1, $p=2$, $M=3$]{MR4165650}; approximation in the quadratic integrability norms gives the general case as in \cref{lem:remainder-energy-entrance}.
Since $$ V\in L^2(s,\theta;H^2),\qquad \mathcal A_h(V)\in L^2(s,\theta;L^2) \quad\text{almost surely}, $$ the drift pairings converge in $L^1(s,\theta)$ to the displayed integrand.
After stopping the relevant integrability quantities, the It\^o isometries and martingale maximal inequalities give convergence of the real-valued martingales; the compensation formula gives $L^1$ convergence of the nonnegative Poisson terms.
The c\`adl\`ag versions and a countable dense set of times give the identity simultaneously on the closed interval.
Increasing the additional stopping times proves the stated local identity.
The case $\theta=s$ is immediate.

\emph{Step 2:} \cref{lem:shifted-h1-drift} gives, on $A$, $$ 2\bigl(\mathcal A_h(V),(I-\Delta)V\bigr)_{L^2} \le-\mathcal D_h(V)+C_{s,K}(1+X).
$$ By \cref{hyp:G}, heat contraction, and $H^1\hookrightarrow L^3$,
\begin{equation}
\label{eq:6.3_h1_growth}
q^{(1)}(r) \le C_T^{(1)}(1+\|h(r)+V(r-)\|_3^2) \le C_{T,K}(1+X(r-)) \quad\text{on }A.
\end{equation}
The processes $X_-$ and $X$ agree almost everywhere with respect to Lebesgue time measure.

Let $\sigma_m$ be the increasing exit times from the truncated-drift construction in \cref{lem:shifted-local-h1}, and set $\theta_m:=T\wedge\sigma_m$.
For $s<r\le\theta_m$, one has $\|V(r-)\|_{H^1}\le m$.
If $\sigma_m=s$, the corresponding stochastic integrals are zero.
Consequently, $$
\begin{aligned}
\E[\mathds1_AJ_{\theta_m}^{(1)}(T)] &\le C_{T,K}(1+m^2)(T-s)\P(A), \\
\E\bigl[\langle\mathds1_AM_{\theta_m}^{(1)}\rangle_T\bigr] &\le4m^2\E\left[\mathds1_A \int_s^{\theta_m}q^{(1)}(r)\dd r\right]<\infty.
\end{aligned}
$$ Here the martingale is viewed with initial time $s$.
Thus $\mathds1_AM_{\theta_m}^{(1)}$ is square integrable.
The energy inequality and these bounds imply finiteness of the expected stopped supremum and dissipation on $A$.
This includes $X(\theta_m)$, which may exceed $m^2$ after a jump.

For $\theta=\theta_m$, define $$ X_\theta^*(t):=\sup_{s\le r\le t\wedge\theta}X(r), \qquad Q_{W,\theta}(t):=\int_s^{t\wedge\theta}q_W^{(1)}(r)\dd r.
$$ Orthogonality of the continuous and purely discontinuous martingale components gives $$ [M_\theta^{(1)}]_t \le4X_\theta^*(t) \bigl(Q_{W,\theta}(t)+J_\theta^{(1)}(t)\bigr).
$$ Since $A\in\mathcal F_s$, the Burkholder--Davis--Gundy and Young inequalities yield $$
\begin{aligned}
\E\left[\mathds1_A\sup_{s\le r\le t} |M_\theta^{(1)}(r)|\right] \le\varepsilon\E[\mathds1_AX_\theta^*(t)]+C_\varepsilon\E\left[\mathds1_A \bigl(Q_{W,\theta}(t)+J_\theta^{(1)}(t)\bigr)\right].
\end{aligned}
$$ Moreover, $$ \E[\mathds1_AJ_\theta^{(1)}(t)] =\E\left[\mathds1_A\int_s^{t\wedge\theta} q_J^{(1)}(r)\dd r\right], $$ because $\mathds1_A\mathds1_{\{s<r\le\theta\}}$ is predictable.
Taking the supremum in the energy inequality, also evaluating it at $t\wedge\theta$, and choosing $\varepsilon$ sufficiently small give $$
\begin{aligned}
F_m(t):=\E\left[\mathds1_A\left\{ X_{\theta_m}^*(t) +\int_s^{t\wedge\theta_m} \mathcal D_{h(r)}(V(r))\dd r\right\}\right] \le C\E[\mathds1_A\|\xi\|_{H^1}^2] +C(t-s)\P(A)+C\int_s^tF_m(r)\dd r.
\end{aligned}
$$ The constant $C=C_{s,T,K}$ is independent of $m$ and $A$.
Gronwall's inequality therefore gives $$ F_m(T)\le Ce^{C(T-s)} \left(\E[\mathds1_A\|\xi\|_{H^1}^2] +(T-s)\P(A)\right).
$$ The squared Poisson amplitudes have been retained in the nondecreasing process $J_\theta^{(1)}$; their fourth moments are not used.

\emph{Step 3:} At each finite exit time, $\|V_m(\sigma_m)\|_{H^1}\ge m$, also when $\sigma_m=s$.
Consistency at the closed exit times gives $$ \P(A\cap\{\sigma_m\le T\}) \le\frac{Ce^{C(T-s)}}{m^2} \left(\E[\mathds1_A\|\xi\|_{H^1}^2] +(T-s)\P(A)\right).
$$ Since $\rho=\lim_m\sigma_m$, $\{\rho\le T\}\subset\{\sigma_m\le T\}$ for every $m$, and letting $m\to\infty$ gives $\P(A\cap\{\rho\le T\})=0$.
Apply this conclusion to $$ \{\|u_0\|_3\le K,\ \|\xi\|_{H^1}\le M\}, \qquad K,M\in\N, $$ and to integer terminal times exceeding $s$.
These events exhaust the probability space up to a null set, so $\rho=\infty$ almost surely.
Monotone convergence along the compatible closed stopping intervals proves \eqref{eq:6.3_h1_est_1}.
Together with $$ \|v\|_{H^2}^2 \le C\bigl(\|v\|_{H^1}^2+\|\Delta v\|_2^2\bigr), $$ this gives \eqref{eq:6.3_global_h1}.

The jump identity of \cref{lem:shifted-local-h1} gives $$ \sum_{s<t\le T}\|\tilde\Delta V(t)\|_{H^1}^2 =\int_{(s,T]\times Z} \|\Psi(t,z)\|_{H^1}^2\,\PN(\dd t,\dd z).
$$ Indeed, partition $Z$ into countably many disjoint sets of finite $\mu$-measure.
On bounded time intervals, the restrictions of the Poisson random measure have finitely many atoms and, almost surely, distinct atom times, both within and between the sets of the partition.
This follows from the diffuse time intensity and independence on disjoint mark sets.
Thus no cross term from simultaneous atoms occurs.
The compensation formula, \eqref{eq:6.3_h1_growth}, and \eqref{eq:6.3_h1_est_1} prove \eqref{eq:6.3_h1_jump}.

\emph{Step 4:} Fix $j,K,M$, and set $s=s_j$ and $A=A_{j,K,M}$.
On $A$, $$ V^{(j)}(s)=R(s),\qquad h(s)+V^{(j)}(s)=u(s),\qquad \|u(s)\|_3\le K+CM.
$$ Fix a deterministic $T>s$.
For each $n$, set $$ B:=A\cap\{s<\tau_n\}\in\mathcal F_s, \qquad \theta:=
\begin{cases}
T\wedge\tau_n,&\text{on }B, \\
T,&\text{on }B^c.
\end{cases}
$$ Then $\theta$ is a bounded stopping time with $\theta>s$.
For every $t\ge s$, define $$ u_B^\theta(t):=
\begin{cases}
u(t\wedge\theta),&\text{on }B, \\
0,&\text{on }B^c,
\end{cases}
\qquad U_B^\theta(t):= \mathds1_B\bigl(h(t\wedge\theta)+V^{(j)}(t\wedge\theta)\bigr).
$$ Both processes are frozen after $\theta$ and vanish on $B^c$.
Their common initial datum is $$ \eta_B:=\mathds1_B(h(s)+\xi_s),\qquad \|\eta_B\|_3\le(K+CM)\mathds1_B, $$ so $\eta_B\in L^3(\Omega,\mathcal F_s;E)$.
Their common stochastic coefficients on $(s,\infty)$ are $\Sigma_B=\mathds1_B\Sigma$ and $\mathcal G_B=\mathds1_B\mathcal G$.

For $u_B^\theta$, the restarted closed mild formulation follows from \cref{lem:stopping-time-restart} and \eqref{eq:3.6_restart_locality}.
For $U_B^\theta$, it follows from the compatibility assertion in \cref{lem:shifted-local-h1}, since part \textup{(i)} gives the required Sobolev regularity on the whole interval $[s,T]$.
The deterministic nonlinearities respect multiplication by $\mathds1_B$, because $\mathds1_B^2=\mathds1_B$ and $\TN(0)=0$.
All stochastic integrals are restricted to $(s,\theta]$, including the terminal jump; no critical integrability of the frozen continuations is required.
\cref{prop:critical-uniqueness-stability}, with initial time $s$, therefore gives $u_B^\theta=U_B^\theta$ simultaneously at all times.

Intersect the probability-one events over $n$ and integer terminal times $T>s$.
Since $\tau_n\uparrow\tau_*$, the resulting intervals exhaust $[s,\tau_*)$ on $A$.
Thus $u=h+V^{(j)}$ there simultaneously in time.
The $L^2$-valued remainder constructed in \cref{Section 5} consequently agrees with $V^{(j)}$ in distributions and hence in $L^2$.
This supplies its $H^1_\sigma$-valued version.
Since the events $A_{j,K,M}$ exhaust $\{s_j<\tau_*\}$, a countable intersection of the probability-one events over $K,M$ proves \eqref{eq:6.3_sob_iden}.

On $A$, the functional $\mathcal H_s^*(T)$ is bounded by the corresponding functional of $V^{(j)}$ on $[s,T]$.
Using $\|\xi_s\|_{H^1}\le M$ on $A$ in \eqref{eq:6.3_h1_est_1} gives \eqref{eq:6.3_maximal_h1}.

\emph{Step 5:} Take a common probability-one event for all indices $j,K,M$ and all integer terminal times.
Fix a sample point in this event and numbers $0<\delta<T$ with $\delta<\tau_*$.
Choose $s_j<\delta$ and then $K,M$ such that the sample point belongs to $A_{j,K,M}$.
On $[\delta,T]\cap[0,\tau_*)$, the remainder is the restriction of the global $H^1_\sigma$-valued c\`adl\`ag process $V^{(j)}$.
Its finite energy on an integer interval containing $[s_j,T]$, the definition of $\mathcal D_h$, and the preceding $H^2$ inequality prove \eqref{eq:6.3_time_max}.
The same argument applies to every such $\delta,T$ on this fixed probability-one event, including when $\tau_*\le T$.

Finally, on $\{s_j\le s_k<\tau_*\}$, the two shifted solutions agree at $s_k$ by \eqref{eq:6.3_sob_iden}.
Uniqueness in \cref{lem:shifted-local-h1}, localized on this $\mathcal F_{s_k}$-measurable event, identifies them for every $t\ge s_k$.
Consequently, at a finite $\tau_*$ they have the same $H^1$ left limit and the same value at $\tau_*$.
No terminal value of the original critical solution is used in this argument.
\end{proof}
The continuation criterion is expressed in the uniform $L^3$ and space-time $L^5$ norms of the full velocity.
The embedding $H^1(\R^3)\hookrightarrow L^3(\R^3)\cap L^5(\R^3)$, together with the heat estimates, converts the preceding bounds for $R$ into estimates for these quantities.
Combined with local critical regularity near zero, these estimates control the continuation quantities on every finite horizon, with the solution evaluated only before its maximal lifetime.

\begin{lemma}
[Transfer to the critical continuation quantities] \label{lem:remainder-critical-transfer} Assume \cref{hyp:T}, \cref{hyp:L}, and \cref{hyp:G}.
Let $u_0\in L^0(\Omega,\mathcal F_0;E)$, let $(u,\{\tau_n\},\tau_*)$ be the maximal critical solution of \cref{thm:main-local}, and write $h=S(\cdot)u_0$ and $R=u-h$ on $[0,\tau_*)$.
At positive times, use the $H^1_\sigma$-valued representative of $R$ supplied by \cref{prop:positive-time-h1-persistence}.

For $0<s<T<\infty$, define on $\{s<\tau_*\}$ $$
\begin{aligned}
B_{s,T} &:=\sup_{\substack{s\le t\le T \\
t<\tau_*}} \|R(t)\|_{H^1},& \mathcal C _{s,T}^* &:=\sup_{\substack{s\le t\le T \\
t<\tau_*}}\|u(t)\|_3 +\left( \int_s^T\mathds1_{\{t<\tau_*\}}\|u(t)\|_5^5\dd t \right)^{1/5},
\end{aligned}
$$ and set both quantities to zero on $\{s\ge\tau_*\}$.
All time integrands are extended by zero outside the interval of existence.

There exists a numerical constant $C$ such that, for every $0<s<T<\infty$, almost surely on $\{s<\tau_*\}$,
\begin{equation}
\label{eq:6.4_time_crit}
\mathcal C _{s,T}^* \le C(1+\nu^{-1/5})\|u_0\|_3 +C\bigl(1+(T-s)^{1/5}\bigr)B_{s,T} <\infty.
\end{equation}
For the deterministic times and events of \cref{lem:admissible-h1}, and every $T>s_j$,
\begin{equation}
\label{eq:6.4_tail_moment}
\E\!\left[ \mathds1_{A_{j,K,M}}(\mathcal C _{s_j,T}^*)^2 \right] \le C_{s_j,T,K}(1+M^2)\P(A_{j,K,M}).
\end{equation}
Here $C_{s_j,T,K}$ may also depend on $\nu,C_g,c_{N,\nu}$ and $C_T^{(1)}$, but is independent of $M$ and of all auxiliary stopping thresholds.

Moreover, outside a single null set, for every $0<T<\infty$,
\begin{equation}
\label{eq:6.4_finite}
\mathcal C _{0,T}^* :=\sup_{\substack{0\le t\le T \\
t<\tau_*}}\|u(t)\|_3 +\left( \int_0^T\mathds1_{\{t<\tau_*\}}\|u(t)\|_5^5\dd t \right)^{1/5} <\infty.
\end{equation}
\end{lemma}

\begin{proof}
For every $v\in H^1(\R^3;\R^3)$, interpolation and the Sobolev embedding into $L^6$ give $$
\begin{aligned}
\|v\|_3 &\le\|v\|_2^{1/2}\|v\|_6^{1/2} \le C\|v\|_{H^1}, \qquad \|v\|_5 \le\|v\|_2^{1/10}\|v\|_6^{9/10} \le C\|v\|_{H^1}.
\end{aligned}
$$ Consequently, on $\{s<\tau_*\}$, $$
\begin{aligned}
\sup_{\substack{s\le t\le T \\
t<\tau_*}}\|R(t)\|_3 &\le CB_{s,T},\quad \left( \int_s^T\mathds1_{\{t<\tau_*\}}\|R(t)\|_5^5\dd t \right)^{1/5} \le C(T-s)^{1/5}B_{s,T}.
\end{aligned}
$$ \cref{prop:positive-time-h1-persistence} ensures that $B_{s,T}<\infty$ almost surely on this event, including when $\tau_*\le T$.
\cref{lem:heat-extension} also gives $$ \sup_{t\ge0}\|h(t)\|_3 +\|h\|_{L^5(Q_\infty)} \le C(1+\nu^{-1/5})\|u_0\|_3.
$$ The triangle inequality in the uniform $L^3$ norm and Minkowski's inequality in the restricted space-time $L^5$ norm, applied to $u=h+R$, prove \eqref{eq:6.4_time_crit}.

For the moment estimate, set $A=A_{j,K,M}$ and $\ell=T-s_j$.
Since $A\subset\Omega_K$, squaring \eqref{eq:6.4_time_crit} yields $$
\begin{aligned}
\E[\mathds1_A(\mathcal C _{s_j,T}^*)^2] \le C(1+\nu^{-2/5})K^2\P(A)+C(1+\ell^{2/5})\E[\mathds1_A B_{s_j,T}^2].
\end{aligned}
$$ By \eqref{eq:6.3_maximal_h1}, $$ \E[\mathds1_A B_{s_j,T}^2] \le\Lambda_{s_j,T,K}(M^2+\ell)\P(A).
$$ Using $M^2+\ell\le(1+\ell)(1+M^2)$ proves \eqref{eq:6.4_tail_moment}, with the stated dependence of the constant.
Only the localized second moment of $B_{s_j,T}$ is required.

Finally, choose a common probability-one event on which $\tau_*>0$, $\tau_n\uparrow\tau_*$, all closed restrictions have finite critical norms, and \eqref{eq:6.4_time_crit} holds for every pair $(s_j,N)$ with $N\in\N$ and $N>s_j$.
Only countably many events are involved.

Fix a sample point in this event and an integer $N\ge1$.
Since $s_j\downarrow0$, choose $j$ such that $s_j<\min\{N,\tau_*\}$, and then choose $n$ such that $s_j<\tau_n$.
The closed critical solution on $[0,\tau_n]$ gives $$ \|u\|_{\mathbb X[0,s_j]}<\infty, $$ while \eqref{eq:6.4_time_crit} gives $\mathcal C _{s_j,N}^*<\infty$.
Splitting the supremum and the time integral at $s_j$, and using $(a+b)^{1/5}\le a^{1/5}+b^{1/5}$ for $a,b\ge0$, we obtain $$ \mathcal C _{0,N}^* \le\|u\|_{\mathbb X[0,s_j]} +\mathcal C _{s_j,N}^* <\infty.
$$ For any finite $T>0$, choose an integer $N\ge T$.
Monotonicity gives $\mathcal C _{0,T}^*\le\mathcal C _{0,N}^*<\infty$ on the same probability-one event, proving \eqref{eq:6.4_finite}.
The argument uses only the countable family of auxiliary solutions constructed at the deterministic times $s_j$.
\end{proof}

\begin{proof}
[Proof of \cref{thm:global-critical-wellposedness}] Assume first that \eqref{eq:initial} holds.
\emph{Step 1.} The maximal solution exists by \cref{thm:maximal-critical-solution}.
Fix $m\in\N$.
On $\{\tau_*\le m\}$, \cref{lem:remainder-critical-transfer} gives, outside a null set, $$ \sup_{0\le t<\tau_*}\|u(t)\|_3<\infty, \qquad \int_0^{\tau_*}\|u(t)\|_5^5\dd t<\infty.
$$ The continuation criterion \eqref{eq:4.4_comb_crit} therefore implies $\P(\tau_*\le m)=0$.
Taking the union over $m$ proves \eqref{eq:6.5_global}.
This argument uses only the solution strictly before its lifetime.

The compatible closed local solutions, with $\tau_n\uparrow\infty$, define a global adapted process satisfying \eqref{eq:6.5_global_2}.
Its mild formulation is, almost surely for all $t\ge0$, $$
\begin{aligned}
u(t)&=S(t)u_0 -\int_0^tS(t-r)\LP\diver(u(r)\otimes u(r))\dd r-\int_0^tS(t-r)\LP\TN(u(r))\dd r \\
&+\int_0^tS(t-r)\Sigma(r,u(r-))\dd W(r)+\int_{(0,t]\times Z}S(t-r)\mathcal G(r,u(r-),z) \,\tPN(\dd r,\dd z).
\end{aligned}
$$ Here $u(0-)=u_0$, and the stochastic integrals are understood by localization.
The deterministic convolutions have continuous $E$-valued representatives by \cref{lem:3.2_critical-convection} and \cref{lem:critical-taming}.
On every finite interval, \cref{hyp:L} and boundedness of the c\`adl\`ag $E$-valued paths imply almost sure finiteness of the Wiener quadratic and the Poisson quadratic and cubic integrability quantities.
After localization of these quantities, \cref{thm:critical-stokes} and \cref{lem:stopping-time-restart} identify the stochastic convolutions with those of the closed local solutions.
Their c\`adl\`ag representatives and a countable intersection over the localizing indices and integer terminal times give the displayed identity simultaneously for all times, including the terminal values of the closed intervals.

\cref{prop:critical-uniqueness-stability}, applied on common closed stochastic intervals, gives equality with any other critical mild solution having the same initial datum, coefficient fields, and driving noises, throughout its interval of existence.
The construction takes place on the prescribed stochastic basis.

\emph{Step 2: } The $L^2_\sigma$-valued remainders of \cref{lem:remainder-energy-entrance} agree in distributions on overlapping closed stochastic intervals.
They therefore determine a global adapted $L^2_\sigma$-valued c\`adl\`ag version of $R$, with $R(0)=0$.
On $\{\tau_*=\infty\}$ we have $\mathcal E_R^*(T)=\mathcal E_R(T)$ for every $T>0$.
Thus \cref{prop:uniform-energy} gives \eqref{eq:6.5_global_6}, the fractional moment of $\mathcal E_R(T)$, and the gradient and $L^4$ integrability in \eqref{eq:6.5_global_3}.
By \eqref{eq:3.1_heat-critical} and H\"older's inequality in time, $$ \|h\|_{L^4(Q_T)}^3 \le T^{3/8}\|h\|_{L^8(0,T;L^4)}^3 \le C\nu^{-3/8}T^{3/8}\|u_0\|_3^3 .
$$ Since $\|R\|_{L^4(Q_T)}^3\le\nu^{3/4}(\mathcal E_R(T))^{3/4}$, the decomposition $u=h+R$ proves \eqref{eq:6.5_global_fractional}.
Furthermore, $$ \int_0^T\|R(t)\|_2^2\dd t \le T\sup_{0\le t\le T}\|R(t)\|_2^2<\infty \quad\text{almost surely}.
$$ This proves the full $L^2(0,T;H^1_\sigma)$ assertion.
Right continuity at zero in $L^3_\sigma$ and $L^2_\sigma$, respectively, gives \eqref{eq:6.5_global_4}.
\cref{prop:positive-time-h1-persistence} gives \eqref{eq:6.5_global_5}.
Its common probability-one event, together with the countable exhaustion by closed solution intervals, gives these pathwise assertions simultaneously for every $0<\delta<T<\infty$.

\emph{Step 3.} Fix $0<\delta<T$, $K\ge1$, and $A\in\mathcal F_0$ with $A\subseteq\Omega_K$.
Put $a=\delta/2$, $X(t)=\|R(t)\|_{H^1}^2$, and $w(t)=\min\{(t-a)/(\delta-a),1\}$ on $[a,T]$.
Step 2 and \cref{prop:positive-time-h1-persistence} identify $R$ with a shifted solution started strictly before $a$.
Thus its localized $H^1$ energy identity holds on $[a,T]$.

Adding \eqref{eq:6.2_shift_coercivity} to the deterministic $L^2$ drift estimate in the proof of \cref{prop:uniform-energy} gives, for almost every $t\in[a,T]$, $$ 2(\mathcal A_h(R),(I-\Delta)R)_{L^2}+\mathcal D_h(R) \le c_{N,\nu}X+b_0\|h\|_4^4+\kappa_{\nu,g}B_h(R), \qquad b_0=18/\nu+3\nu^{1/3}C_g^{4/3}.
$$ Here $B_h(R)=\int |u|^2|\nabla h|^2\dd x$, and hence $$ B_h(R)\le\|u\|_4^2\|\nabla h\|_4^2 \le\tfrac12\bigl(\|u\|_4^4+\|\nabla h\|_4^4\bigr).
$$ On $A$, \cref{hyp:G} yields $q^{(1)}\le C_T(1+K^2+X)$ almost everywhere in time, with $C_T$ independent of $K$.

For $a\le t\le T$, put $$ L(t)=\int_a^t\bigl(q^{(1)}(r)+\mathcal D_{h(r)}(R(r))\bigr)\dd r, \qquad \theta_m=T\wedge\inf\{t\in[a,T]: X(t)>m\ \text{or}\ L(t)\ge m\}, $$ with $\inf\varnothing=\infty$.
These stopping times increase and eventually equal $T$ almost surely by Step 2 and \cref{hyp:G}.
For $a<r\le\theta_m$ one has $X(r-)\le m$, and $L(\theta_m)\le m$ by continuity.

Apply the product rule to $wX$ up to $\theta_m$.
The initial term vanishes and the derivative of the weight contributes $\int w'X\dd t$.
Set $U_m=\sup_{a\le t\le\theta_m} w(t)X(t)$ and $\Psi(r,z)=\mathcal G(r,u(r-),z)$.
For the weighted, stopped real martingale $M_m$, continuity of $w$ and the pre-jump bound give $$ [M_m]_T\le 4U_m\left( \int_a^{\theta_m}wq_W^{(1)}\dd t+ \int_{(a,\theta_m]\times Z}w(t)\|\Psi(t,z)\|_{H^1}^2\, \PN(\dd t,\dd z)\right).
$$ The predictable quadratic variation has expectation at most $4m^2$, and the jump-square integral in this display has expectation at most $m$.
On $A$, the positive drift terms up to $\theta_m$ are integrable by the pre-exit bound and the positive-time heat estimates.
Apply the Burkholder--Davis--Gundy and Young inequalities to $\mathds1_A M_m$, absorb the supremum term, and compensate the jump-square integral.
Using $0\le w\le1$ and $|w'|\le2/\delta$ gives $$
\begin{aligned}
&\E\!\left[\mathds1_A\left(U_m+ \int_a^{\theta_m}w\mathcal D_h(R)\dd t\right)\right] \le C_{\delta,T}\E\!\left[\mathds1_A\int_a^T \bigl(X+1+K^2+\|h\|_4^4+\|u\|_4^4+\|\nabla h\|_4^4\bigr)\dd t\right].
\end{aligned}
$$ The left-hand terms increase to their values on $[a,T]$.
Monotone convergence and $w=1$ on $[\delta,T]$ give the same upper bound for $\E[\mathds1_A\mathcal H_{\delta,T}]$.

We have $\int_a^T X\dd t\le(T+\nu^{-1})\mathcal E_R(T)$ and $\int_a^T\|R\|_4^4\dd t\le\nu\mathcal E_R(T)$.
Also $\|u\|_4^4\le 8(\|R\|_4^4+\|h\|_4^4)$.
On $A$, the heat estimates give $$ \int_a^T\|h\|_4^4\dd t\le C\nu^{-1/2}T^{1/2}K^4,\qquad \|\nabla h(t)\|_4^4\le C(\nu a)^{-5/2}K^4 \quad(a\le t\le T).
$$ It follows from \eqref{eq:6.5_global_6} and $\Gamma_{T,K}\le C_T(1+K^4)$ that \eqref{eq:6.5_global_positive} holds with $C_{\delta,T}$ independent of $A$ and $K$.

Finally let $D_0=\{\|u_0\|_3\le1\}$ and $D_j=\{2^{j-1}<\|u_0\|_3\le2^j\}$ for $j\ge1$.
Applying \eqref{eq:6.5_global_positive} on $D_j$ and using H\"older's inequality gives $$ \E[\mathds1_{D_j}\mathcal H_{\delta,T}^{3/4}] \le \bigl(\E[\mathds1_{D_j}\mathcal H_{\delta,T}]\bigr)^{3/4} \P(D_j)^{1/4} \le C_{\delta,T}(1+2^{3j})\P(D_j).
$$ Since $2^{3j}\le8\|u_0\|_3^3$ on $D_j$ for $j\ge1$, summation proves \eqref{eq:6.5_positive_fractional}.

\emph{Step 4:} \cref{lem:stopping-time-restart} gives the $L^3_\sigma$ jump identity on each closed stochastic interval of localization.
A countable exhaustion gives it simultaneously for all $t>0$.
Since $h$ is continuous in $L^3_\sigma$, $\tilde \Delta u(t)=\tilde \Delta R(t)$ in that space.

For the $H^1_\sigma$ identity, use the deterministic sequence $s_j\downarrow0$ of \cref{lem:admissible-h1}.
On a common probability-one event, the global shifted solutions $V^{(j)}$ satisfy $R=V^{(j)}$ on $[s_j,\infty)$ and obey \eqref{eq:6.1_h1_jump} on every finite interval.
Given $t>0$, choose $j$ with $s_j<t$.
The identification holds also for the left limits at $t$, and continuity of $h$ in $L^3_\sigma$ gives $u(t-)=h(t)+V^{(j)}(t-)$.
Thus $$ \tilde \Delta R(t) =\tilde \Delta V^{(j)}(t) =\int_Z\mathcal G(t,u(t-),z)\,\PN(\{t\},\dd z) \quad\text{in }H^1_\sigma.
$$ The countable choice of $j$ and integer terminal times makes this identity simultaneous for every $t>0$.
The continuous embedding $H^1_\sigma\hookrightarrow L^3_\sigma$ makes it consistent with the $L^3_\sigma$ identity and proves \eqref{eq:6.5_global_7}.

\emph{Step 5.
Extension to measurable initial data.} Let $u_0\in L^0(\Omega,\mathcal F_0;E)$ and use the events $D_k$, the data $a_k$, and the coefficients $(\Sigma_k,\mathcal G_k)$ from the proof of \cref{thm:main-local}.
Multiplication by $\mathds1_{D_k}$ also preserves the compatible $H^1_\sigma$-valued versions in \cref{hyp:G} and all their Lipschitz and growth bounds.
Since $\|a_k\|_3\le k$, Steps 1--4 apply to each of these equations with the same structural and coefficient bounds.
In particular, $\zeta_k=\infty$ almost surely for every $k$.
The construction in \cref{thm:main-local} gives $\tau_*=\zeta_k$ and $u=u^k$ on $D_k$.
It follows that $\tau_*=\infty$ almost surely, with pathwise uniqueness on the prescribed stochastic basis.
Since $\tau_n\uparrow\tau_*=\infty$ almost surely, the closed local mild identities of the \cref{thm:main-local} on $[0,\tau_n]$ yield global mild formulation on $[0,\infty)$ exactly as in step 1. 

Put $h_k=S(\cdot)a_k$ and $R^k=u^k-h_k$.
On $D_k$ one has $h=h_k$ and $R=R^k$ at all times, including the left limits at positive times.
Consequently \eqref{eq:6.5_global_2}--\eqref{eq:6.5_global_5} and the jump identity \eqref{eq:6.5_global_7} transfer from $u^k$.
For the jump term, $\mathcal G_k(t,u^k(t-),z)$ equals $\mathcal G(t,u(t-),z)$ on $D_k$.
Intersecting the probability-one events for the countably many $k$ gives the simultaneous pathwise assertions with the time quantifiers in the statement.

Fix $T>0$, $K\ge1$ and $A\in\mathcal F_0$ with $A\subseteq\{\|u_0\|_3\le K\}$.
Set $A_k=A\cap D_k$.
Then $A_k\subseteq\{\|a_k\|_3\le K\}$ for every $k$.
On $A_k$ the energy agrees with $\mathcal E_{R^k}(T)$ and, for $0<\delta<T$, the positive-time functional agrees with $\mathcal H^k_{\delta,T}$ computed from $(h_k,R^k)$.
Hence $$
\begin{aligned}
\E[\mathds1_A\mathcal E_R(T)] &=\sum_{k\ge1}\E[\mathds1_{A_k}\mathcal E_{R^k}(T)] \le \Gamma_{T,K}\sum_{k\ge1}\P(A_k) =\Gamma_{T,K}\P(A), \\
\E[\mathds1_A\mathcal H_{\delta,T}] &=\sum_{k\ge1}\E[\mathds1_{A_k}\mathcal H^k_{\delta,T}] \le C_{\delta,T}(1+K^4)\P(A).
\end{aligned}
$$ Both estimates use the original threshold $K$ on each $A_k$.
Under \eqref{eq:initial}, the solution agrees by pathwise uniqueness with the one in Steps 1--4, for which \eqref{eq:6.5_global_fractional} and \eqref{eq:6.5_positive_fractional} have already been proved.
\end{proof}

\section{Stability and examples}
\label{Section 7}

Throughout this section, we assume \cref{hyp:T}, \cref{hyp:L}, and \cref{hyp:G}, and compare solutions on a fixed stochastic basis with the same coefficient fields and driving noises.
We first establish continuous dependence on the initial datum for the global critical solutions of \cref{thm:global-critical-wellposedness}.
The convergence is measured by the uniform $L^3$ norm and the space--time $L^5$ norm on finite deterministic intervals.
We then give examples of initial velocities with infinite kinetic energy and of spatially smoothing mixed multiplicative noise.
These constructions make explicit the initial-data and coefficient classes covered by the global theorem.

\begin{proof}
[Proof of \cref{prop:global-critical-stability}] The running supremum in $\mathcal C_{v,w}$ is adapted and c\`adl\`ag, and its integral term is adapted and continuous.
Since $\mathcal C_{v,w}$ is nondecreasing and right-continuous, $\{\sigma_{M,T}(v,w)\le t\}=\{\mathcal C_{v,w}(t)\ge M\}$ for $0\le t<T$.
Thus $\sigma_{M,T}(v,w)$ is a stopping time.
Freeze each global solution at $T$.
The resulting processes, paired with $T$, satisfy \cref{def:critical-solution} by restriction of the mild identities.
Their common stopping time in \cref{prop:critical-uniqueness-stability} is $\sigma_{M,T}(v,w)$ and the stopped critical norm of their difference is $D_{\sigma_{M,T}(v,w)}(v,w)$.
The conditional estimate of that proposition gives \eqref{eq:7.1_closed-1}.

Fix $R,\varepsilon>0$, put $M=2R+1$ and $\sigma=\sigma_{M,T}(v,u)$.
The spatial triangle inequality for the supremum term and Minkowski's inequality for $(v,u)=(u,u)+(v-u,0)$ in the product $L^5$ norm give $$ \mathcal C_{v,u}(t)\le 2\mathcal C_u(t)+D_t(v,u),\qquad 0\le t\le T .
$$ Here the contribution of $(u,u)$ to the integral term is $2^{1/5}\|u\|_{L^5(0,t;L^5)}$.
On $\{\sigma<T\}$, right continuity gives $\mathcal C_{v,u}(\sigma)\ge M$, including when $\sigma=0$ or the crossing occurs at a jump.
Thus $D_\sigma(v,u)\ge1$ on $\{\sigma<T,\ \mathcal C_u(T)\le R\}$.
On $\{\sigma=T\}$, $D_\sigma(v,u)=D_T(v,u)$.
It follows that $$ \{D_T(v,u)>\varepsilon\}\subseteq \{\mathcal C_u(T)>R\}\cup\{D_\sigma(v,u)\ge a_\varepsilon\}.
$$ Conditional Markov's inequality and \eqref{eq:7.1_closed-1} give $$ \P(D_\sigma(v,u)\ge a_\varepsilon\mid\mathcal F_0) \le \min\!\left\{1,\frac{C_{T,M}\|v_0-u_0\|_3^3}{a_\varepsilon^3}\right\}.
$$ Taking expectations proves \eqref{eq:7.1_prob_est}.

Finally, set $d_n=\|u_0^n-u_0\|_3$.
For fixed $T,R,\varepsilon$, let $c=C_{T,\,2R+1}/a_\varepsilon^3$.
For every $\eta>0$, $$ \E[\min\{1,cd_n^3\}]\le\P(d_n>\eta)+c\eta^3 .
$$ Convergence of $d_n$ in probability, followed by $\eta\downarrow0$, makes this expectation tend to zero.
In \eqref{eq:7.1_prob_est} it remains to let $R\uparrow\infty$, since the global critical regularity of $u$ gives $\mathcal C_u(T)<\infty$ almost surely.
This proves \eqref{eq:7.1_global-convergence}.
\end{proof}

The first example shows that the global theory includes solutions whose kinetic energy remains infinite at every finite time.
The rotational field $(-x_2,x_1,0)$ is tangent to spheres, so multiplication by a smooth radial factor preserves the divergence constraint.

\begin{example}
[Explicit infinite-energy initial data] \label{ex:infinite-energy-data} For $1<\alpha\le3/2$, define $$ v_\alpha(x) :=\frac{(-x_2,x_1,0)} {(1+|x|^2)^{(\alpha+1)/2}}, \qquad x\in\R^3.
$$ Then $$ v_\alpha\in C^\infty(\R^3;\R^3)\cap L^3_\sigma, \qquad v_\alpha\notin L^2.
$$ Every deterministic scalar multiple of $v_\alpha$ is therefore an admissible initial datum, without a restriction on its size.

More generally, let $\zeta$ be any almost surely finite real $\mathcal F_0$-measurable random variable, and set $u_0=\zeta v_\alpha$.
Under the standing probabilistic and coefficient measurability assumptions, \cref{hyp:T}, \cref{hyp:L}, and \cref{hyp:G}, \cref{thm:global-critical-wellposedness} applies.
The corresponding solution satisfies, outside a single null set,
\begin{equation}
\label{eq:7.2_finite_eng}
u(\cdot)-u_0\in D([0,T];L^2_\sigma) \quad\text{for every }0<T<\infty, \qquad \lim_{t\downarrow0}\|u(t)-u_0\|_2=0.
\end{equation}
On $\{\zeta\ne0\}$, it also satisfies
\begin{equation}
\label{eq:7.2_inf_eng}
u(t)\notin L^2 \qquad\text{for every finite }t\ge0.
\end{equation}
\end{example}

\begin{proof}
Set $f(x)=(1+|x|^2)^{-(\alpha+1)/2}$.
Since $f$ is radial, $$ \diver v_\alpha =-x_2\partial_1f+x_1\partial_2f=0.
$$ The field is smooth.
For every $1\le q<\infty$, spherical coordinates give $$ \|v_\alpha\|_q^q =A_q\int_0^\infty r^{q+2}(1+r^2)^{-(\alpha+1)q/2}\dd r, \qquad A_q:=\int_{\mathbb S^2} (\omega_1^2+\omega_2^2)^{q/2}\dd\omega, $$ where $0<A_q<\infty$.
The radial integrand behaves as $r^{q+2}$ near zero and as $r^{2-\alpha q}$ at infinity.
Thus $$ v_\alpha\in L^q \quad\Longleftrightarrow\quad \alpha q>3, \qquad 1\le q<\infty.
$$ At $\alpha q=3$, the radial integral diverges logarithmically.
In particular, $v_\alpha\in L^3_\sigma\setminus L^2$ for $1<\alpha\le3/2$.
The map $z\mapsto z v_\alpha$ is continuous into $E$, so $\zeta v_\alpha\in L^0(\Omega,\mathcal F_0;E)$.
If $\E|\zeta|^3<\infty$, then $\E\|\zeta v_\alpha\|_3^3=\|v_\alpha\|_3^3\,\E|\zeta|^3$, and the unconditional fractional estimates of \cref{thm:global-critical-wellposedness} also apply.

For $\rho=|x|^2$ and $F(\rho)=(1+\rho)^{-(\alpha+1)/2}$, direct differentiation gives $$ \Delta\bigl(x_iF(|x|^2)\bigr) =x_i\bigl(4|x|^2F''(|x|^2)+10F'(|x|^2)\bigr).
$$ Consequently, $$ \Delta v_\alpha(x) =(\alpha+1) \frac{\bigl((\alpha-2)|x|^2-5\bigr)(-x_2,x_1,0)} {(1+|x|^2)^{(\alpha+5)/2}}.
$$ This field is smooth, is divergence-free, and decays at least as $|x|^{-\alpha-2}$ at infinity.
Hence $\Delta v_\alpha\in L^2_\sigma\cap L^3_\sigma$.
More generally, $$ |\partial^\gamma v_\alpha(x)| \le C_{\alpha,\gamma}|x|^{-\alpha-|\gamma|}, \qquad |x|\ge1,\quad |\gamma|\le2.
$$ Since $\alpha>1$, these bounds imply $v_\alpha\in W^{2,3}$.

The generator identity for $S(t)=e^{\nu t\Delta}$ on $L^3$ therefore yields $$ S(t)v_\alpha-v_\alpha =\nu\int_0^t S(r)\Delta v_\alpha\dd r.
$$ Since $\Delta v_\alpha\in L^2\cap L^3$, the integral is also well defined in $L^2$, and the two Bochner integrals agree as distributions.
Contractivity and strong continuity of $S(t)$ on $L^2$ give $$ \|S(t)v_\alpha-v_\alpha\|_2 \le\nu t\|\Delta v_\alpha\|_2, \qquad S(\cdot)v_\alpha-v_\alpha\in C([0,T];L^2_\sigma).
$$ In particular, the heat increment vanishes at $t=0$ in $L^2$, although $v_\alpha$ itself does not belong to $L^2$.

Let $R=u-S(\cdot)u_0$.
By \cref{thm:global-critical-wellposedness}, $R\in D([0,T];L^2_\sigma)$ and $R(0)=0$ almost surely.
Thus $$ u(t)-u_0 =\zeta\bigl(S(t)v_\alpha-v_\alpha\bigr)+R(t) $$ defines a c\`adl\`ag $L^2_\sigma$-valued increment vanishing at $t=0$.
A countable intersection over positive integer terminal times proves \eqref{eq:7.2_finite_eng} on a single probability-one event.

Fix a sample point in this event with $\zeta\ne0$.
If $u(t)\in L^2$ at any finite time $t$, then $$ u_0=u(t)-\bigl(u(t)-u_0\bigr)\in L^2, $$ contrary to $u_0=\zeta v_\alpha$ and $v_\alpha\notin L^2$.
This proves \eqref{eq:7.2_inf_eng} simultaneously for every finite time.
\end{proof}
The next example verifies that the coefficient hypotheses allow nontrivial state-dependent forcing for the infinite-energy data just constructed.
Multiplication by a compactly supported function converts an $L^3$ state into an $L^2$ field, convolution then provides $H^1$ regularity, and the Leray projection restores the divergence constraint.
To obtain infinite rank at a suitable state, we use disjoint spatial localizations and Gaussian convolution, whose injectivity preserves linear independence.
Square-summable coefficients ensure the required Hilbert--Schmidt bound.

The L\'evy density is chosen independently near zero and at infinity: the first part gives infinite activity with finite quadratic intensity, whereas the second gives finite second and third moments and an infinite fourth moment.
This separates the spatial regularity of the coefficients from the integrability of their jump amplitudes.

\begin{example}
[Smoothing mixed multiplicative noise] \label{ex:smoothing-mixed-noise} Let $\U=\ell^2(\N)$ with canonical orthonormal basis $(e_k)_{k\ge1}$.
Choose $a\in C_c^\infty(B(0,2);\R)$ such that $0\le a\le1$ and $a=1$ on $B(0,1)$.
For $k=0,1,\ldots$, define $$ x_k=(6k,0,0),\qquad a_k(x)=a(x-x_k),\qquad \kappa(x)=(4\pi)^{-3/2}e^{-|x|^2/4}, $$ and $$ Q_k(v):=\LP\bigl[\kappa*(a_kv)\bigr], \qquad v\in E.
$$ Fix $1\le\vartheta<2$ and $3<\beta<4$.
On $Z=\R\setminus\{0\}$, equipped with its Borel $\sigma$-field, let $$ \mu(\dd z) =\mathds1_{\{0<|z|\le1\}}|z|^{-1-\vartheta}\dd z +\mathds1_{\{|z|>1\}}|z|^{-1-\beta}\dd z.
$$ Define the deterministic coefficients $$ \Sigma(v)\eta :=\sum_{k=1}^\infty 2^{-k}(\eta,e_k)_\U Q_k(v), \qquad \mathcal G(v,z):=zQ_0(v).
$$

These coefficients satisfy \cref{hyp:L} and \cref{hyp:G}, including their measurability and compatibility requirements.
There exists $v\in E$ such that $\Sigma(v)$ has infinite rank and $Q_0(v)\ne0$.
The measure $\mu$ is a L\'evy measure with $\mu(Z)=\infty$.
For every $v\in E$ with $Q_0(v)\ne0$,
\begin{equation}
\label{eq:7.3_fourth_moment}
\int_Z\|\mathcal G(v,z)\|_{H^1}^4\,\mu(\dd z)=\infty,
\end{equation}
whereas all second- and third-order jump integrability conditions in \cref{hyp:L} and \cref{hyp:G} hold.

For the fields $v_\alpha$ of \cref{ex:infinite-energy-data}, one also has $Q_0(v_\alpha)\ne0$.
The coefficients are compatible with a smooth taming function satisfying \cref{hyp:T} and with the initial data $u_0=\zeta v_\alpha$ considered in that example.
\end{example}

\begin{proof}
\emph{Step 1:} The Gaussian kernel satisfies $$ \|\kappa\|_1=1, \qquad \|\partial_j\kappa\|_1=\pi^{-1/2}, \qquad j=1,2,3.
$$ Set $$ A:=\sqrt{1+3/\pi}\,\|a\|_6.
$$ For $v,w\in E$, H\"older's inequality gives $$ \|a_k(v-w)\|_2 \le\|a\|_6\|v-w\|_3.
$$ The Leray projection is contractive on $L^2$ and commutes with spatial derivatives, hence is contractive on $H^1$.
Young's convolution inequality therefore yields

$$\|Q_k(v)-Q_k(w)\|_{H^1}^2 \le \left(\|\kappa\|_1^2+ \sum_{j=1}^3\|\partial_j\kappa\|_1^2\right) \|a_k(v-w)\|_2^2 \le A^2\|v-w\|_3^2.$$ In particular, $Q_k:E\to H^1_\sigma$ is a bounded linear map, and
\begin{equation}
\label{eq:7.3_lipschitz}
\|Q_k(v)-Q_k(w)\|_{H^1} \le A\|v-w\|_3, \qquad Q_k(0)=0, \qquad k\ge0.
\end{equation}
The bound is uniform in $k$ and requires only the $L^3$ norm of $v-w$.

Since $\sum_{k\ge1}4^{-k}=1/3$, the partial sums defining $\Sigma(v)$ converge in $\HS(\U;H^1_\sigma)$, uniformly for $v$ in bounded subsets of $E$.
Furthermore, $$\|\Sigma(v)-\Sigma(w)\|_{\HS(\U;H^1_\sigma)}^2 =\sum_{k=1}^\infty 4^{-k}\|Q_k(v)-Q_k(w)\|_{H^1}^2 \le \frac{A^2}{3}\|v-w\|_3^2.$$

Let $\phi:H^1_\sigma\hookrightarrow E$ denote the continuous embedding, and set $\tilde C:=\|\phi\|$.
The Gaussian definition of the $\gamma$-radonifying norm and the Hilbert-space Gaussian isometry give $$ \|\phi\Phi\|_{\gamma(\U;E)} \le \tilde C\|\Phi\|_{\HS(\U;H^1_\sigma)}, $$ first for finite-rank $\Phi:\U\to H^1_\sigma$ and then for every $\Phi\in\HS(\U;H^1_\sigma)$ by completion.
Identifying $\Sigma(v)$ with its composition with $\phi$ when regarded as an $E$-valued operator, we obtain
\begin{equation*}
\|\Sigma(v)-\Sigma(w)\|_{\gamma(\U;E)} \le\frac{\tilde CA}{\sqrt3}\|v-w\|_3.
\end{equation*}

\emph{Step 2:} The measure $\mu$ is $\sigma$-finite, has infinite total mass, and satisfies $$ \int_Z(1\wedge|z|^2)\,\mu(\dd z) =\frac2{2-\vartheta}+\frac2\beta<\infty.
$$ Thus it is a L\'evy measure.
For $\vartheta<p<\beta$, $$ m_p:=\int_Z|z|^p\,\mu(\dd z) =2\int_0^1r^{p-1-\vartheta}\dd r +2\int_1^\infty r^{p-1-\beta}\dd r =\frac2{p-\vartheta}+\frac2{\beta-p}.
$$ In particular, $m_2,m_3<\infty$, whereas $$ \int_Z|z|^4\,\mu(\dd z)=\infty, \qquad \int_{\{0<|z|\le1\}}|z|\,\mu(\dd z)=\infty.
$$ The latter divergence does not obstruct the compensated Poisson integral, which is defined by quadratic integrability, with localization when necessary.

By \eqref{eq:7.3_lipschitz}, $$ \|\mathcal G(v,\cdot)-\mathcal G(w,\cdot)\|_{L^2(Z;H^1_\sigma)} \le A m_2^{1/2}\|v-w\|_3.
$$ For $p=2,3$, the embedding $H^1_\sigma\hookrightarrow E$ also gives
\begin{equation*}
\|\mathcal G(v,\cdot)-\mathcal G(w,\cdot)\|_{L^p(Z;E)} \le \tilde CA m_p^{1/p}\|v-w\|_3.
\end{equation*}
 The constants in \cref{hyp:L} and \cref{hyp:G} may be chosen as $$
\begin{aligned}
L_T &=\tilde CA \bigl(3^{-1/2}+m_2^{1/2}+m_3^{1/3}\bigr), &K_T&=0, \\
L_T^{(1)} &=A\bigl(3^{-1/2}+m_2^{1/2}\bigr), &K_T^{(1)}&=0.
\end{aligned}
$$ These choices are independent of $T$.
In particular, $$ \|\Sigma(v)\|_{\HS(\U;H^1_\sigma)}^2 +\int_Z\|\mathcal G(v,z)\|_{H^1}^2\,\mu(\dd z) \le A^2(1/3+m_2)\|v\|_3^2.
$$ Similarly, $$
\begin{aligned}
&\|\Sigma(v)\|_{\gamma(\U;E)}^3 +\|\mathcal G(v,\cdot)\|_{L^2(Z;E)}^3 +\|\mathcal G(v,\cdot)\|_{L^3(Z;E)}^3 \le (\tilde CA)^3 \bigl(3^{-3/2}+m_2^{3/2}+m_3\bigr)\|v\|_3^3.
\end{aligned}
$$ All required growth estimates therefore hold.
For $Q_0(v)\ne0$, one has $$ \int_Z\|\mathcal G(v,z)\|_{H^1}^4\,\mu(\dd z) =\|Q_0(v)\|_{H^1}^4 \int_Z|z|^4\,\mu(\dd z) =\infty, $$ which proves \eqref{eq:7.3_fourth_moment}.

\emph{Step 3:} The preceding bounds show that $\Sigma$ is continuous from $E$ into $\HS(\U;H^1_\sigma)$ and into $\gamma(\U;E)$.
The map $v\mapsto\mathcal G(v,\cdot)$ is continuous into $L^2(Z;H^1_\sigma)$ and into $L^2(Z;E)\cap L^3(Z;E)$.
Moreover, $(v,z)\mapsto zQ_0(v)$ is jointly continuous from $E\times Z$ into $H^1_\sigma$.
These spaces are separable.
Since the coefficients are deterministic and independent of time, they have the predictable measurability required by \cref{hyp:L} and \cref{hyp:G}.
The $H^1_\sigma$- and $E$-valued definitions agree under $\phi$, and all estimates hold for every $v,w\in E$.

Choose a nonzero field $w\in C_c^\infty(B(0,1);\R^3)$ with $\diver w=0$; for example, $$ w=(\partial_2\varphi,-\partial_1\varphi,0) $$ for a suitable $\varphi\in C_c^\infty(B(0,1))$.
Set $$ w_k=w(\cdot-x_k), \qquad v_*:=\sum_{k\ge1}2^{-k}w_k.
$$ The supports of the $w_k$ are pairwise disjoint, so $$ \|v_*\|_3^3 =\sum_{k\ge1}2^{-3k}\|w\|_3^3<\infty.
$$ The partial sums are solenoidal, hence $v_*\in E$.
The support conditions on $a_k$ imply $a_kv_*=2^{-k}w_k$ for $k\ge1$.
Since Gaussian convolution preserves solenoidality, $$ Q_k(v_*)=2^{-k}\kappa*w_k, \qquad \Sigma(v_*)e_k=4^{-k}\kappa*w_k.
$$

Gaussian convolution is injective on $L^2$, because its Fourier multiplier is everywhere nonzero.
A finite linear relation among the fields $\kappa*w_k$ therefore implies the same relation among the $w_k$.
Their disjoint supports and nonzero norms imply that all coefficients vanish.
Thus $\Sigma(v_*)$ has infinite rank, both as an $H^1_\sigma$-valued operator and as an $E$-valued operator.

Let, $\tilde v:=w+v_*$  and it satisfies $$ \Sigma(\tilde v)=\Sigma(v_*), \qquad Q_0(\tilde v)=\kappa*w\ne0.
$$ Indeed, $a_kw=0$ for $k\ge1$, $a_0v_*=0$, and $a_0w=w$.
This proves the simultaneous nontriviality asserted in the statement.

For completeness, consider the infinite-energy field $v_\alpha$ of \cref{ex:infinite-energy-data}.
Since $a_0=1$ near the origin, $$ \bigl(\nabla\times(a_0v_\alpha)\bigr)(0) =\bigl(\nabla\times v_\alpha\bigr)(0) =(0,0,2).
$$ Thus $\nabla\times(a_0v_\alpha)$ is a nonzero compactly supported smooth field.
The Leray projection preserves distributional curl and commutes with Gaussian convolution, so $$ \nabla\times Q_0(v_\alpha) =\kappa*\bigl(\nabla\times(a_0v_\alpha)\bigr).
$$ Injectivity of Gaussian convolution on $L^2$ shows that the right-hand side is nonzero.
Hence $Q_0(v_\alpha)\ne0$.
By linearity, \eqref{eq:7.3_fourth_moment} also holds at $u_0=\zeta v_\alpha$ on $\{\zeta\ne0\}$.

The condition $\mu(Z)=\infty$ concerns the intensity of the Poisson random measure.
It does not assert nonzero jumps along every solution; in particular, $\Sigma(0)=0$ and $\mathcal G(0,z)=0$.

\emph{Step 4:} For the fixed taming parameter $N>0$, choose $\chi\in C^\infty(\R)$ such that $0\le\chi\le1$, $\chi=0$ on $(-\infty,0]$, and $\chi=1$ on $[1,\infty)$.
Define $$ g_N(r):=\nu^{-1}\int_0^r\chi(s-N)\dd s, \qquad r\ge0.
$$ Then $g_N\in C^\infty([0,\infty))$ and $$ g_N(0)=0,\qquad 0\le g_N'\le\nu^{-1},\qquad g_N(r)\ge(r-N-1)_+/\nu.
$$ Consequently, $$ \frac r\nu-2g_N(r) \le
\begin{cases}
(N+1)/\nu,&0\le r\le N+1, \\
(2N+2-r)/\nu,&r\ge N+1.
\end{cases}
$$ Since the expression on the left vanishes at $r=0$, $$ 0\le c_{N,\nu} =\sup_{r\ge0}\left(\frac r\nu-2g_N(r)\right) \le\frac{N+1}{\nu}.
$$ Thus \cref{hyp:T} holds with $C_g=\nu^{-1}$.

On any stochastic basis satisfying the standing assumptions and carrying the prescribed independent Wiener process and Poisson random measure, these coefficients and this taming function satisfy the hypotheses of \cref{thm:global-critical-wellposedness} and \cref{prop:global-critical-stability}.
In particular, both results apply to the initial data $u_0=\zeta v_\alpha$ of \cref{ex:infinite-energy-data}.
\end{proof}

\begin{remark}
[Nonlinear coefficients with predictable amplitudes] The preceding construction extends to nonlinear spatial maps and random predictable amplitudes.
Let $e_1,\ldots,e_m$ be an orthonormal family in $\U$ and set $I=\{1,\ldots,m,\mathrm G\}$.
For $j\in I$, take real-valued functions $a_j\in L^2\cap L^6$ and $\kappa_j\in W^{1,1}$, and a globally Lipschitz map $F_j:\R^3\to\R^3$ with Lipschitz constant $\ell_j$.
Define $$ Q_j(v):=\LP\bigl[\kappa_j*(a_jF_j(v))\bigr], \qquad v\in E, $$ where $F_j(v)(x)=F_j(v(x))$.
H\"older's and Young's inequalities, together with the boundedness of $\LP$ on $H^1$, give $$
\begin{aligned}
\|Q_j(v)-Q_j(w)\|_{H^1} &\le C\|\kappa_j\|_{W^{1,1}}\ell_j\|a_j\|_6\|v-w\|_3, \\
\|Q_j(0)\|_{H^1} &\le C\|\kappa_j\|_{W^{1,1}}|F_j(0)|\|a_j\|_2.
\end{aligned}
$$ Thus $Q_j:E\to H^1_\sigma$ is globally Lipschitz. 

Let $b_1,\ldots,b_m$ be real predictable processes and let $\lambda:\Omega\times[0,\infty)\times Z\to\R$ be $\mathscr P\otimes\mathcal Z$-measurable.
Assume that, for each $T<\infty$, there is a finite deterministic constant $B_T$ such that $$ \max_{1\le k\le m}|b_k(t)| +\|\lambda(t,\cdot)\|_{L^2(Z,\mu)} +\|\lambda(t,\cdot)\|_{L^3(Z,\mu)} \le B_T $$ for $\P\otimes\mathrm dt$-almost every $(\omega,t)\in\Omega\times[0,T]$.
Then the coefficients $$ \Sigma(t,v)\eta :=\sum_{k=1}^m b_k(t)(\eta,e_k)_\U Q_k(v), \qquad \mathcal G(t,v,z):=\lambda(t,z)Q_{\mathrm G}(v) $$ satisfy \cref{hyp:L} and \cref{hyp:G}.
Indeed, the spatial estimates and the amplitude bounds give the required Hilbert--Schmidt and Bochner-space estimates; the continuous embedding $H^1_\sigma\hookrightarrow E$ and the ideal property of $\gamma$-radonifying operators give the corresponding critical-space bounds.
Joint measurability, Bochner Fubini, and continuity of the maps $Q_j$ provide the required strongly predictable versions, which are compatible under this embedding.
\end{remark}

The local well-posedness theorem \cref{thm:main-local}, the global continuation of \cref{thm:global-critical-wellposedness}, and \cref{prop:global-critical-stability} establish the three principal statements.
The examples verify that their hypotheses admit arbitrary-size critical initial data with infinite energy and nontrivial mixed multiplicative noise.

\section*{CRediT authorship contribution statement}
Bikram Podder : Conceptualization, Formal analysis, Investigation, Methodology, Writing- original draft.

Surendra Kumar : Supervision, Validation, Writing- review and editing.
\section*{ Declaration of competing interest }
The authors declare that they have no known competing financial interest or personal relationship that could have appeared to influence the work reported in this paper.
\section*{ Data availability}
No data has been used to prepare this manuscript.
\bibliographystyle{elsarticle-num} \bibliography{Reference}
\end{document}